\documentclass[paper=a4,abstracton,numbers=noenddot]{scrartcl}

\usepackage[english]{babel}  
\usepackage{hyperref}
\usepackage[intlimits]{amsmath}
\usepackage{amsthm}
\usepackage{amssymb} 
\usepackage{units} 
\usepackage{amsfonts,mathrsfs,bbm}
\allowdisplaybreaks[1]
\usepackage[numbers]{natbib}
\usepackage{graphicx}
\usepackage[space]{grffile} 
\usepackage{tikz}
\usepackage[T1]{fontenc}
\usepackage{mathptmx}
\usepackage[scaled=.92]{helvet}
\usepackage{courier}

\usepackage[a4paper,left=2.5cm,right=2.5cm,top=3cm,bottom=2.5cm]{geometry}
\usepackage{authblk}
\usepackage{amsmath,bm}
\usepackage{graphicx}
\usepackage{mathrsfs}
\usepackage{amssymb}
\usepackage{physics}
\usepackage{tensor}
\usepackage{xcolor}

\title{Nonlinear static isogeometric analysis of arbitrarily curved Kirchhoff-Love shells }

\author[1,2]{G. Radenkovi\'{c}}
\author[2,3]{A. Borkovi\'{c}}
\author[3]{B. Marussig}

\affil[1]{Faculty of Civil Engineering, University of Belgrade, Bulevar kralja Aleksandra 73, 11000 Belgrade, Serbia}
\affil[2]{University of Banja Luka, Faculty of Architecture, Civil Engineering and Geodesy, Department of Mechanics and Theory of Structures, 78000 Banja Luka, Bosnia and Herzegovina, aleksandar.borkovic@aggf.unibl.org, aborkovic@tugraz.at}
\affil[3]{Institute of Applied Mechanics, Graz University of Technology, Technikerstra\ss e 4/II, 8010 Graz, Austria}
\date{}                     
\begin{document}

	\newcommand{\ssub}[2]{{#1}_{#2}} 
	\newcommand{\vsub}[2]{\textbf{#1}_{#2}} 
	\newcommand{\ssup}[2]{{#1}^{#2}} 
	\newcommand{\vsup}[2]{\textbf{#1}^{#2}} 
	\newcommand{\ssupsub}[3]{{#1}^{#2}_{#3}} 
	\newcommand{\vsupsub}[3]{\textbf{#1}^{#2}_{#3}} 

	\newcommand{\veq}[1]{\bar{\textbf{#1}}} 
	\newcommand{\seq}[1]{\bar{#1}} 
	\newcommand{\ve}[1]{\textbf{#1}} 
	\newcommand{\sdef}[1]{#1^*} 
	\newcommand{\vdef}[1]{{\textbf{#1}}^*} 
	\newcommand{\vdefeq}[1]{{\bar{\textbf{#1}}}^*} 
	\newcommand{\trans}[1]{\textbf{#1}^\mathsf{T}} 
	\newcommand{\transmd}[1]{\dot{\textbf{#1}}^\mathsf{T}} 
	\newcommand{\mdvdef}[1]{\dot{\textbf{#1}}^*} 
	\newcommand{\mdsdef}[1]{\dot{#1}^*} 
	\newcommand{\mdv}[1]{\dot{\textbf{#1}}} 
	\newcommand{\mds}[1]{\dot{#1}} 
	
	\newcommand{\loc}[1]{\hat{#1}} 
	\newcommand{\md}[1]{\dot{#1}} 

	\newcommand{\ii}[3]{{#1}^{#2}_{#3}} 
	\newcommand{\iv}[3]{\textbf{#1}^{#2}_{#3}} 
	\newcommand{\idef}[3]{{#1}^{* #2}_{#3}} 
	\newcommand{\ivdef}[3]{\textbf{#1}^{* #2}_{#3}} 
	\newcommand{\iloc}[3]{\hat{#1}^{#2}_{#3}} 
	\newcommand{\ieq}[3]{\bar{#1}^{#2}_{#3}} 
	\newcommand{\iveq}[3]{\bar{\textbf{#1}}^{#2}_{#3}} 
	\newcommand{\ieqdef}[3]{\bar{#1}^{* #2}_{#3}} 
	\newcommand{\iveqdef}[3]{\bar{\textbf{#1}}^{* #2}_{#3}} 

	\newcommand{\ieqmddef}[3]{\dot{\bar{#1}}^{* #2}_{#3}} 
	\newcommand{\iveqmddef}[3]{\dot{\bar{\textbf{#1}}}^{* #2}_{#3}} 
	
	\newcommand{\ieqmd}[3]{\dot{\bar{#1}}^{#2}_{#3}} 
	\newcommand{\iveqmd}[3]{\dot{\bar{\textbf{#1}}}^{#2}_{#3}} 
	
	\newcommand{\imddef}[3]{\dot{#1}^{* #2}_{#3}} 
	\newcommand{\ivmddef}[3]{\dot{\textbf{#1}}^{* #2}_{#3}} 
	
	\newcommand{\imd}[3]{\dot{#1}^{#2}_{#3}} 
	\newcommand{\ivmd}[3]{\dot{\textbf{#1}}^{#2}_{#3}} 
	
	\newcommand{\iii}[5]{^{#2}_{#3}{#1}^{#4}_{#5}} 
	\newcommand{\iiv}[5]{^{#2}_{#3}{\textbf{#1}}^{#4}_{#5}} 
	\newcommand{\iivn}[5]{^{#2}_{#3}{\tilde{\textbf{#1}}}^{#4}_{#5}} 
	\newcommand{\iiieq}[5]{^{#2}_{#3}{\bar{#1}}^{#4}_{#5}} 
	
	\newcommand{\eqqref}[1]{Eq.~\eqref{#1}} 
	\newcommand{\fref}[1]{Fig.~\ref{#1}} 

\maketitle

\section*{Abstract}
The geometrically rigorous nonlinear analysis of elastic shells is considered in the context of finite, but small, strain theory. The research is focused on the introduction of the full shell metric and examination of its influence on the nonlinear structural response. \textcolor{black}{The exact relation between the reference and equidistant strains is employed and the complete analytic elastic constitutive relation between energetically conjugated forces and strains is derived.} Utilizing these strict relations, the geometric stiffness matrix is derived explicitly by the variation of the unknown metric. Moreover, a compact form of this matrix is presented. Despite the linear displacement distribution due to the Kirchhoff-Love hypothesis, a nonlinear strain distribution arises along the shell thickness. This fact is sometimes disregarded for the nonlinear analysis of thin shells based on the initial geometry, thereby ignoring the strong curviness of a shell at some subsequent configuration. We show that the curviness of a shell at each configuration determines the appropriate shell formulation. For shells that become strongly curved at some configurations during deformation, the nonlinear distribution of strain throughout the thickness must be considered in order to obtain accurate results. We investigate four computational models: one based on the full analytical constitutive relation, and three simplified ones. \textcolor{black}{Robustness, accuracy and relative efficiency of the presented formulation are examined via selected numerical experiments.} Our main finding is that the employment of the full metric is often required when the complete response of the shells is sought, even for the initially thin shells. Finally, the simplified model that provided the best balance between efficiency and accuracy is suggested for the nonlinear analysis of strongly curved shells.

\vspace{5mm}
	
\noindent
\textbf{Keywords}: strongly curved Kirchhoff-Love shell; nonlinear analysis; reciprocal shift tensor; analytical constitutive relation
	
\section{Introduction}
\label{intro}
	
Thin shells are essential components of various modern-day engineering structures and their analysis is a classical topic of applied mechanics, \cite{1973naghdi}. Due to permanent advances in computational capabilities and increasing requirements from industry, the development of more accurate mathematical and mechanical models is an everlasting topic. The nonlinear analysis of thin shells is often performed assuming the linear elastic stress-strain relation which is appropriate when these structures exhibit finite displacements and rotations but small strains. The analysis of arbitrarily curved elastic shells is often simplified by neglecting higher-order terms in the metric and the strain-displacement relation. 

The seminal work by Naghdi, \cite{1973naghdi}, is recommended as a primary source for shell theory, while a comprehensive overview of contemporary achievements in shell theories and finite element applications can be found in \cite{2004bischoff}. Besides the classical finite element method, \cite{2015jeon, 2017ko}, several other efficient techniques exist, such as the meshfree approaches, \cite{2020sadamoto, 2020dai}, the generalized quadrature method, \cite{2015kurtaran}, and the finite strip method, \cite{2017borkovic}. Also, there is ongoing work on the enhancements of shell theories, see e.g.~\cite{2019wu}. In the following, we will focus on the contributions within the scope of isogeometric analysis (IGA) \cite{2005hughes}, since it is the discretization technique employed herein. IGA utilizes CAD basis functions (e.g.~NURBS) for finite element analysis. They are employed for the description of the geometry as well as the physical fields. Many useful computational properties stem from this approach, the most prominent one is increased interelement continuity, which provided new impetus for the development of thin shell formulations. 

The first IGA Kirchhoff-Love (KL) shell element was developed in \cite{2009kiendl} where the required continuity between KL elements is easily handled within the scope of IGA. Since the second-order derivatives of displacements are present in the variational equations of the KL theory, the approximative functions must have at least \textit{C}\textsuperscript{1} interelement continuity. Using NURBS basis functions enables \textit{C}\textsuperscript{$p$-1} interelement continuity where $p$ is the order of the basis. This alleviates locking issues which often occur due to low order discretizations. Furthermore, higher continuity allows the implementation of the KL shell models without rotational DOFs. A large deformation rotation-free formulation was presented in \cite{2011benson} where different choices for defining the shell normal vector are discussed. The formulation for compressible and incompressible hyperelastic thin shells that can use general 3D constitutive models is developed in \cite{2015kiendl}. Numerical tests, including structural dynamics simulations of a bioprosthetic heart valve, emphasize the good performance and applicability of the present method. An extended IGA for the analysis of through-the-thickness cracks is developed in \cite{2015nguyen-thanh} while the application of the subdivision algorithm is discussed in \cite{2016riffnaller-schiefer}. In all these papers, a linear strain distribution is assumed.
	
Nonlinear buckling of thin shells is considered in \cite{2016luo} where a reduced-order model based on Koiter perturbation is utilized as a predictor, leading to an improved efficiency. An isogeometric KL shell formulation is extended to the large strain elasto-plastic analysis using a stress-based approach in \cite{2018ambati}. Furthermore, a reduced model with a simplified strain measure is considered in \cite{2019leonetti} and no significant loss of accuracy is detected, even for very large displacements and composite structures. The authors in \cite{2017duong} analyze shell models extracted from 3D continua and those directly formulated in a 2D manifold form in which the 2D Koiter model gives excellent performance. \textcolor{black}{Notably, a new type of instability named \textit{snap-backward} is detected in \cite{2020nguyen} while examining the postbuckling behavior of imperfect shells.}
	
\textcolor{black}{Although single-patch NURBS surfaces are of primary interest here, let us note that the construction of more complex geometries is an on-going topic in IGA. For example, the assembly of structures consisting of multiple shell patches requires special attention. The bending strip method is introduced in \cite{2010kiendl} while the coupling without penalty or stabilization parameters is suggested in \cite{2017coox}. Furthermore, other CAD basis functions, besides NURBS, may also be used. For instance, the application of subdivision surfaces \cite{2011cirak} or T-splines \cite{2017casquero} was discussed. The latter shows that various trimmed geometries can be converted to T-splines and used in the analysis. Trimmed geometries are omnipresent in most engineering CAD models, but they also impose several computational challenges \cite{2018marussig}. In the context of isogeometric KL analysis, the adequate integration of trimmed surfaces and the application of boundary conditions using penalty approaches \cite{2012schmidt}, \cite{2015breitenberger} or Nitsche’s method \cite{2015guo,2017guo,2018guo} were of particular interest.}
	
While this research deals exclusively with KL shells, we draw attention to the fact that the first IGA-based thick shell formulation within the scope of Mindlin-Reissner theory is presented in \cite{2010benson} while implicit dynamics are considered in \cite{2017sobota}. Since the rotation field is independent of the displacement field, Lagrange-based finite elements with \textit{C}\textsuperscript{0} interelement continuity are well-suited for these shells, \cite{2004lee}. Nevertheless, IGA provides important enhancements for the modeling of Mindlin-Reissner shells, such as the exact initial geometry and increased smoothness of both rotational and displacement fields, \cite{2013dornisch,2017kiendl}. Interesting work is given in \cite{2017oesterle}, where two hierarchic finite element formulations for geometrically nonlinear shell analysis including the effects of transverse shear are presented. Both methods combine a fully nonlinear KL shell model with hierarchically added linearized transverse shear components. The underlying assumption is that in most practical applications the transverse shear angles are small, which is confirmed by various numerical experiments. The hierarchic construction results in an additive strain decomposition into two parts: one resulting from the membrane and bending deformation and the other from transverse shear. The authors conclude that large rotations of Mindlin-Reissner shells can be modeled without parametrizing the corresponding rotation tensor.
	
It is often reckoned that KL theory is appropriate for thin shells with  $ h/R \leq 0.05 $, where $ R $ is the minimum shell radius of curvature and $ h $ is the thickness, \cite{2018li}. This limit is not strict, however, and it is a subject to be further investigated, \cite{1973naghdi}. Unlike strongly curved beams, \cite{2015cazzani,2018borkovic,2019borkovic,2018radenkovic,2020radenkovic}, strongly curved shells are rarely considered in the literature. The influences of the approximation of shift tensor, nonlinear strain distribution and integration through the thickness of thick shells are discussed in \cite{2004bischoff}. The classification of thin shell models with emphasis on strongly curved shells is considered in \cite{2003hamdouni} using the asymptotic expansion technique. The same method is utilized in \cite{1988gamby} for the study of thick KL shells. This relative lack of literature is attributed to the fact that the KL theory is appropriate for the thin shells. The ratio $h/R$ is already mentioned as the one type of a measure of shell slenderness. We suggest a more comprehensive measure of slenderness, the product $Kh$, where $K$ is the absolute value of the trace of the curvature tensor. We will refer to this quantity as the curviness of a shell, analogously to beams, \cite{2018radenkovic}. Therefore, the construction \textit{arbitrarily curved} refers to the curviness of a shell and not merely to its curvature. Since we are dealing with nonlinear analysis, this parameter will be considered during the whole deformation process, \cite{2004bischoff}. As strict limits do not exist, we will refer to the shells with curviness $Kh \leq 0.05$ as small-curvature shells, while all the others belong to the class of strongly curved, i.e.~large-curvature shells. It is emphasized that this measure has a strictly local character. 
	
\textcolor{black}{The present research is based on \cite{2017radenkovic} but with significant theoretical improvements and thorough numerical investigation. In order to study arbitrarily curved KL shells, the complete metric of the shell is considered and exact relations between strain at \textit{the equidistant surface} and reference strains are employed. Here, \textit{the equidistant surface} refers to the surface which consists of points that are equally distanced from the shell midsurface, measured along the midsurface normal.} Analytical integration of the virtual power along the thickness is complex due to the presence of higher order rational functions in the constitutive relation, as noted in \cite{2017duong,2004bischoff}, although the mathematics is well-posed, \cite{1973naghdi}. This integration is the main theoretical contribution of the present research. We propose to utilize the analytical expression for the reciprocal shift tensor. The full metric can be considered numerically, as well. Analytical integration, however, allows us to define reduced constitutive models and compare them with the complete one. The simplification is not applied in the strain-displacement relation for these reduced models, but in the constitution itself. We show that the simplest constitutive plate model that decouples bending and membrane actions gives reasonably accurate results for small-curvature shells, whereas more involved coupling is required for strongly curved shells. \textcolor{black}{By utilizing this novel integration and equations of motion derived in \cite{2017radenkovic}, the geometric stiffness matrix is obtained explicitly by the variation of the current, unknown metric.} This includes the variation of Christoffel symbols, the determinant of the metric tensor and the normal. Considering the complexity of these expressions, elegant and compact form of the geometric stiffness matrix is obtained. Finally, the validity of the present formulation is confirmed via several numerical examples. By considering the efficiency analysis and accuracy of results, the most involved reduced model is recommended as the favorable one.
	
The paper is structured as follows. The next section deals with the metric of a shell with the emphasis on the reciprocal shift tensor. Subsequently, the strict definition of the shell’s kinematics is introduced and the rigorous relation between equidistant and reference strains is highlighted. The virtual power of an arbitrarily curved shell is postulated and discretized in the fourth section. It is shown that through-the-thickness integration can be performed analytically. Thorough numerical experiments are given in the fifth section and these are followed by the conclusions.

\section{Metric of the shell continuum}
	
A detailed and rigorous definition of the shell metric is presented in this section. The classical KL assumption states that a straight line element is rigid and remains perpendicular to the shell’s midsurface in the deformed configuration. This assumption leads to the degeneration of a 3D continuum model into an arbitrarily shaped surface. Note that another approach also exists in which the surface is primarily observed and the so-called directors are introduced at each point on the surface, \cite{1973naghdi,2004bischoff}. The material law that follows from this approach, however, neglects the influence of curvature of the shell, restricting its applicability solely to small-curvature shells, \cite{2004bischoff}. The present analysis is performed with respect to the convective frame of reference while the complete shell kinematics is defined by the translation of shell midsurface.
	
\textbf{Remark}. It is interesting to note that the convective coordinate frame is not an arbitrary choice among the spatial coordinate frames. We must find coordinates and base vectors in order to fully describe the metric of the spatial configuration. This is a significant issue since we are essentially dealing with one equation with two unknowns. The solution is to keep one quantity constant, either the coordinates or the base vectors. It is a standard practice in computational mechanics to keep the coordinates constant, which gives rise to the convective description, where every particle at each configuration has the same coordinates.
	
Regarding the notation, boldface lowercase letters are used for vectors, while uppercase is used for tensors and matrices. An overbar designates the quantity at the equidistant surface of the shell while the asterisk sign is used to designate deformed configuration. Finally, the hat symbol designates the local component of a vector, with respect to the curvilinear coordinates. The direct and index notations are applied simultaneously, depending on the context, and the standard summation convention is adopted. The Greek index letters take values of 1 and 2 while the Latin indices take values of 1, 2, and 3. Partial and covariant derivatives with respect to the \textit{m}\textsuperscript{th} coordinate are designated with $( )_{,m}$ and $( )_{\vert m}$, respectively.
	
The elaboration on the NURBS-based IGA surface modeling is excluded for brevity since it is readily available in the literature. For a detailed insight into the IGA, references \cite{2005hughes} and \cite{2009kiendl} are recommended. 
	
\subsection{Metric of the midsurface}
	
\begin{figure}
	\includegraphics[width=\linewidth]{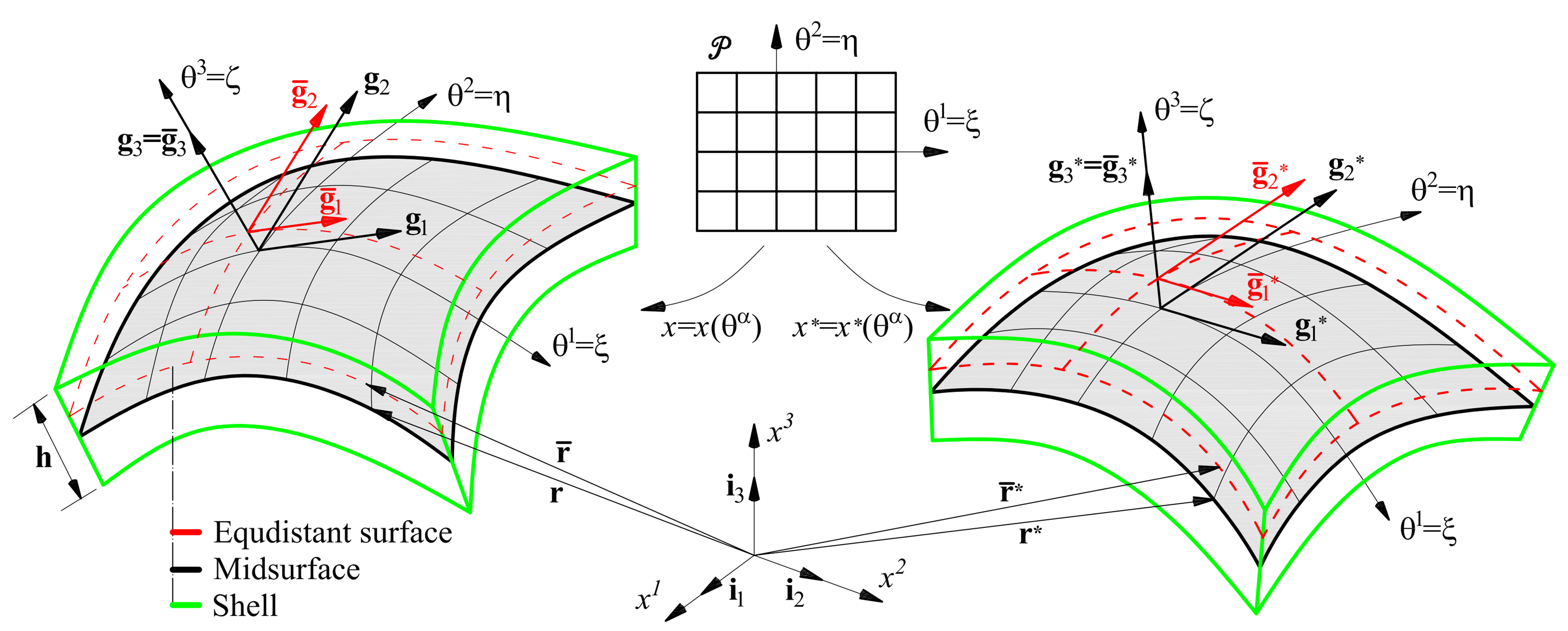}
	\caption{The mapping between (center) the parameter domain $\mathscr{P} $, (left) reference and (right) deformed configurations of a Kirchhoff-Love shell. The shell boundaries are colored with solid green lines, the equidistant surface is indicated by dashed red lines, and the midsurface is colored in gray with black boundaries. Quantities at the deformed configuration are denoted with an asterisk while those at the equidistant surfaces have an overbar.}
	\label{fig:shell}
\end{figure}

The midsurface of a shell with respect to the Cartesian coordinates is determined by its position vector $\textbf{r}=\{x=x^1, y=x^2, z=x^3\}$ which will be described as: 
	\begin{equation}
	\label{eq:def:r}
	\textbf{r} = x^m \textbf{i}_m=\sum\limits_{I =1}^{N} \sum\limits_{J =1}^{M} R_{IJ} (\xi,\eta) \textbf{r}_{IJ},  \quad   
	x^m = \sum\limits_{I = 1}^{N} \sum\limits_{J = 1}^{M} R_{IJ} (\theta^\alpha) x^m_{IJ},
	\end{equation}
\noindent where $R_{IJ}$ are bivariate NURBS basis functions, $\textbf{r}_{IJ} = \{x_{IJ}, y_{IJ}, z_{IJ}\}$ are the position vectors of the control points $IJ$, while $N$ and $M$ are the total numbers of control points along the $\xi$ and $\eta$ directions, respectively, \cite{2005hughes}. Furthermore, $\textbf{i}^m = \textbf{i}_m$ are the base vectors of the Cartesian coordinate system while $\theta^\alpha$ are coordinates $\xi$ and $\eta$ such that $\theta^1=\xi$ and $\theta^2=\eta$, see \fref{fig:shell}. 
	
Tangent vectors of the midsurface are given by 
	\begin{equation}
	\label{eq:def:tangvec}
	\textbf{g}_\alpha =\frac{\partial \textbf{r}}{\partial \theta^\alpha}= x^m_{,\alpha} \textbf{i}_m,
	\quad
	x^m_{,\alpha} = \sum\limits_{I=1}^{N} \sum\limits_{J =1}^{M} R_{IJ,\alpha} x^m_{IJ},
	\end{equation}
\noindent where $x^m_{,\alpha}$ are the partial derivatives of the Cartesian components of the midsurface with respect to the $\theta^\alpha$ coordinates. In order to describe the complete geometry of a shell, the third coordinate $ \theta^3 = \zeta $ is required. Due to the KL hypothesis, this coordinate line is rigid, straight and perpendicular to the $\theta^\alpha$ coordinate lines. It follows that the base vector $\textbf{g}_3$ and its reciprocal counterpart $\textbf{g}^3$ are the same. They have the unit length and represent the normal at each point of the shell midsurface: 
\begin{equation}
	\label{eq:def:normal}
	\textbf{g}^3 = \textbf{g}_3 = \frac{1}{\sqrt{g}} (\textbf{g}_1 \times \textbf{g}_2)= x^m_{,3}\textbf{i}_m,
	\quad
	x^m_{,3} = x_{m,3} = \frac{1}{\sqrt{g}} x^k_{,1} x^l_{,2} e_{klm},
\end{equation}
\noindent where $e_{klm}$ is permutation symbol while $g$ is the determinant of the metric tensor of shell midsurface: 
\begin{equation}
	\label{eq:def:met_tensor}
	g_{ij} = 
	\begin{bmatrix}
		x^m_{,1} x_{m,1} & x^m_{,1} x_{m,2} & 0 \\
		x^m_{,2} x_{m,1} & x^m_{,2} x_{m,2} & 0 \\
		0 & 0 & 1
	\end{bmatrix}
	=
	\begin{bmatrix}
		g_{11} & g_{12} & 0 \\
		g_{21} & g_{22} & 0 \\
		0 & 0 & 1
	\end{bmatrix},
	\; 
	g_{12} = g_{21} , \; g = \det g_{ij} = g_{11} g_{22}-{g_{12}}^2.
\end{equation}
\noindent The reciprocal metric tensor of midsurface is:	
\begin{equation}
	\label{eq:def:rec_met_tensor}
	g^{ij} = {g_{ij}}^{-1} = \frac{1}{g}
	\begin{bmatrix}
		g_{22} & -g_{12} & 0 \\
		-g_{12} & g_{11} & 0 \\
		0 & 0 & g
	\end{bmatrix}
	=
	\begin{bmatrix}
		g^{11} & g^{12} & 0 \\
		g^{21} & g^{22} & 0 \\
		0 & 0 & 1
	\end{bmatrix},
	\; 
	\det g^{ij} = 1/g,
\end{equation}
\noindent while the reciprocal tangent vectors are:
\begin{equation}
	\label{eq:def:rec_tang_vec}
	\begin{aligned}
		\textbf{g}^{1} &= \frac{1}{\sqrt{g}} (\textbf{g}_{2} \times \textbf{g}_{3}) = x^{,1}_n \textbf{i}^n,
		\quad
		x^{,1}_n = \frac{1}{\sqrt{g}} x^{m}_{,2} x^{k}_{,3} e_{mkn}, \\
		\textbf{g}^{2} &= \frac{1}{\sqrt{g}} (\textbf{g}_{3} \times \textbf{g}_{1}) = x^{,2}_n \textbf{i}^n,
		\quad
		x^{,2}_n = \frac{1}{\sqrt{g}} x^{m}_{,3} x^{k}_{,1} e_{mkn}.
	\end{aligned}
\end{equation}
\noindent In order to completely define the metric of the midsurface, the Christoffel symbols are required. They stem from the differentiation of the base vectors with respect to the curvilinear coordinates, \cite{1973naghdi}:
\begin{equation}
	\label{eq:def:christ}
	\textbf{g}_{\alpha,\beta} = x^m_{,\alpha\beta} \textbf{i}_m = \Gamma_{\alpha \beta k} \textbf{g}^k = \Gamma^{n}_{\alpha \beta} \textbf{g}_n 
	\Rightarrow 
	\begin{cases}
	\begin{aligned}
	\Gamma_{\alpha \beta k} &= x^m_{,\alpha\beta} x_{m,k}\\
	\Gamma^n_{\alpha \beta} &= x^m_{,\alpha\beta} x^{,n}_m 
	\end{aligned}
	\end{cases} 
	\Rightarrow \Gamma^n_{\alpha \beta} = g^{nk} \Gamma_{\alpha \beta k}.
\end{equation}
\noindent Here, $\Gamma_{\alpha \beta k}$ and $\Gamma^n_{\alpha \beta}$ are the Christoffel symbols of the first and the second kind, respectively. The components of the differentiated base vectors with respect to the Cartesian coordinates are:
\begin{equation}
	\label{eq:def:diff_base_vec}
	x^m_{,\alpha \beta} = \Gamma^n_{\alpha \beta} x^m_{,n} = \Gamma_{\alpha \beta n}x^{m,n}.
\end{equation}

The Christoffel symbols of the second kind related to the third coordinate line $\Gamma^3_{\alpha \beta}$ are of special importance since they represent curvatures of the surface with respect to the coordinate system $\theta^\alpha$. They determine the curvature tensor and it is a common practice to designate a particular symbol for them, i.e., $\Gamma^3_{\alpha \beta}=b_{\alpha\beta}$. The components of the curvature tensor of a surface can be determined by differentiating the normal with respect to the curvilinear coordinates:
\begin{equation}
	\label{eq:def:diff_base_v3}
	\textbf{g}_{3,\alpha} = - \Gamma^\mu_{3 \alpha} \textbf{g}_\mu = - b^\mu_\alpha \textbf{g}_\mu , \; ( \Gamma^3_{3 \alpha} = 0 ),
\end{equation}
\noindent where $b^\mu_\alpha$ are mixed components of the curvature tensor which relate to the covariant ones by $b^\mu_\alpha = g^{\mu \nu} b_{\nu \alpha}$, \cite{1973naghdi}.

\subsection{Metric of an equidistant surface}

\textcolor{black}{The position vector of an equidistant surface and its appropriate tangent vectors, \fref{fig:shell}, are:} 
\begin{equation}
	\label{eq:def:r_eq}
	\begin{aligned}
	\veq{r} &= \ve{r} + \zeta \ve{g}_3, \\
	\veq{g}_\alpha &= \frac{\partial \veq{r}}{\partial \theta^\alpha} = \ve{g}_\alpha - \zeta b^\nu_\alpha \iv{g}{}{\nu} = ( \delta^\nu_\alpha - \zeta b^\nu_\alpha ) \ve{g}_\nu = \seq{C}^\nu_\alpha \ve{g}_\nu , \; 
	\seq{C}^\nu_\alpha = \delta^\nu_\alpha - \zeta b^\nu_\alpha,
	\end{aligned}
\end{equation}
\noindent where $\seq{C}^\nu_\alpha$ are components of the shift tensor, \cite{1990naomis}, while $\delta^\nu_\alpha$ are components of the Kronecker tensor. This expression reveals the fact that the normal $\ve{g}_3 = \veq{g}_3$ and the tangent vectors of an equidistant surface are mutually orthogonal: 
\begin{equation}
	\label{eq:def:ortho_tan_eq}
	\ve{g}_3 \cdotp \veq{g}_{\alpha} = \seq{C}^\nu_\alpha g_{3 \nu} = 0.
\end{equation}
\noindent Furthermore, the metric tensor and its determinant at an equidistant surface are: 
\begin{equation}
	\label{eq:def:metric_ten_eq}
	\begin{aligned}
	\seq{g}_{\alpha \beta} &= \veq{g}_\alpha \cdotp \veq{g}_\beta = \seq{C}^\mu_\alpha \seq{C}^\nu_\beta g_{\mu \nu}, \\
	\seq{g} &= \det \veq{g}_{\alpha \beta} = \left( \det \seq{C}^\mu_\alpha \right) \left( \det \seq{C}^\nu_\beta \right) g.
	\end{aligned}
\end{equation}
\noindent Using the relations
\begin{equation}
	\label{eq:def:det_shift}
	\det \seq{C}^\mu_\alpha  = \det ( \delta^\mu_\alpha - \zeta b^\mu_\alpha) = 1 - \zeta b^\nu_\nu + \zeta^2 b, \quad b=\det b^\mu_\alpha, \quad b^\nu_\nu = \Tr b^\mu_\alpha, 
\end{equation}
\noindent the determinant of the metric tensor at an equidistant surface reduces to:
\begin{equation}
	\label{eq:def:det_metric_eq}
	\seq{g} = \left( 1- \zeta b^\nu_\nu + \zeta^2 b \right)^2 g = {g_0}^2 g, \quad g_0 = \det \seq{C}^\mu_\alpha,
\end{equation}
\noindent where $ g_0 $ is the determinant of the shift tensor, also known as \textit{the shifter}, \cite{2017duong}, while $ b^\nu_\nu $ and $ b $ are the trace and determinant of the curvature tensor, respectively. Note that the trace of the curvature tensor is equal to the sum of the principal curvatures, $ b^\nu_\nu = K_1 + K_2 = 2H $, where $ H $ is \textit{the mean curvature}, \cite{1990naomis}. 

For the explicit derivation of the constitutive relation which follows, it is crucial to establish the analytical relation between the tangent vectors of the equidistant and middle surfaces. From \eqqref{eq:def:r_eq} it follows that the reciprocal shift tensor $ C^\alpha_\nu = \left( \seq{C}^\nu_\alpha \right) ^{-1} $ could be introduced as well:
\begin{equation}
	\label{eq:def:basis_rec_eq}
	\veq{g}^\alpha = C^\alpha_\nu \ve{g}^\nu,
\end{equation}
\noindent and it can be easily numerically calculated. However, we are interested in its analytical form. Let us utilize the fact that the base vector $ \ve{g}_3 $ is independent of $ \zeta $ coordinate and observe the vector product:
\begin{equation}
	\label{eq:def:vec_prod_eq}
	\veq{g}_3 \times \veq{g}_\mu = \ve{g}_3 \times \seq{C}^\gamma_\mu \ve{g}_\gamma = \sqrt{g} \seq{C}^\gamma_\mu e_{3 \gamma \nu} \ve{g}^\nu.
\end{equation}
\noindent If we multiply this expression with $ e^{3 \mu \alpha} $:
\begin{equation}
	\label{eq:hlp:shift}
	\sqrt{\seq{g}} \delta^\alpha_\beta \veq{g}^\beta = \sqrt{g} \seq{C}^\gamma_\mu \delta^{\mu \alpha}_{\gamma \nu} \ve{g}^\nu,
\end{equation}
\noindent the analytical form of the reciprocal shift tensor follows as:
\begin{equation}
	\label{eq:def:rec_shift_tens}
	C^\alpha_\nu = \frac{1}{g_0} \delta^{\alpha \mu}_{\nu \gamma} \seq{C}^\gamma_\mu = \frac{1}{g_0} \left[ \left( 1 - \zeta b^\gamma_\gamma \right) \delta^\alpha_\nu + \zeta b^\alpha_\nu \right].
\end{equation}
\noindent The same result can be found in \cite{1973naghdi} and it allows us to represent the reciprocal metric tensor at an arbitrary point as an analytic function of the metric and curvature of the midsurface: 
\begin{equation}
	\label{eq:rel:met_tens_eq_via_ref}
	\seq{g}^{\alpha \beta} = \veq{g}^\alpha \veq{g}^\beta = C^\alpha_\mu C^\beta_\nu g^{\mu \nu}.
\end{equation}
\noindent The tensor $ C^\alpha_\mu C^\beta_\nu $ is a rational function of the $ \zeta $ coordinate where the numerator and denominator are polynomials of the $2^{nd}$ and $4^{th}$ degree, respectively. The relation between the metrics of the equidistant and reference surfaces directly influences the equivalent strain relations, \cite{2004bischoff}. In order to simplify the computational model, many types of approximations are readily introduced by researchers, mostly based on first- or second-order Taylor approximations of the shift tensor, \cite{1968baker}. The exact analytical form of the reciprocal shift tensor, on the other hand, is often disregarded. An interesting discussion on the influence of the exact shift tensor is given in \cite{2004bischoff}. 

Finally, let us define the relations between a differential area at the equidistant and middle surfaces which follow from \eqref{eq:def:det_metric_eq}. If we notice that these areas are:
\begin{equation}
	\label{eq:def:diff area}
	\begin{aligned}
	\dd{\veq{a}_\alpha} &= \dd{\seq{a}_\alpha \veq{g}^\alpha} = \dd{\zeta} \ve{g}_3 \times \veq{g}_\beta = \dd{\zeta} \sqrt{\seq{g}} e_{3 \beta \alpha} \veq{g}^\alpha, \quad \dd{\seq{a}_\alpha} = \sqrt{\seq{g}} \dd{\zeta} e_{3 \beta \alpha}, \\
	\dd{\ve{a}_\alpha} &= \dd{a_\alpha \ve{g}^\alpha} = \dd{\zeta} \ve{g}_3 \times \ve{g}_\beta = \dd{\zeta} \sqrt{g} e_{3 \beta \alpha} \ve{g}^\alpha, \quad \dd{a_\alpha} = \sqrt{g} \dd{\zeta} e_{3 \beta \alpha},
	\end{aligned}
\end{equation}
\noindent the required relation follows as:
\begin{equation}
	\label{eq:def:diff area eq}
	\dd{\seq{a}_\alpha} = \sqrt{\frac{\seq{g}}{g}} \dd{a_\alpha} = g_0 \dd{a}_\alpha.
\end{equation}
\noindent Additionally, the differential volume of a shell is: 
\begin{equation}
	\label{eq:def:diff volume}
	\dd{v} = \left( \dd{\xi} \ve{g}_1 \times \dd{\eta} \ve{g}_2 \right) \cdotp \dd{\zeta} \ve{g}_3 = \sqrt{\seq{g}} \dd{\xi} \dd{\eta} \dd{\zeta} = g_0 \dd{a} \dd{\zeta},
\end{equation}
\noindent where $\dd{a} = \sqrt{g} \dd{\xi} \dd{\eta} $ is the differential area of the midsurface.

\section{Kirchhoff-Love shell theory}

After the metric of the shell continuum is defined, the next step is to introduce some strain measure. For the convective coordinate frame, the Lagrange strain equals the difference between the current and reference metrics:
\begin{equation}
	\label{eq:def:strain}
	\ii{\epsilon}{}{\alpha\beta} = \frac{1}{2} \left( \idef{g}{}{\alpha \beta} - \ii{g}{}{\alpha \beta} \right).
\end{equation}
\noindent Due to the KL hypothesis, the components of the shear strain rate along the normal direction vanish, i.e.~$ d_{\alpha 3}=0 $. Consequently, the position of the line element is completely determined by the translation of the midsurface, which is the only generalized coordinate of the KL shell. This fact gives rise to the so-called \textit{rotation-free} shell theories, \cite{2009kiendl}. These formulations are especially efficient for the finite rotation analysis since the parametrization of the finite rotation variables is not required.

\subsection{Kinematics of the shell}

As indicated in \fref{fig:shell}, the position vector of an equidistant surface in the deformed configuration is:
\begin{equation}
	\label{eq:def:r equidistant def}
	\ieqdef{\ve{r}}{}{} = \idef{\ve{r}}{}{} (\zeta) = \idef{\ve{r}}{}{} + \zeta \idef{\ve{g}}{}{3},
\end{equation}
\noindent where $ \vdef{r} $ is the position vector of the deformed midsurface while the base vectors $ \vdef{g}_m $ are calculated analogously to the Eqs.~\eqref{eq:def:tangvec} and \eqref{eq:def:normal}. The position vector of the deformed midsurface is:
\begin{equation}
	\label{eq:def:r def2}
	\vdef{r} = \ve{r} + \ve{u},
\end{equation}
\noindent where the displacement $ \ve{u} $ is discretized with NURBS, in the same way as the geometry, \eqqref{eq:def:r}:
\begin{equation}
	\label{eq:def:u}
	\ve{u} = \sum\limits_{I =1}^{N} \sum\limits_{J =1}^{M} R_{IJ} (\theta^\alpha) \ve{u}_{IJ}.
\end{equation}
\noindent $ \ve{u}_{IJ} $ is the vector that consists of displacement components of a control point $IJ$ with respect to the Cartesian system. For the sake of the seamless transition to the discrete equation of motion which follows in Section 4, let us introduce the matrix of basis functions $ \ve{N} $ such that \eqqref{eq:def:u} can be written as:
\begin{equation}
	\label{eq:def:u via matrices}
	\ve{u} = \ve{N} \ve{q},
\end{equation}
\noindent where:
\begin{equation}
	\label{eq:def:matrices for u}
	\begin{aligned}
		\trans{q} &=
		\begin{bmatrix}
			\trans{u}_{11} & \trans{u}_{12} & ... & \trans{u}_{IJ} & ... & \trans{u}_{NM}
		\end{bmatrix}, \quad
		\ve{u}_{IJ} = 
		\begin{bmatrix}
			u^1_{IJ} & u^2_{IJ} & u^3_{IJ}
		\end{bmatrix}, 
		\\
		\trans{N} &= 
		\begin{bmatrix}
			\trans{R} & \trans{R} & \trans{R}
		\end{bmatrix}, \quad
		\ve{R} = 
		\begin{bmatrix}
			R_{11} & R_{12} & ... & R_{IJ} & ... & R_{NM}
		\end{bmatrix}.
	\end{aligned}
\end{equation}
\noindent The velocity field is obtained as the material derivative of the displacement field:
\begin{equation}
	\label{eq:def:velocity}
	\ve{v} = \mdvdef{r} = \mdv{u} = \md{u}^n \ve{i}_n = v^n \ve{i}_n = \md{x}^{*n} \ve{i}_n,
\end{equation}
\noindent and it is represented analogously: 
\begin{equation}
	\label{eq:def:matrices for velocity}
	\ve{v} = \mdv{u} = \ve{N} \mdv{q}, \;\;
	\trans{$\mdv{q}$} = 
		\begin{bmatrix}
		\trans{$\mdv{u}_{11}$} & \trans{$\mdv{u}_{12}$} & ... & \trans{$\mdv{u}_{IJ}$} & ... & \trans{$\mdv{u}_{NM}$}
		\end{bmatrix}, \;\;
	\mdv{u}_{IJ} = 
		\begin{bmatrix}
		v^1_{IJ} & v^2_{IJ} & v^3_{IJ}
		\end{bmatrix}.
\end{equation}

The strain rate tensor is the symmetric part of the velocity gradient and it equals the material derivative of \eqref{eq:def:strain}: 
\begin{equation}
	\label{eq:def:strain rate}
	\ii{d}{}{\alpha \beta} = \imd{\epsilon}{}{\alpha \beta} = \frac{1}{2} \imddef{g}{}{\alpha \beta} = \frac{1}{2} \left( \iloc{v}{}{\alpha \vert \beta} + \iloc{v}{}{\beta \vert \alpha} \right) = \frac{1}{2} \left( \iloc{v}{}{\alpha , \beta} + \iloc{v}{}{\beta , \alpha} - 2 \iloc{v}{}{m} \idef{\Gamma}{m}{\alpha \beta} \right),
\end{equation}
\noindent where $ \iloc{v}{}{m} $ are the components of the velocity vector with respect to the local curvilinear coordinates, which are related to the global ones as:
\begin{equation}
	\label{eq:def:relation global local}
		\ve{v} = \iloc{v}{}{m} \ivdef{g}{m}{} = \iloc{v}{m}{} \ivdef{g}{}{m} = \ii{v}{}{k} \iv{i}{k}{} = \ii{v}{k}{} \iv{i}{}{k} \Rightarrow
		\begin{cases}
			\iloc{v}{}{m} &= \ii{v}{}{k} \iv{i}{k}{} \cdotp \ivdef{g}{}{m} = \idef{x}{k}{,m} \ii{v}{}{k} \\ 
			\iloc{v}{m}{} &= \ii{v}{k}{} \iv{i}{}{k} \cdotp \ivdef{g}{m}{} = \idef{x}{,m}{k} \ii{v}{k}{}
		\end{cases}.
\end{equation}
\noindent The covariant components of the strain rate with respect to the Cartesian coordinates can be written as:
\begin{equation}
	\label{eq:def: strain rate wrt Cartesian}
	\begin{aligned}
	\ii{d}{}{\alpha \beta} &= \frac{1}{2} \left( \idef{x}{k}{,\alpha} \ii{v}{}{k,\beta} + \idef{x}{k}{,\beta} \ii{v}{}{k,\alpha} + 2 \idef{x}{k}{,\alpha\beta} \ii{v}{}{k} - 2 \idef{\Gamma}{m}{\alpha\beta} \idef{x}{n}{,m} \ii{v}{}{n} \right) \\
	&= \frac{1}{2} \left( \idef{x}{k}{,\alpha} \ii{v}{}{k,\beta} + \idef{x}{k}{,\beta} \ii{v}{}{k,\alpha} \right) = \frac{1}{2} \left( \ivdef{g}{}{\alpha} \cdotp \iv{v}{}{\beta} + \ivdef{g}{}{\beta} \cdotp \iv{v}{}{\alpha} \right),
	\end{aligned}
\end{equation}
\noindent while the trace of the strain rate tensor is: 
\begin{equation}
	\label{eq:trace of strain rate tens}
	\ii{d}{\alpha}{\alpha} = \idef{g}{\alpha\beta}{} \ii{d}{}{\alpha \beta} = \frac{1}{2} \left( \idef{x}{,\beta}{n} \ii{v}{n}{,\beta} + \idef{x}{,\alpha}{n} \ii{v}{n}{,\alpha} \right) = \idef{x}{,\alpha}{n} \ii{v}{n}{,\alpha} = \ivdef{g}{\alpha}{} \cdotp \iv{v}{}{\alpha}.
\end{equation}
\noindent The material derivative of \eqref{eq:def:r equidistant def} gives the velocity vector of an equidistant surface:
\begin{equation}
	\label{eq:velocity at eq surface}
	\iveq{v}{}{} = \ve{v} (\zeta) = \ve{v} + \zeta \iv{v}{}{3},
\end{equation}
\noindent where the following notation is introduced:
\begin{equation}
	\label{eq:notation for veloc at eq surf}
	\iveq{v}{}{} = \iveqmddef{r}{}{} = \ieq{v}{n}{} \iv{i}{}{n} = \ieqmddef{x}{n}{} \iv{i}{}{n}, \; \iv{v}{}{m} = \ivmddef{g}{}{m} = \ii{v}{n}{,m} \iv{i}{}{n} = \imddef{x}{n}{,m} \iv{i}{}{n}.
\end{equation}
\noindent Thus, the vector $ \iv{v}{}{3} $ is the material derivative of the deformed normal:
\begin{equation}
	\label{eq:def v3}
	\begin{aligned}
	\iv{v}{}{3} &= \ivmddef{g}{}{3} = \ii{v}{}{k,3} \iv{i}{k}{} = \frac{1}{\sqrt{\idef{g}{}{}}} \left( \iv{v}{}{1} \times \ivdef{g}{}{2} + \ivdef{g}{}{1} \times \iv{v}{}{2} - \frac{\imddef{g}{}{}}{2\sqrt{\idef{g}{}{}}} \ivdef{g}{}{3} \right), \\
	\ii{v}{}{k,3} &= \frac{1}{\sqrt{\idef{g}{}{}}} \left[ \left( \ii{v}{m}{,1} \idef{x}{n}{,2} + \idef{x}{m}{,1} \ii{v}{n}{,2} \right) \ii{e}{}{mnk} - \frac{\imddef{g}{}{}}{2\sqrt{\idef{g}{}{}}} \idef{x}{}{k,3} \right],
	\end{aligned}
\end{equation}
\noindent and its derivation requires special attention. Firstly, the material derivative of the determinant of the metric tensor, utilizing Eqs.~\eqref{eq:def:met_tensor}, \eqref{eq:def:strain rate}, and \eqref{eq:def: strain rate wrt Cartesian}, is computed: 
\begin{equation}
	\label{eq:calc: mat der of a determ of metric tesor}
	\begin{aligned}
	\imddef{g}{}{} &= 2 \left( \ii{d}{}{11} \idef{g}{}{22} + \ii{d}{}{22} \idef{g}{}{11} - 2 \ii{d}{}{12} \idef{g}{}{12} \right) \\
	&= 2 \left[ \left( \ivdef{g}{}{1} \cdotp \iv{v}{}{1} \right) \idef{g}{}{22} + \left( \ivdef{g}{}{2} \cdotp \iv{v}{}{2} \right) \idef{g}{}{11} - 2 \left( \ivdef{g}{}{1} \cdotp \iv{v}{}{2} + \ivdef{g}{}{2} \cdotp \iv{v}{}{1}\right) \idef{g}{}{12} \right],
	\end{aligned}
\end{equation}
\noindent which leads to:
\begin{equation}
	\label{eq:ratio g/2g equals trace of strain rate tensor}
	\frac{\imddef{g}{}{}}{2\idef{g}{}{}} = \ii{d}{}{11} \idef{g}{11}{} + \ii{d}{}{12} \idef{g}{12}{} + \ii{d}{}{21} \idef{g}{21}{} + \ii{d}{}{22} \idef{g}{22}{} = \ii{d}{1}{1} + \ii{d}{2}{2} = \ii{d}{\alpha}{\alpha}.
\end{equation}
\noindent If we represent the tangent vectors of the deformed configuration as:
\begin{equation}
	\label{eq:represent tang vectors of def config}
	\begin{aligned}
	\ivdef{g}{}{1} &= \idef{x}{n}{,1} \iv{i}{}{n} = \sqrt{\idef{g}{}{}} \left( \ivdef{g}{2}{} \times \ivdef{g}{}{3}\right), \quad \idef{x}{n}{,1} = \sqrt{\idef{g}{}{}} \idef{x}{,2}{k} \idef{x}{}{l,3} \ii{e}{kln}{},\\
	\ivdef{g}{}{2} &= \idef{x}{n}{,2} \iv{i}{}{n} = \sqrt{\idef{g}{}{}} \left( \ivdef{g}{}{3} \times \ivdef{g}{1}{}\right), \quad \idef{x}{n}{,2} = \sqrt{\idef{g}{}{}} \idef{x}{}{k,3} \idef{x}{,1}{l} \ii{e}{kln}{}, 
	\end{aligned}	
\end{equation}
\noindent the first part of \eqqref{eq:def v3} can be rewritten as: 
\begin{equation}
	\label{eq:first part of v3}
	\begin{aligned}
	\frac{1}{\sqrt{\idef{g}{}{}}} \left( \iv{v}{}{1} \times \ivdef{g}{}{2} + \ivdef{g}{}{1} \times \iv{v}{}{2} \right) &= \iv{v}{}{1} \times \left( \ivdef{g}{}{3} \times \ivdef{g}{1}{} \right) - \iv{v}{}{2} \times \left( \ivdef{g}{2}{} \times \ivdef{g}{}{3} \right) \\
	&= \left( \ivdef{g}{\alpha}{} \cdotp \iv{v}{}{\alpha} \right) \ivdef{g}{}{3} - \left( \ivdef{g}{}{3} \cdotp \iv{v}{}{\alpha} \right) \ivdef{g}{\alpha}{}.
	\end{aligned}
\end{equation}
\noindent By the insertion of Eqs.~\eqref{eq:ratio g/2g equals trace of strain rate tensor}, \eqref{eq:trace of strain rate tens}, and \eqref{eq:first part of v3} into \eqref{eq:def v3}, we can finally obtain the velocity of the normal:
\begin{equation}
	\label{eq:v3 final}
	\begin{aligned}
	\iv{v}{}{3} &= \left( \ivdef{g}{\alpha}{} \cdotp \iv{v}{}{\alpha} \right) \ivdef{g}{}{3} - \left( \ivdef{g}{}{3} \cdotp \iv{v}{}{\alpha} \right) \ivdef{g}{\alpha}{} - \left( \ivdef{g}{\alpha}{} \cdotp \iv{v}{}{\alpha} \right) \ivdef{g}{}{3} = - \left( \ivdef{g}{}{3} \cdotp \iv{v}{}{\alpha} \right) \ivdef{g}{\alpha}{}, \\
	\ii{v}{}{k,3} &= \idef{x}{,\alpha}{m} \ii{v}{m}{,\alpha} \idef{x}{}{k,3} - \idef{x}{,\alpha}{k} \ii{v}{m}{,\alpha} \idef{x}{}{m,3} - \idef{x}{,\alpha}{m} \ii{v}{m}{,\alpha} \idef{x}{}{k,3} = - \idef{x}{,\alpha}{k} \ii{v}{m}{,\alpha} \idef{x}{}{m,3},
	\end{aligned}	
\end{equation}
\noindent which is the same result as calculated in \cite{1973naghdi}. Finally, the components of velocity at an equidistant surface \eqref{eq:velocity at eq surface} can be represented as:
\begin{equation}
	\label{eq:velocity at eq surface final}
	\iveq{v}{}{} = \ieq{v}{}{n} \iv{i}{n}{} = \ve{v} - \zeta \left( \ivdef{g}{}{3} \cdotp \iv{v}{}{\alpha} \right) \ivdef{g}{\alpha}{}, \; \ieq{v}{}{n} = \ii{v}{}{n} - \zeta \idef{x}{,\alpha}{n} \ii{v}{m}{,\alpha} \idef{x}{}{m,3}.
\end{equation}

\noindent \textbf{Remark.} Proving that the line element remains perpendicular to the midsurface at each configuration is now straightforward. That is, multiplying \eqqref{eq:v3 final} with $ \idef{x}{k}{,\alpha} $ gives:

\begin{equation}
	\label{eq:prove ortho line elem each config}
	\idef{x}{k}{,\alpha} \ii{v}{}{k,3} = - \idef{x}{k}{,\alpha} \idef{x}{,\beta}{k} \ii{v}{m}{,\beta} \idef{x}{}{m,3} = - \ii{\delta}{\beta}{\alpha} \ii{v}{m}{,\beta} \idef{x}{}{m,3} = - \ii{v}{m}{,\alpha} \idef{x}{}{m,3},
\end{equation}

\noindent and the shear strain rates in the planes perpendicular to the tangent planes of the midsurface vanish:

\begin{equation}
	\label{eq:prove shear strain perp to tang palnes vanishes}
	2 \ii{d}{}{\alpha 3} = \idef{x}{n}{,\alpha} \ii{v}{}{n,3} + \idef{x}{n}{,3} \ii{v}{}{n,\alpha}=0.
\end{equation}

\noindent Additionally, the strain rate along the normal also vanishes at each configuration:

\begin{equation}
	\label{eq:prove shear strain along normal vanishes}
	\ii{d}{}{33}= \idef{x}{n}{,3} \ii{v}{}{n,3}= - \idef{x}{n}{,3} \idef{x}{,\alpha}{n} \ii{v}{k}{,\alpha} \idef{x}{}{k,3} = -\ii{\delta}{\alpha}{3} \ii{v}{k}{,\alpha} \idef{x}{}{k,3}=0.
\end{equation}

\subsection{Strain rate at an equidistant surface}

The crucial step for the derivation of structural theories is the definition of the strain at each point of a structure via chosen reference quantities. It allows a dimensional reduction from 3D to 2D. The strain rate at an equidistant surface is, analogously to \eqqref{eq:def: strain rate wrt Cartesian}:
\begin{equation}
	\label{eq:def: strain rate equidistant wrt Cartesian}
	\ieq{d}{}{\alpha \beta} = \ii{d}{}{\alpha \beta} (\zeta) = \frac{1}{2} \left( \iveqdef{g}{}{\alpha} \cdotp \iveq{v}{}{\beta} + \iveqdef{g}{}{\beta} \cdotp \iveq{v}{}{\alpha} \right) = \frac{1}{2} \left( \ieqdef{x}{n}{,\alpha} \ieq{v}{}{n,\beta} + \ieqdef{x}{n}{,\beta} \ieq{v}{}{n,\alpha} \right),
\end{equation}
\noindent where the components of the velocity gradient at an arbitrary point, $ \ieq{v}{}{n,\beta} $, follow from \eqqref{eq:velocity at eq surface final}:
\begin{equation}
	\label{eq: velocity gradient equidistant}
	\ieq{v}{}{n,\alpha} = \ii{v}{}{n,\alpha} + \zeta \ii{v}{}{n,3\alpha},
\end{equation}
\noindent while the derivative of the velocity of normal stems from \eqqref{eq:v3 final}:
\begin{equation}
\label{eq: derivative of v3}
\begin{aligned}
	\ii{v}{}{n,3\alpha} &= -\left( \idef{x}{,\mu}{n,\alpha} \ii{v}{}{k,\mu} \idef{x}{k}{,3} + \idef{x}{,\mu}{n} \ii{v}{}{k,\mu \alpha} \idef{x}{k}{,3} + \idef{x}{,\mu}{n} \ii{v}{}{k,\mu} \idef{x}{k}{,3\alpha} \right) \\
	&= \idef{\Gamma}{\mu}{\alpha\nu} \idef{x}{,\nu}{n} \ii{v}{}{k,\mu} \idef{x}{k}{,3} - \idef{b}{\mu}{\alpha} \idef{x}{,3}{n} \ii{v}{}{k,\mu} \idef{x}{k}{,3} - \idef{x}{,\mu}{n} \ii{v}{}{k,\mu \alpha} \idef{x}{k}{,3} + \idef{x}{,\mu}{n} \ii{v}{}{k,\mu} \idef{b}{\nu}{\alpha} \idef{x}{k}{,\nu}.
\end{aligned}
\end{equation}

\noindent \textbf{Remark.} Utilizing the last two expressions, we can show that the shear strain rate $ \ii{d}{}{\alpha 3}$ is zero at an equidistant surface, as well as at the midsurface:
\begin{equation}
\label{eq: prove shear is zero at eq seurface}
\begin{aligned}
	2 \ieq{d}{}{\alpha 3} &= \left( \iveqdef{g}{}{\alpha} \cdotp \iv{v}{}{3} + \ivdef{g}{}{3} \cdotp \iv{v}{}{\alpha} \right) = \ieqdef{x}{n}{,\alpha} \ii{v}{}{n,3} + \idef{x}{n}{,3} \ieq{v}{}{n,\alpha} \\
	&= \left( \ii{\delta}{\mu}{\alpha} - \zeta \idef{b}{\mu}{\alpha} \right) \idef{x}{n}{,\mu} \ii{v}{}{n,3} + \idef{x}{n}{,3} \left( \ii{v}{}{n,\alpha} + \zeta \ii{v}{}{n,3\alpha} \right) \\
	&= \idef{x}{n}{,\alpha} \ii{v}{}{n,3} - \zeta \idef{b}{\mu}{\alpha} \idef{x}{n}{,\mu} \ii{v}{}{n,3} + \idef{x}{n}{,3} \ii{v}{}{n,\alpha} \\
	&\;\;\;\; + \zeta \idef{x}{n}{,3} \left( \idef{\Gamma}{\mu}{\alpha \nu} \idef{x}{,\nu}{n} \ii{v}{}{k,\mu} \idef{x}{k}{,3} - \idef{b}{\mu}{\alpha} \idef{x}{,3}{n} \ii{v}{}{k,\mu} \idef{x}{k}{,3} - \idef{x}{,\mu}{n} \ii{v}{}{k,\mu\alpha} \idef{x}{k}{,3} -  \idef{x}{,\mu}{n} \ii{v}{}{k,\mu} \idef{b}{\nu}{\alpha} \idef{x}{k}{,\nu} \right) \\
	&= \idef{x}{n}{,\alpha} \ii{v}{}{n,3} - \zeta \idef{b}{\mu}{\alpha} \idef{x}{n}{,\mu} \ii{v}{}{n,3} +  \idef{x}{n}{,3} \ii{v}{}{n,\alpha} - \zeta \idef{b}{\mu}{\alpha} \ii{v}{}{k,\mu} \idef{x}{k}{,3} \\
	&= \idef{x}{n}{,\alpha} \ii{v}{}{n,3} + \idef{x}{n}{,3} \ii{v}{}{n,\alpha} - \zeta \idef{b}{\mu}{\alpha} \left( \idef{x}{n}{,\mu} \ii{v}{}{n,3} + \idef{x}{n}{,3} \ii{v}{}{n,\mu}  \right) = 0,
\end{aligned}
\end{equation}
\noindent which proves that the line segment is perpendicular to the equidistant surfaces as well as to the midsurface.

Now we have all the ingredients required for the exact relation between equidistant and reference strain rates. Inserting Eqs.~\eqref{eq: velocity gradient equidistant} and \eqref{eq: derivative of v3} into \eqqref{eq:def: strain rate equidistant wrt Cartesian}, returns the components of the equidistant strain rate:
\begin{equation}
\label{eq: components of eq strain rate final}
	\ieq{d}{}{\alpha\beta} = \ii{A}{\mu\nu}{\alpha\beta} \ii{d}{}{\mu\nu} - \zeta \ii{B}{\mu\nu}{\alpha\beta} \imd{\kappa}{}{\mu\nu},
\end{equation}
\noindent where:
\begin{equation}
\label{eq: A and B in eq strain rate}
\begin{aligned}
	\ii{A}{\mu\nu}{\alpha\beta} &= \frac{1}{2} \left[ \ieqdef{C}{\mu}{\alpha} \left( \ii{\delta}{\nu}{\beta} + \zeta \idef{b}{\nu}{\beta} \right) + \ieqdef{C}{\nu}{\beta} \left( \ii{\delta}{\mu}{\alpha} + \zeta \idef{b}{\mu}{\alpha} \right) \right] = \ii{\delta}{\mu}{\alpha} \ii{\delta}{\nu}{\beta} - \zeta^2 \idef{b}{\mu}{\alpha} \idef{b}{\nu}{\beta},\\
	\ii{B}{\mu\nu}{\alpha\beta} &= \frac{1}{2} \left( \ii{\delta}{\nu}{\alpha} \ieqdef{C}{\mu}{\beta} + \ii{\delta}{\mu}{\beta} \ieqdef{C}{\nu}{\alpha} \right),
\end{aligned}
\end{equation}
\noindent while the tensor of the rate of curvature change of the midsurface, $\imd{\kappa}{}{\alpha\beta}$, is introduced as:
\begin{equation}
\label{eq: tensor of curvature change -rate}
\begin{aligned}
	\imd{\kappa}{}{\alpha\beta}&=\imddef{b}{}{\alpha\beta}=\ivmddef{g}{}{\alpha,\beta} \cdotp \ivdef{g}{}{3} + \ivdef{g}{}{\alpha,\beta} \cdotp \ivmddef{g}{}{3} = \iv{v}{}{\alpha,\beta} \cdotp \ivdef{g}{}{3} + \ivdef{g}{}{\alpha,\beta} \cdotp \iv{v}{}{3} = \ii{v}{}{k,\alpha\beta} \idef{x}{k}{,3} + \idef{\Gamma}{n}{\alpha \beta} \idef{x}{k}{,n} \ii{v}{}{k,3}\\
	&= \idef{x}{k}{,3} \ii{v}{}{k,\alpha\beta} - \idef{\Gamma}{\nu}{\alpha\beta} \idef{x}{k}{,\nu} \idef{x}{,\mu}{k} \ii{v}{m}{,\mu} \idef{x}{}{m,3} - \idef{b}{}{\alpha\beta} \idef{x}{k}{,3} \idef{x}{,\mu}{k} \ii{v}{m}{,\mu} \idef{x}{}{m,3} = \idef{x}{k}{,3} \left( \ii{v}{}{k,\alpha\beta} - \idef{\Gamma}{\mu}{\alpha\beta} \ii{v}{}{k,\mu} \right).
\end{aligned}
\end{equation}
\noindent The detailed derivation of Eqs.~\eqref{eq: components of eq strain rate final} and \eqref{eq: A and B in eq strain rate} is given in Appendix A. Equations \eqref{eq: components of eq strain rate final} and \eqref{eq: A and B in eq strain rate} define the exact relation of the equidistant and reference strain rates of a KL shell within the scope of finite, but small, strain theory. These expressions can be found in \cite{1973naghdi} where they are derived as a special case of the general shell theory. 

As already noted, many existing IGA shell formulations use a simplified form of the expression \eqref{eq: components of eq strain rate final}. This is justified by the fact that the influence of the second-order terms with respect to the thickness coordinate can be neglected for thin shells, \cite{2015kiendl, 2015nguyen-thanh, 2017casquero}. This assumption corresponds to the following approximation:
\begin{equation}
\label{eq:approximate A and B}
	\ii{A}{\mu\nu}{\alpha\beta} = \ii{B}{\mu\nu}{\alpha\beta} \approx \ii{\delta}{\mu\nu}{\alpha\beta},
\end{equation}
\noindent which reduces \eqqref{eq: components of eq strain rate final} to the linear form:
\begin{equation}
\label{eq:linear form of eq strain}
	\ieq{d}{}{\alpha\beta} = \ii{d}{}{\alpha\beta} - \zeta \imd{\kappa}{}{\alpha\beta}.
\end{equation}
\noindent A notable exception is the reference \cite{2017duong} where one type of simplified, but nonlinear, strain distribution is utilized. However, the comparison of different models is not considered. 

It is interesting to note that the physical components of strain are distributed nonlinearly due to the metric, although the previous equation imposes a linear strain distribution:
\begin{equation}
\label{eq:physical strain}
	\ieq{d}{}{(\alpha\beta)} = \frac{\ieq{d}{}{\alpha\beta}}{\sqrt{\ieq{g}{}{\alpha\alpha} \ieq{g}{}{\beta\beta}}}.
\end{equation}
\noindent The linear analysis performed in \cite{2004bischoff} shows that for thick shells, the utilization of a nonlinear strain distribution and an exact shifter is recommended, along with exact integration. It is also noted that for $ Kh<0.1 $, the error introduced by the assumption of a linear strain distribution is not significant.

Equations \eqref{eq:def: strain rate wrt Cartesian}, \eqref{eq:prove shear strain perp to tang palnes vanishes}, \eqref{eq:prove shear strain along normal vanishes}, \eqref{eq: components of eq strain rate final}, and \eqref{eq: A and B in eq strain rate} are the basis for the formulation of the isogeometric element for the finite strain and finite rotation analysis of an arbitrarily curved KL shell.

\section{Finite element formulation}

\textcolor{black}{In line with the previous derivation, we will formulate isogeometric finite elements of the KL shell using the principle of virtual power. The present formulation follows the approach given in \cite{2017radenkovic}.} The generalized coordinates are the components of the velocities of the control points.

Let us start from the generalized Hooke law for the linear elastic material, also known as the Saint Venant-Kirchhoff material model. This is well-suited for the small strain and large rotation analysis, \cite{2004bischoff}. The components of the stress and strain rates are related as:
\begin{equation}
\label{eq:stress strain relation}
	\imd{\sigma}{i}{j} = 2 \mu \ii{d}{i}{j} + \lambda \ii{\delta}{i}{j} \ii{d}{k}{k},
\end{equation}
\noindent where $\mu$ and $\lambda$ are the Lamé material parameters. The plane stress assumption states that the normal component of the stress rate in the direction of the shell thickness is negligible compared to the other components:
\begin{equation}
\label{eq:stress 3 equals zero}
	\imd{\sigma}{3}{3} = 2 \mu \ii{d}{3}{3} + \lambda \left( \ii{d}{\alpha}{\alpha} + \ii{d}{3}{3} \right) = 0,
\end{equation}
\noindent and the corresponding strain rate becomes:
\begin{equation}
\label{eq:strain 3 equals zero}
	\ii{d}{3}{3} = -\frac{\nu}{1-\nu} \ii{d}{\alpha}{\alpha}.
\end{equation}
\noindent The existence of the $\ii{d}{3}{3}$ component violates the KL hypothesis and it is neglected here, together with the $\imd{\sigma}{3}{3}$. This is a known inconsistency of the KL theory caused by its approximate nature. Nevertheless, a great quantity of literature accounts for this component of strain rate since it increases the accuracy of the mechanical model and is easily introduced via the static condensation of the material tensor, \cite{2015kiendl}. 

Here, the plane stress and the plane strain conditions are enforced and the stress-strain relation takes its final form:
\begin{equation}
\label{eq:stress strain final form}
	\imd{\sigma}{\alpha\beta}{} = 2 \mu \left( \idef{g}{\alpha\nu}{} \idef{g}{\beta\gamma}{} + \frac{\nu}{1-\nu} \idef{g}{\alpha\beta}{} \idef{g}{\nu\gamma}{}\right) \ii{d}{}{\nu\gamma} = \ii{D}{\alpha\beta\nu\gamma}{} \ii{d}{}{\nu\gamma},
\end{equation}
\noindent where $\ii{D}{\alpha\beta\nu\gamma}{}$ are the components of the constitutive, material, tensor.

\subsection{Principle of virtual power}

The principle of virtual power states that at any time, the total power of the external, internal and inertial forces is zero for any admissible virtual state of motion. If the inertial effects are neglected and only the surface loads are considered, it can be written as: 
\begin{equation}
\label{eq:virtual power}
	\delta P = \int_{v}^{} \ii{\sigma}{ij}{} \delta \ii{d}{}{ij} \dd{v} - \int_{a}^{} \iloc{p}{i}{} \delta \iloc{v}{}{i} \dd{a} = \int_{v}^{} \bm{\sigma} : \delta \ve{d} \dd{v} - \int_{a}^{} \ve{p} \delta \ve{v} \dd{a} = 0,
\end{equation}
\noindent where $\bm{\sigma} $ is the Cauchy stress tensor, \ve{d} is the strain rate tensor, and \ve{p} is the vector of external surface loads. All these quantities are measured with respect to the current configuration. This configuration is unknown and its determination is the main subject of the nonlinear analysis. For problems with deformation-independent loads, it is only required to linearize the Cauchy stress. At the current $(n+1)$ configuration this stress is approximated by:
\begin{equation}
\label{eq:linearization of stress}
\begin{aligned}
	\iii{\bm{\sigma}}{(n+1)}{}{}{} &= \iii{\bm{\sigma}}{(n)}{(n+1)}{}{} + \frac{{\operatorname{d}}^t \: \iii{\bm{\sigma}}{(n)}{}{}{}} {\operatorname{d}t} \Delta \ii{t}{}{(n+1)}, \\ 
	\iiieq{\sigma}{(n)}{(n+1)}{}{\alpha\beta} &= \sqrt{\frac{\iiieq{g}{(n)}{}{}{}}{\iiieq{g}{(n+1)}{}{}{}}} \: \iiieq{\sigma}{(n)}{}{}{\alpha\beta} = \frac{\iii{g}{(n)}{}{}{0}}{\iii{g}{(n+1)}{}{}{0}} \sqrt{\frac{\iii{g}{(n)}{}{}{}}{\iii{g}{(n+1)}{}{}{}}}  \: \iiieq{\sigma}{(n)}{}{}{\alpha\beta},
\end{aligned}
\end{equation}
\noindent where $\iiieq{\sigma}{(n)}{(n+1)}{}{\alpha\beta}$ are the components of the equidistant stress from the previous $(n)$ configuration expressed with respect to the metric of the current $(n+1)$ configuration. Similarly, $\iiieq{\sigma}{(n)}{}{}{\alpha\beta}$ is the same quantity but expressed with respect to the previous $(n)$ configuration. The additive decomposition of stress, as in \eqqref{eq:linearization of stress}, is only valid if all of the terms are expressed with respect to the same metric. However, the rigorous integration of this expression is far from trivial due to the ratio of shifters in two different configurations. If the incremental change of curvature is small, this ratio can be approximated as $\iii{g}{(n)}{}{}{0}/\iii{g}{(n+1)}{}{}{0} \approx 1 $. An appropriate time derivative is designated with $\operatorname{d}^t () / \operatorname{d}t$ while $\Delta \ii{t}{}{(n+1)} = \ii{t}{}{(n+1)} - \ii{t}{}{(n)}$ is the time increment. 

The rate form of stress-strain relations requires an objective time derivative. Note that the material derivative is not objective, while the corotational (Jaumann) and convective time derivatives fulfill this property, \cite{2008wriggers, 1984johnson}. The relation \eqref{eq:linearization of stress} is rewritten as: 
\begin{equation}
\label{eq:linearization of stress rewritten}
	\iii{\bm{\sigma}}{(n+1)}{}{}{} = \iii{\bm{\sigma}}{(n)}{}{}{} + \left( \ii{L}{}{v} \: \iii{\bm{\sigma}}{(n)}{}{}{}\right) \Delta \ii{t}{}{(n+1)}
\end{equation}
\noindent where $\ii{L}{}{v}$ designates Lie time derivative which matches the convective derivative, for convective coordinates. With respect to the covariant metric base it is:
\begin{equation}
\label{eq:Lie derivative wrt covariant metric}
	\ii{L}{}{v} \bm{\sigma} = \md{\bm{\sigma}} - \ve{L} \bm{\sigma} - \bm{\sigma} \trans{L}, \quad \ve{L} = \iloc{v}{}{i \vert j} \iv{g}{i}{} \otimes \iv{g}{j}{},
\end{equation}
\noindent where $\ve{L}$ is the spatial velocity gradient. On the other hand, since the material time derivative of the stress tensor reads:
\begin{equation*}
\label{eq:material derivative of stress tensor}
	\md{\bm{\sigma}} = \imd{\sigma}{ij}{} \iv{g}{}{i} \otimes \iv{g}{}{j} + \ve{L} \bm{\sigma} + \bm{\sigma} \trans{L},
\end{equation*}
\noindent relation \eqref{eq:Lie derivative wrt covariant metric} returns: 
\begin{equation}
\label{eq:Lie derivative equals material}
	\ii{L}{}{v} \bm{\sigma} = \imd{\sigma}{ij}{} \iv{g}{}{i} \otimes \iv{g}{}{j} \Rightarrow \left(\ii{L}{}{v} \bm{\sigma}\right)^{ij} = \imd{\sigma}{ij}{}.
\end{equation}
\noindent Therefore, the components of the stress rate tensor equal the material derivative of the components of the stress tensor. This is an important fact which will be utilized further.

After the insertion of Eqs.~\eqref{eq:linearization of stress rewritten} and \eqref{eq:Lie derivative equals material} into \eqqref{eq:virtual power}, the linearized form of the principle of virtual power at the current configuration is obtained:
\begin{equation}
\label{eq:linearized virtual power}
\begin{aligned}
	&\int_{v}^{} \md{\bm{\sigma}} : \delta \ve{d} \dd{v} \Delta t + \int_{v}^{} \bm{\sigma} : \delta \ve{d} \dd{v} = \int_{a}^{} \ve{p} \delta \ve{v} \dd{a}, \\
	&\int_{v}^{} \ieqmd{\sigma}{\alpha\beta}{} \delta \ieq{d}{}{\alpha\beta} \dd{v} \Delta t + \int_{v}^{} \ieq{\sigma}{\alpha\beta}{} \delta \ieq{d}{}{\alpha\beta} \dd{v} = \int_{a}^{} \iloc{p}{i}{} \delta \iloc{v}{}{i} \dd{a}.
\end{aligned}
\end{equation}
\noindent From this equation and further in the text, the notation is slightly simplified for the sake of readability by removing the time indices and asterisks. This can be done without introducing any notational ambiguity since (i) the stress and strain rates are instantaneous quantities, while the known stress $\bm{\sigma}$ is calculated at the previous configuration and transformed to the metric of the current configuration, see Eqs.~\eqref{eq:linearization of stress} and \eqref{eq:linearized virtual power}, and (ii) all integrations are performed with respect to the metric of the current configuration, in accordance with the updated Lagrangian procedure, \cite{2007bathe}.

In order to reduce the dimension from 3D to 2D, it is necessary to integrate the left-hand side of \eqqref{eq:linearized virtual power} along the thickness. Thus, integrals over the area of midsurface are obtained: 
\begin{equation}
\label{eq:from 3D to 2D}
\begin{aligned}
	\int_{v}^{} \ieqmd{\sigma}{\alpha\beta}{} \delta \ieq{d}{}{\alpha\beta} \dd{v} \Delta t + \int_{v}^{} \ieq{\sigma}{\alpha\beta}{} \delta \ieq{d}{}{\alpha\beta} \dd{v} = &\int_{a}^{} \left( \imd{N}{\mu\nu}{} \delta \ii{d}{}{\mu\nu} + \imd{M}{\mu\nu}{} \delta \imd{\kappa}{}{\mu\nu} \right) \dd{a} \Delta t \\
	&+ \int_{a}^{} \left( \ii{N}{\mu\nu}{} \delta \ii{d}{}{\mu\nu} + \ii{M}{\mu\nu}{} \delta \imd{\kappa}{}{\mu\nu} \right) \dd{a},
\end{aligned}
\end{equation}
\noindent where $\ii{N}{\mu\nu}{}$ and  $\ii{M}{\mu\nu}{}$ are stress resultants and stress couples, which are energetically conjugated with the reference strain rates of the midsurface, $\ii{d}{}{\mu\nu}$ and $\imd{\kappa}{}{\mu\nu}$, while $\imd{N}{\mu\nu}{}$ and $\imd{M}{\mu\nu}{}$ are their respective rates: 
\begin{equation}
\label{eq: rates of section forces}
\begin{aligned}
	\imd{N}{\mu\nu}{} &= \int_{-h/2}^{h/2} \ieqmd{\sigma}{\alpha\beta}{} \ii{A}{\mu\nu}{\alpha\beta} \ii{g}{}{0} \dd{\zeta} = \int_{-h/2}^{h/2} \ieqmd{\sigma}{\alpha\beta}{} \left( \ii{\delta}{\mu}{\alpha} \ii{\delta}{\nu}{\beta} - \zeta^2 \ii{b}{\mu}{\alpha} \ii{b}{\nu}{\beta} \right) \ii{g}{}{0} \dd{\zeta},\\
	\imd{M}{\mu\nu}{} &= - \frac{1}{2} \int_{-h/2}^{h/2} \zeta \ieqmd{\sigma}{\alpha\beta}{} \ii{B}{\mu\nu}{\alpha\beta} \ii{g}{}{0} \dd{\zeta} = - \frac{1}{2} \int_{-h/2}^{h/2} \zeta \ieqmd{\sigma}{\alpha\beta}{} \left(\ii{\delta}{\nu}{\alpha} \ieq{C}{\mu}{\beta} + \ii{\delta}{\mu}{\beta} \ieq{C}{\nu}{\alpha}\right) \ii{g}{}{0} \dd{\zeta}.
\end{aligned}
\end{equation}
\noindent These expressions are the same as those obtained in \cite{1973koiter} but they are rarely applied for the nonlinear analysis due to the complexity. If we introduce the vectors of generalized section forces and strain rates of the shell midsurface: 
\begin{equation}
\label{eq:vectors of section forces and stran rates}
\begin{aligned}
	\trans{f} &= 
	\begin{bmatrix}
	\ii{N}{11}{} & \ii{N}{22}{} & \ii{N}{12}{} & \ii{M}{11}{} & \ii{M}{22}{} & \ii{M}{12}{}
	\end{bmatrix},\\
	\trans{e} &= 
	\begin{bmatrix}
	\ii{d}{}{11} & \ii{d}{}{22} & 2 \ii{d}{}{12} & \imd{\kappa}{}{11} & \imd{\kappa}{}{22} & 2 \imd{\kappa}{}{12}
	\end{bmatrix},
\end{aligned}
\end{equation}	
\noindent \eqqref{eq:linearized virtual power} can be written in compact matrix form as:
\begin{equation}
\label{matrix form of linearized virtual power}
	\int_{a}^{} \transmd{f} \delta \ve{e} \dd{a} \Delta t + \int_{a}^{} \trans{f} \delta \ve{e} \dd{a} = \int_{a}^{} \ve{p} \delta \ve{v} \dd{a}.
\end{equation}

The geometrically exact relations \eqref{eq:from 3D to 2D} and \eqref{eq: rates of section forces} are crucial for the accurate formulation of structural shell theories. In particular, energetically conjugated pairs are defined rigorously and the appropriate constitutive matrices are guaranteed to be symmetric, as will be shown in the next section.

\subsection{Analytical relation between energetically conjugated pairs}

The integration represented with Eqs.~\eqref{eq:from 3D to 2D} and \eqref{eq: rates of section forces} is an essential step for the development of the structural theory of shells, since it allows the calculation of virtual power by the reference quantities of the midsurface. This requires the constitutive relation between the stress and strain rates at an arbitrary point of the shell continuum. Analogously to \eqqref{eq:stress strain final form}, this relation is:
\begin{equation}
\label{eq:constitutive relation at equidistant surface}
	\ieqmd{\sigma}{\alpha\beta}{} = 2 \mu \left( \ieq{g}{\alpha\nu}{} \ieq{g}{\beta\gamma}{} + \frac{\nu}{1-\nu} \ieq{g}{\alpha\beta}{} \ieq{g}{\nu\gamma}{} \right) \ieq{d}{}{\nu\gamma} = \ieq{D}{\alpha\beta\nu\gamma}{} \ieq{d}{}{\nu\gamma}.
\end{equation}
\noindent Using \eqqref{eq:rel:met_tens_eq_via_ref}, it is evident that the following transformation is valid:
\begin{equation}
\label{eq: relation of constitutuive tensor at equidistant surface and the reference one}
\begin{aligned}
	\ieq{D}{\mu\nu\epsilon\varphi}{} &= 2 \mu \left( \ii{C}{\mu}{\rho} \ii{C}{\epsilon}{\delta} \ii{g}{\rho \delta}{} \ii{C}{\nu}{\alpha} \ii{C}{\varphi}{\beta} \ii{g}{\alpha\beta}{} + \frac{\nu}{1-\nu} \ii{C}{\mu}{\rho} \ii{C}{\nu}{\alpha} \ii{g}{\rho \alpha}{} \ii{C}{\epsilon}{\delta} \ii{C}{\varphi}{\beta} \ii{g}{\delta\beta}{} \right) \\
	&= \ii{C}{\mu}{\rho} \ii{C}{\nu}{\alpha} \ii{C}{\epsilon}{\delta}\ii{C}{\varphi}{\beta} 2 \mu\left(\ii{g}{\rho \delta}{} \ii{g}{\alpha \beta}{} + \frac{\nu}{1-\nu} \ii{g}{\rho \alpha}{} \ii{g}{\delta\beta}{} \right) = \ii{C}{\mu}{\rho} \ii{C}{\nu}{\alpha} \ii{C}{\epsilon}{\delta}\ii{C}{\varphi}{\beta} \ii{D}{\rho\alpha\delta\beta}{}.
\end{aligned}
\end{equation}
\noindent This relation explicitly relates the components of constitutive tensor at an arbitrary point with the ones at the midsurface via the reciprocal shift tensor. In this way, the constitutive relations between the energetically conjugated rates of generalized section forces and reference strain rates of the midsurface are:
\begin{equation}
\label{eq: section forces analytic}
\begin{aligned}
	\imd{N}{\alpha\beta}{} &= \ii{D}{\alpha\beta\gamma\lambda}{M} \ii{d}{}{\gamma\lambda} + \ii{D}{\alpha\beta\gamma\lambda}{MB} \imd{\kappa}{}{\gamma\lambda}, \\
	\imd{M}{\alpha\beta}{} &= \ii{D}{\alpha\beta\gamma\lambda}{BM} \ii{d}{}{\gamma\lambda} + \ii{D}{\alpha\beta\gamma\lambda}{B} \imd{\kappa}{}{\gamma\lambda},
\end{aligned}
\end{equation}
\noindent where: 
\begin{equation}
\label{eq: definition of parts of constitutive tensor}
\begin{aligned}
	\ii{D}{\alpha\beta\gamma\lambda}{M} &= \int_{-h/2}^{h/2} \ii{A}{\alpha\beta}{\mu\nu} \ieq{D}{\mu\nu\delta\varphi}{} \ii{A}{\gamma\lambda}{\delta\varphi} \ii{g}{}{0} \dd{\zeta},\\
		\ii{D}{\alpha\beta\gamma\lambda}{MB} &= - \int_{-h/2}^{h/2} \zeta \ii{A}{\alpha\beta}{\mu\nu} \ieq{D}{\mu\nu\delta\varphi}{} \ii{B}{\gamma\lambda}{\delta\varphi} \ii{g}{}{0} \dd{\zeta},\\
			\ii{D}{\alpha\beta\gamma\lambda}{BM} &= - \int_{-h/2}^{h/2} \zeta \ii{B}{\alpha\beta}{\mu\nu} \ieq{D}{\mu\nu\delta\varphi}{} \ii{A}{\gamma\lambda}{\delta\varphi} \ii{g}{}{0} \dd{\zeta} = \ii{D}{\gamma\lambda\alpha\beta}{MB},\\
				\ii{D}{\alpha\beta\gamma\lambda}{B} &= \int_{-h/2}^{h/2} \zeta^2 \ii{B}{\alpha\beta}{\mu\nu} \ieq{D}{\mu\nu\delta\varphi}{} \ii{B}{\gamma\lambda}{\delta\varphi} \ii{g}{}{0} \dd{\zeta}.
\end{aligned}
\end{equation}
\noindent Now, we can represent \eqqref{eq: section forces analytic} in a compact matrix form as:
\begin{equation}
\label{eq: stress strain analytic compact matrix form}
	\ivmd{f}{}{} = \ve{D} \ve{e},
\end{equation}
\noindent where $\ve{D}$ is the symmetric constitutive tensor at the current configuration:
\begin{equation}
\label{eq:constitutive tensor final}
	\ve{D} = 
	\begin{bmatrix}
		\iv{D}{}{M} & \iv{D}{}{MB} \\
		\trans{$\; \iv{D}{}{MB}$} & \iv{D}{}{B}
	\end{bmatrix} = 
	\begin{bmatrix}
		\ii{D}{\alpha\beta\gamma\lambda}{M} & \ii{D}{\alpha\beta\gamma\lambda}{MB} \\
		\ii{D}{\gamma\lambda\alpha\beta}{MB} & \ii{D}{\alpha\beta\gamma\lambda}{B}
	\end{bmatrix} 
	\iv{g}{}{\alpha} \otimes \iv{g}{}{\beta} \otimes \iv{g}{}{\gamma} \otimes \iv{g}{}{\lambda}.
\end{equation}
\noindent Complete expressions for the components of the material tensor $\ve{D}$ are given in Appendix B along with the closed-form solutions of the appropriate integrals. To the best of the authors’ knowledge, this is the first complete analytical representation of the constitutive relation between the energetically conjugated pairs for the elastic KL shell.

\subsection{Variation of strains}

Since the strain rate is a function of the generalized coordinates as well as the metric, we must vary it with respect to both arguments. By noting that the variation of the tangent vectors can be expressed as:
\begin{equation}
\label{variation of base vector g_alpha}
	\delta \iv{g}{}{\alpha} = \delta \iv{v}{}{\alpha} \Delta t,
\end{equation}
\noindent the variations of the reference strains are given by:
\begin{equation}
\label{eq: variations of reference strains}
	\begin{aligned}
	\delta \ii{d}{}{\alpha\beta} &= \delta \frac{1}{2} \left( \iv{g}{}{\alpha} \cdotp \iv{v}{}{\beta} + \iv{g}{}{\beta} \cdotp \iv{v}{}{\alpha}\right) = \frac{1}{2} \left( \delta \iv{v}{}{\alpha} \cdotp \iv{v}{}{\beta} + \delta \iv{v}{}{\beta} \cdotp \iv{v}{}{\alpha}\right) \Delta t + \frac{1}{2} \left( \iv{g}{}{\alpha} \cdotp \delta \iv{v}{}{\beta} + \iv{g}{}{\beta} \cdotp \delta \iv{v}{}{\alpha}\right), \\
	\delta \imd{\kappa}{}{\alpha \beta} &= \delta \left[ \iv{g}{}{3} \cdotp \left( \iv{v}{}{\alpha , \beta} - \ii{\Gamma}{\mu}{\alpha\beta} \iv{v}{}{\mu} \right) \right] \\
	&= \left[ \delta \iv{v}{}{3} \cdotp \left( \iv{v}{}{\alpha,\beta} - \ii{\Gamma}{\lambda}{\alpha\beta} \iv{v}{}{\lambda}\right) - \delta \ii{\Gamma}{\lambda}{\alpha\beta} \left(\iv{g}{}{3} \cdotp \iv{v}{}{\lambda}\right)\right] \Delta t + \iv{g}{}{3} \cdotp \left( \delta \iv{v}{}{\alpha.\beta} - \ii{\Gamma}{\mu}{\alpha\beta} \delta \iv{v}{}{\mu}\right).
	\end{aligned}
\end{equation}
\noindent Most of the terms in the previous expression can be easily computed since the variation is explicitly performed with respect to the unknown variables. However, the first addend of $\imd{\kappa}{}{\alpha\beta}$ is an exception since it requires the variation of the normal and the Christoffel symbols. These variations must be represented via the variations of generalized coordinates. 

The variation of the normal can be performed analogously as the material differentiation given in \eqqref{eq:v3 final}:
\begin{equation}
\label{eq: variation of normal}
	\delta \iv{v}{}{3} = - \left( \iv{g}{}{3} \cdotp \delta \iv{v}{}{\alpha} \right) \iv{g}{\alpha}{}.
\end{equation}
\noindent For the variation of the Christoffel symbols we need the variations of reciprocal base vectors \eqref{eq:def:rec_tang_vec}:
\begin{equation}
\label{eq:variation of reciprocal base vectors}
\begin{aligned}
	\iv{g}{1}{} &=  \left(\iv{g}{}{2} \times \iv{g}{}{3}\right) / \sqrt{g} \Rightarrow \delta \iv{g}{1}{} = - \left[ \left( \delta \iv{v}{}{2} \cdotp \iv{g}{1}{} \right) \iv{g}{2}{} + \left(\delta \iv{v}{}{3} \cdotp \iv{g}{1}{} \right) \iv{g}{3}{} + \left(\delta \iv{v}{}{1} \cdotp \iv{g}{1}{} \right) \iv{g}{1}{}\right] \Delta t, \\
	\iv{g}{2}{} &=  \left(\iv{g}{}{3} \times \iv{g}{}{1}\right) / \sqrt{g} \Rightarrow \delta \iv{g}{2}{} = - \left[ \left( \delta \iv{v}{}{3} \cdotp \iv{g}{2}{} \right) \iv{g}{3}{} + \left( \delta \iv{v}{}{1} \cdotp \iv{g}{2}{} \right) \iv{g}{1}{} + \left(\delta \iv{v}{}{2} \cdotp \iv{g}{2}{} \right) \iv{g}{2}{}\right] \Delta t, 
\end{aligned}
\end{equation}
\noindent which allows us to write: 
\begin{equation}
\label{eq: variation of christoffel}
	\delta \ii{\Gamma}{\lambda}{\alpha\beta} = \delta \left( \iv{g}{}{\alpha,\beta} \cdotp \iv{g}{\lambda}{}\right)	 = \delta \iv{g}{}{\alpha,\beta} \cdotp \iv{g}{\lambda}{} + \iv{g}{}{\alpha,\beta} \cdotp \delta \iv{g}{\lambda}{} = \left[ \iv{g}{\lambda}{} \cdotp \delta \iv{v}{}{\alpha,\beta} - \ii{\Gamma}{\nu}{\alpha\beta} \left(\iv{g}{\lambda}{} \cdot \delta \iv{v}{}{\nu}\right)\right] \Delta t.
\end{equation}
\noindent By inserting Eqs.~\eqref{eq: variation of normal} and \eqref{eq: variation of christoffel} into \eqqref{eq: variations of reference strains}, we obtain the final expression for the variation of the rate of the change of curvature: 
\begin{equation}
\label{eq: variation of curvature change}
\begin{aligned}
	\delta \imd{\kappa}{}{\alpha\beta} &= \iv{g}{}{3} \cdotp \left( \delta \iv{v}{}{\alpha,\beta} - \ii{\Gamma}{\mu}{\alpha\beta} \delta \iv{v}{}{\mu}\right) \\
	&\;\;\;\; - \left\{ \left( \iv{g}{}{3} \cdotp \delta \iv{v}{}{\mu} \right) \iv{g}{\mu}{} \cdotp \left( \iv{v}{}{\alpha,\beta} - \ii{\Gamma}{\nu}{\alpha\beta} \iv{v}{}{\nu}\right) + \left[ \iv{g}{\nu}{} \cdotp \delta \iv{v}{}{\alpha,\beta} - \ii{\Gamma}{\mu}{\alpha\beta} \left(\iv{g}{\nu}{} \cdotp \delta \iv{v}{}{\mu}\right)\right] \left( \iv{g}{}{3} \cdotp \iv{v}{}{\nu}\right)\right\} \Delta t.
\end{aligned}
\end{equation}

\subsection{Discrete equation of motion}

Let us find the matrix $\iv{B}{}{L}$ which relates the reference strain rates of the midsurface with the velocities of control points:
\begin{equation}
\label{eq: e=BL q}
	\ve{e} = \iv{B}{}{L} \ivmd{q}{}{}.
\end{equation}
\noindent The vector of reference strain rates, using Eqs.~\eqref{eq:def: strain rate wrt Cartesian} and \eqref{eq: tensor of curvature change -rate}, can be represented as:
\begin{equation}
\label{eq: vector of reference strains matrix form}
\begin{gathered}
	\ve{e} = \ve{H} \ve{w} = \sum_{m=1}^{3} \iv{H}{}{m} \iv{w}{m}{}, \quad \ve{H} = 
	\begin{bmatrix}
	\iv{H}{}{1} & \iv{H}{}{2} & \iv{H}{}{3} 
	\end{bmatrix},\\
	\iv{H}{}{m} = 
	\begin{bmatrix}
		\ii{x}{}{m,1} & 0 & 0 & 0 & 0 \\
		0 & \ii{x}{}{m,2} & 0 & 0 & 0 \\
		\ii{x}{}{m,2} & \ii{x}{}{m,1} & 0 & 0 & 0 \\
		-\ii{\Gamma}{1}{11} \ii{x}{}{m,3} & -\ii{\Gamma}{2}{11} \ii{x}{}{m,3} & \ii{x}{}{m,3} & 0 & 0 \\
		-\ii{\Gamma}{1}{22} \ii{x}{}{m,3} & -\ii{\Gamma}{2}{22} \ii{x}{}{m,3} & 0 & \ii{x}{}{m,3} & 0 \\
		-2\ii{\Gamma}{1}{12} \ii{x}{}{m,3} & -2 \ii{\Gamma}{2}{12} \ii{x}{}{m,3} & 0 & 0 & 2\ii{x}{}{m,3}
	\end{bmatrix},  \quad 
	\iv{w}{m}{} = 
	\begin{bmatrix}
		\ii{v}{m}{,1} \\
		\ii{v}{m}{,2} \\
		\ii{v}{m}{,11} \\
		\ii{v}{m}{,22} \\
		\ii{v}{m}{,12} \\
	\end{bmatrix},
	\end{gathered}
\end{equation}
\noindent where the vector $\ve{w}$ is:
\begin{equation}
\label{eq:w definition}
	\ve{w} = \ve{B} \ivmd{q}{}{}, \;  \ve{w} = 
	\begin{bmatrix}
	\iv{w}{1}{} \\
	\iv{w}{2}{} \\
	\iv{w}{3}{} 
	\end{bmatrix}, \quad \ve{B} = 
	\begin{bmatrix}
	\iv{B}{1}{11} & ... & \iv{B}{1}{IJ} & ... & \iv{B}{1}{NM} \\
	\iv{B}{2}{11} & ... & \iv{B}{2}{IJ} & ... & \iv{B}{2}{NM} \\
	\iv{B}{3}{11} & ... & \iv{B}{3}{IJ} & ... & \iv{B}{3}{NM} 
	\end{bmatrix}.
\end{equation}
\noindent The submatrices $\iv{B}{m}{IJ}$ for an arbitrary control point $IJ$ consists of the derivatives of basis functions:
\begin{equation}
\label{eq: submatrices BIJ}
	\iv{B}{1}{IJ} = 
	\begin{bmatrix}
	\ii{R}{}{IJ,1} & 0 & 0 \\
	\ii{R}{}{IJ,2} & 0 & 0 \\
	\ii{R}{}{IJ,11} & 0 & 0 \\
	\ii{R}{}{IJ,22} & 0 & 0 \\
	\ii{R}{}{IJ,12} & 0 & 0 \\
	\end{bmatrix}, \quad
	\iv{B}{2}{IJ} = 
	\begin{bmatrix}
	0 & \ii{R}{}{IJ,1} & 0 \\
	0 & \ii{R}{}{IJ,2} & 0 \\
	0 & \ii{R}{}{IJ,11} & 0 \\
	0 & \ii{R}{}{IJ,22} & 0 \\
	0 & \ii{R}{}{IJ,12} & 0 \\
	\end{bmatrix}, \quad
	\iv{B}{3}{IJ} = 
	\begin{bmatrix}
	0 & 0 & \ii{R}{}{IJ,1} \\
	0 & 0 & \ii{R}{}{IJ,2} \\
	0 & 0 & \ii{R}{}{IJ,11} \\
	0 & 0 & \ii{R}{}{IJ,22} \\
	0 & 0 & \ii{R}{}{IJ,12} \\
	\end{bmatrix}.
\end{equation}
\noindent The matrix $\iv{B}{}{L}$ now follows as:
\begin{equation}
\label{eq: e=Hw, BL=HB}
	\ve{e} = \ve{H} \ve{w} = \ve{H} \ve{B} \ivmd{q}{}{} = \iv{B}{}{L} \ivmd{q}{}{}, \quad \iv{B}{}{L} = \ve{H} \ve{B}.
\end{equation}

Furthermore the requirement, in order to rigorously define the geometric stiffness term, is to find the explicit matrix form of the virtual power generated by the known stress and the variation of strain rate with respect to the metric, Eqs.~\eqref{matrix form of linearized virtual power} and \eqref{eq: variations of reference strains}:
\begin{equation}
\label{eq: part of vp generated by known stress and variation of strain rate}
\begin{aligned}
	\int_{a}^{} \trans{f} \delta \iv{B}{}{L} \ivmd{q}{}{} \dd{a} \Delta t &= \int_{a}^{} \sum_{m=1}^{3} \sum_{n=1}^{3}  \trans{$(\iv{w}{m}{})$} \iv{G}{n}{m} \delta \iv{w}{}{m} \dd{a} \Delta t \\
	&= \frac{1}{2} \int_{a}^{} \ii{N}{\alpha\beta}{} \left( \ii{v}{m}{,\beta} \delta \ii{v}{}{m,\alpha} + \ii{v}{m}{,\alpha} \delta \ii{v}{}{m,\beta} \right) \dd{a} \Delta t \\
	&\;\;\;\; + \int_{a}^{} \ii{M}{\alpha\beta}{} \left[ \ii{v}{k}{,\nu} \left( \ii{\Gamma}{\nu}{\alpha\beta} \ii{x}{n}{,3} \ii{x}{,\mu}{k} + \ii{\Gamma}{\mu}{\alpha\beta} \ii{x}{n,\nu}{} \ii{x}{}{k,3}\right) \delta \ii{v}{}{n,\mu} \right. \\ 
	&\left. \;\;\;\; - \ii{v}{k}{,\alpha\beta} \ii{x}{n}{,3} \ii{x}{,\mu}{k} \delta \ii{v}{}{n,\mu} - \ii{v}{k}{,\nu} \ii{x}{n,\nu}{} \ii{x}{}{k,3} \delta \ii{v}{}{n,\alpha\beta} \right] \dd{a} \Delta t,
	\end{aligned}
\end{equation}
\noindent where the matrix $\iv{G}{n}{m}$ is: 
\begin{equation}
\label{eq:Gnm def}
	\iv{G}{n}{m} = 
	\begin{bmatrix}
	\left(\ii{G}{n}{m}\right)_{11} & \left(\ii{G}{n}{m}\right)_{12} & -\ii{M}{11}{} \ii{x}{n,1}{} \ii{x}{}{m,3} & - \ii{M}{22}{} \ii{x}{n,1}{} \ii{x}{}{m,3} & - 2 \ii{M}{12}{} \ii{x}{n,1}{} \ii{x}{}{m,3} \\
	\left(\ii{G}{n}{m}\right)_{21} & \left(\ii{G}{n}{m}\right)_{22} & -\ii{M}{11}{} \ii{x}{n,2}{} \ii{x}{}{m,3} & - \ii{M}{22}{} \ii{x}{n,2}{} \ii{x}{}{m,3} & - 2 \ii{M}{12}{} \ii{x}{n,2}{} \ii{x}{}{m,3} \\
	-\ii{M}{11}{} \ii{x}{n}{,3} \ii{x}{,1}{m} & -\ii{M}{11}{} \ii{x}{n}{,3} \ii{x}{,2}{m} & 0 & 0 & 0 \\
	-\ii{M}{22}{} \ii{x}{n}{,3} \ii{x}{,1}{m} & -\ii{M}{22}{} \ii{x}{n}{,3} \ii{x}{,2}{m} & 0 & 0 & 0 \\
	- 2 \ii{M}{12}{} \ii{x}{n}{,3} \ii{x}{,1}{m} & - 2 \ii{M}{12}{} \ii{x}{n}{,3} \ii{x}{,2}{m} & 0 & 0 & 0
	\end{bmatrix},
\end{equation}
\noindent with:
\begin{equation}
\label{eq:Gnm elements def}
\begin{aligned}
	\left(\ii{G}{n}{m}\right)_{11} &= \ii{\Gamma}{1}{M} \left(\ii{x}{n}{,3} \ii{x}{,1}{m} + \ii{x}{n,1}{} \ii{x}{}{m,3} \right) + \ii{\delta}{n}{m} \ii{N}{11}{},\\
	\left(\ii{G}{n}{m}\right)_{12} &= \ii{\Gamma}{1}{M} \ii{x}{n}{,3} \ii{x}{,2}{m} + \ii{\Gamma}{2}{M} \ii{x}{n,1}{} \ii{x}{}{m,3} + \ii{\delta}{n}{m} \ii{N}{12}{},\\
	\left(\ii{G}{n}{m}\right)_{21} &= \ii{\Gamma}{2}{M} \ii{x}{n}{,3} \ii{x}{,1}{m} + \ii{\Gamma}{1}{M} \ii{x}{n,2}{} \ii{x}{}{m,3} + \ii{\delta}{n}{m} \ii{N}{21}{},\\
	\left(\ii{G}{n}{m}\right)_{22} &= \ii{\Gamma}{2}{M} \left(\ii{x}{n}{,3} \ii{x}{,2}{m} + \ii{x}{n,2}{} \ii{x}{}{m,3} \right) + \ii{\delta}{n}{m} \ii{N}{22}{},\\
	\ii{\Gamma}{\nu}{M} &= \ii{\Gamma}{\nu}{11} \ii{M}{11}{} + 2 \ii{\Gamma}{\nu}{12} \ii{M}{12}{} + \ii{\Gamma}{\nu}{22} \ii{M}{22}{}.
\end{aligned}
\end{equation}
\noindent By introducing the total matrix of the generalized section forces:
\begin{equation}
\label{eq: G matrix def}
	\ve{G} = 
	\begin{bmatrix}
	\iv{G}{1}{1} & \iv{G}{2}{1} & \iv{G}{3}{1} \\
	\iv{G}{1}{2} & \iv{G}{2}{2} & \iv{G}{3}{2} \\
	\iv{G}{1}{3} & \iv{G}{2}{3} & \iv{G}{3}{3}
	\end{bmatrix},
\end{equation}
\noindent expression \eqref{eq: part of vp generated by known stress and variation of strain rate} can be written as:
\begin{equation}
\label{eq:transf geometric term of VP to bilinear form}
	\int_{a}^{} \trans{f} \delta \iv{B}{}{L} \ivmd{q}{}{} \dd{a} \Delta t = \int_{a}^{} \trans{w} \ve{G} \delta \ve{w} \dd{a} \Delta t.
\end{equation}
\noindent A careful inspection of Eqs.~\eqref{eq:Gnm def} and \eqref{eq:Gnm elements def} reveals the fact that the matrix of generalized section forces, $\ve{G}$, is symmetric. \textcolor{black}{This form was obtained in \cite{2017radenkovic} first.} Now, the terms in the equation of the virtual power, \eqref{matrix form of linearized virtual power}, reduce to:
\begin{equation}
\label{eq: terms of VP rewritten}
\begin{aligned}
	\int_{a}^{} \transmd{f} \delta \ve{e} \dd{a} &\approx \int_{a}^{} \transmd{q} \trans{$\iv{B}{}{L}$} \ve{D} \iv{B}{}{L} \delta \ivmd{q}{}{} \dd{a}, \quad \left( \delta \ve{e} \approx \frac{1}{2} \left(\iv{g}{}{\alpha} \cdotp \delta \iv{v}{}{\beta} + \iv{g}{}{\beta} \cdotp \delta \iv{v}{}{\alpha}\right) = \iv{B}{}{L} \delta \ivmd{q}{}{} \right)\\
	\int_{a}^{} \trans{f} \delta \ve{e} \dd{a} &= \int_{a}^{} \trans{f} \left( \delta \iv{B}{}{L} \ivmd{q}{}{} + \iv{B}{}{L} \delta \ivmd{q}{}{}\right) \dd{a} = \int_{a}^{} \left( \transmd{q} \trans{B} \ve{G} \ve{B} \delta \ivmd{q}{}{} \Delta t + \trans{f} \iv{B}{}{L} \delta \ivmd{q}{}{} \right) \dd{a},
\end{aligned}
\end{equation}
\noindent and the equation of equilibrium becomes:
\begin{equation}
\label{eq: virtual equilibrium}
	\transmd{q} \int_{a}^{} \left( \trans{$\iv{B}{}{L}$} \ve{D} \iv{B}{}{L} + \trans{B} \ve{G} \ve{B} \right) \dd{a} \delta \ivmd{q}{}{} \Delta t = \int_{a}^{} \trans{p} \ve{N} \dd{a} \delta \ivmd{q}{}{} - \int_{a}^{} \trans{f} \iv{B}{}{L} \dd{a} \delta \ivmd{q}{}{}.
\end{equation}
\noindent \eqqref{eq: virtual equilibrium} reduces to the standard form:
\begin{equation}
\label{eq:standard form of equlibrium}
	\iv{K}{}{T} \Delta \ve{q} = \ve{Q} - \ve{F}, \quad \left( \Delta \ve{q} = \ivmd{q}{}{} \Delta t\right),
\end{equation}
\noindent where:
\begin{equation}
\label{eq: Kt}
	\iv{K}{}{T} = \int_{a}^{} \trans{$\iv{B}{}{L}$} \ve{D} \iv{B}{}{L} \dd{a} + \int_{a}^{} \trans{B} \ve{G} \ve{B} \dd{a}, 
\end{equation}
\noindent is the tangent stiffness matrix and:
\begin{equation}
\label{eq: Q and F}
	\ve{Q} = \int_{a}^{} \trans{N} \ve{p} \dd{a}, \quad \ve{F} = \int_{a}^{} \trans{$\iv{B}{}{L}$} \ve{f} \dd{a},
\end{equation}
\noindent are the vectors of the external and internal forces, respectively. The vector $\Delta \ve{q}$ in \eqqref{eq:standard form of equlibrium} contains increments of displacements of control points with respect to the previous configuration. Due to the approximations introduced with \eqqref{eq: terms of VP rewritten}, the solution of \eqqref{eq:standard form of equlibrium} does not satisfy the principle of virtual power and some additional numerical procedure is required, see Section 4.4.1.

Although the derivation of the geometric stiffness matrix represents a standard procedure for nonlinear shell formulations, \cite{2015kiendl}, the present derivation differs. In particular, the full KL shell metric is incorporated in \eqqref{eq: part of vp generated by known stress and variation of strain rate} and the symmetric nature of the geometric stiffness term is stressed by the elegant and compact form of Eqs.~\eqref{eq: G matrix def} and \eqref{eq: Kt}. This confirms that the formulation is correct and that the adopted force and strain quantities are energetically conjugated.

\subsubsection{Numerical solution procedure}

The incremental form of \eqqref{eq:standard form of equlibrium} emphasizes that the equilibrium path must be considered as a set of discrete equilibrium points $(\iv{q}{i}{}, \iv{Q}{i}{})$ which correspond to a set of fictitious moments in time. In order to find these points, some iterative procedure must be employed. Here, the arc-length method is primarily used for all examples. It is worth noting that the Newton-Raphson method is applicable for the examples that do not exhibit snap-through behavior. We compared the performance of the linearized arc-length method, the cylindrical arc-length method, and the modified Riks method, \cite{1981crisfield, 2008ritto-correa}. For the numerical tests which follow in the next section, all of these variants return similar results. The exception is the last example where only the modified Riks method, as implemented in Abaqus, found the complete structural response, \cite{2009smith}. Hence, this procedure is briefly presented. 

The main idea of the arc-length method is to introduce a load proportionality factor $\lambda^i$ which scales the value of total external load, $\iv{Q}{f}{}$, and allows it to increase or decrease during the deformation process:
\begin{equation}
\label{eq eq for numeric1}
	\iv{Q}{i}{} = \ii{\lambda}{i}{} \iv{Q}{f}{}.
\end{equation}
\noindent This approach results in an indeterminate system of $3NM$ equations with $3NM+1$ unknowns where the additional unknown is the load proportionality factor itself, \cite{1981crisfield}. We start from some known, converged, configuration defined in a point $(\iv{q}{i}{}, \ii{\lambda}{i}{})$ and look for the next equilibrium point $(\iv{q}{i+1}{}, \ii{\lambda}{i+1}{})$. The first step is to find a predictor solution. Therefore, the predictor tangential displacement vector, $\iiv{q}{0}{}{}{T}$, is calculated from the known tangent stiffness and the total external load:
\begin{equation}
\label{eq eq for numeric2}
	\iiv{K}{}{}{i}{T} \: \iiv{q}{0}{}{i+1}{T} = \iv{Q}{f}{}.
\end{equation}
\noindent The predictor load increment $\Delta \: \iii{\lambda}{0}{}{i+1}{}$ follows from the appropriate constraint equation, \cite{1981crisfield, 2009smith}:
\begin{equation}
\label{eq eq for numeric3}
	\Delta \: \iii{\lambda}{0}{}{i+1}{} = \pm \frac{\Delta \: \iii{l}{}{}{i+1}{}}{\sqrt{\iivn{q}{0}{}{i+1}{T} \: \cdotp \: \iivn{q}{0}{}{i+1}{T} +1}}.
\end{equation}
\noindent Here, $\Delta \iii{l}{}{}{i+1}{} $ is the arc-length value which is estimated for the first increment, and calculated for all the other increments as a function of desired number of iterations. For convenience, the displacement vectors are scaled by the maximum displacement component of the corresponding linear solution and designated with tilde, \cite{2009smith}. The sign of the predictor solution is calculated from the condition that the projection of the predictor tangential displacement onto the previously converged displacement increment must be positive, ensuring continuation along the equilibrium path. This condition can be written as:
\begin{equation}
\label{eq eq for numeric6}
	\Delta \: \iii{\lambda}{0}{}{i+1}{} \left( \iivn{q}{0}{}{i+1}{T} \Delta \: \iivn{q}{}{}{i}{} + \Delta \: \iii{\lambda}{}{}{i}{} \right) > 0.
\end{equation}
\noindent In general, the predictor solution does not satisfy the equilibrium and the arc-length iterations act as correctors. For the first iteration, the increments of displacement and load are initialized from the predictor solutions:
\begin{equation}
\label{eq eq for numeric7}
	\Delta \: \iii{\lambda}{1}{}{i+1}{} = \Delta \: \iii{\lambda}{0}{}{i+1}{}, \; \Delta \: \iiv{q}{1}{}{i+1}{} = \Delta \: \iii{\lambda}{0}{}{i+1}{} \: \iiv{q}{0}{}{i+1}{T}.
\end{equation}
\noindent Now, the current tangent stiffness, $\iiv{K}{j}{}{i+1}{T}$, and internal forces, $\iiv{F}{j}{}{i+1}{}$, are calculated and convergence is checked for each iteration $j=1,2,3...$. The convergence criteria can be defined with respect to the current unbalanced forces, 
$ \iii{\boldsymbol{\Psi}}{j}{}{i+1}{} = \iv{Q}{i}{} - \: \iii{\textbf{F}}{j}{}{i+1}{} $, or/and with respect to the previously calculated iterative displacement. If the defined criteria is not met, new tangential, $\iiv{q}{j}{}{i+1}{T}$, and residual, $\delta \: \iiv{q}{j}{}{i+1}{}$, displacements are calculated from:
\begin{equation}
\label{eq eq for numeric8}
\iiv{K}{j}{}{i+1}{T} \: \iiv{q}{j}{}{i+1}{T} = \iiv{Q}{}{}{f}{}, \; \; \iiv{K}{j}{}{i+1}{T} \: \delta \: \iiv{q}{j}{}{i+1}{} = \: \iii{\boldsymbol{\Psi}}{j}{}{i+1}{}.
\end{equation}
\noindent The iterative load factor, $\delta \lambda$, follows from the condition that a new potential equilibrium point is orthogonal to the tangential displacement in the arc-length solution space:
\begin{equation}
\label{eq eq for numeric12}
\delta \lambda = -\frac{\delta \: \iivn{q}{j}{}{i+1}{} \: \cdotp \: \iivn{q}{0}{}{i+1}{T}}{\iivn{q}{j}{}{i+1}{T} \: \cdotp \: \iivn{q}{0}{}{i+1}{T} + 1}.
\end{equation}
\noindent Finally, the new solution is:
\begin{equation}
\label{eq eq for numeric11}
\begin{aligned}
	\iiv{q}{j+1}{}{i+1}{} &= \iiv{q}{}{}{i}{} + \Delta \: \iiv{q}{j}{}{i+1}{} + \delta \: \iiv{q}{j}{}{i+1}{} + \delta \lambda \: \iiv{q}{j}{}{i+1}{T}, \\
	\iii{\lambda}{j+1}{}{i+1}{} &= \iii{\lambda}{}{}{i}{} + \Delta \: \iii{\lambda}{j}{}{i+1}{} + \delta \lambda,
\end{aligned}
\end{equation}
\noindent and the procedure is repeated until the convergence criteria is fulfilled, \cite{2009smith}.

\section{Numerical examples}

There are four main objectives of the present numerical analysis: i) verification and validation of the developed nonlinear formulation and its implementation, ii) examination of the influence of full constitutive relation on the structural response, iii) comparison of the linear and nonlinear relation between reference and equidistant strains, iv) assessment of the influence of the curviness. 

To examine the influence of the different constitutive relations, we introduce four models $D$:
\begin{itemize}
	\item $\ii{D}{a}{}$ is the full analytical constitutive model. It is given in an expanded form in Appendix B.
	\item $\ii{D}{0}{}$ decouples the bending and membrane actions (e.g., \cite{2009kiendl}).
	\item $\ii{D}{1}{}$ employs a membrane-bending coupling based on the first-order Taylor approximation of the determinant of the reciprocal shift tensor $1/\ii{g}{}{0}$ as suggested in \cite{1968baker}, \eqref{C3}.
	\item $\ii{D}{2}{}$ is like $\ii{D}{1}{}$ but introduces an additional effect due to curvature (similar to the Flugge-Lure-Byrne-Truesdell equations, \cite{1973naghdi}), \eqref{C3}.
\end{itemize}

In order to examine the influence of the constitutive relation on the structural response, reference strains are observed and compared for different models. For $\ii{D}{a}{}$ model, the influence of the relation between the reference and equidistant strains is also examined by evaluating the strains on the outer surface of shell using the exact \eqref{eq: components of eq strain rate final} and the linearized relation \eqref{eq:linear form of eq strain}. Furthermore, the stresses are calculated by \eqqref{eq:constitutive relation at equidistant surface}. \textcolor{black}{Since the equilibrium paths of stresses and strains for nonlinear shells are rarely found in the literature, these quantities are here compared with the results from Abaqus simulations, which utilize the standard assumption of the linear strain distribution, \cite{2009smith}.} Since the time is a fictitious quantity in the present static analysis, strain and stress rates are equal to strains and stresses, respectively.

\textcolor{black}{We have analyzed each nonlinear example with various different meshes. For the sake of brevity, only a few results are reported here. Moreover, the solutions have singularities due to point forces and no optimal higher-order convergence can be achieved.}

Due to the $\ii{C}{1}{}$ interelement continuity, reference membrane strains and $\ii{\kappa}{}{12}$ are continuous at the interelement boundaries, Eqs.~\eqref{eq:def: strain rate wrt Cartesian} and \eqref{eq: tensor of curvature change -rate}. However, the curvature changes $\ii{\kappa}{}{\alpha\alpha}$ require at least $\ii{C}{2}{}$ continuous meshes and they are linearly averaged here for $\ii{C}{1}{}$ meshes. Generally, the post-processing of results in IGA is an ongoing topic, \cite{2017stahl}. The models in Abaqus are meshed with S4R shell elements and element nodal results are again linearly averaged. This general-purpose four-node quadrilateral element with reduced integration is readily utilized for the benchmarks of shell formulations, \cite{2004sze}. First, a classical linear shell example is investigated. The other examples address nonlinear responses of the different shell constitutive models.

The boundary conditions are imposed in a well-known manner where rotations are treated with special care since they are not utilized as DOFs, \cite{2009kiendl}. \textcolor{black}{The symmetry boundary conditions are enforced by constraining the displacements of control points that influences rotation.} Numerical integration is performed on an element level using the Gauss quadrature with $p \times p$ integration points. There are more efficient quadrature schemes for IGA \cite{2012auricchio}, but they are not considered in this paper. Furthermore, locking effects are not studied here.

All the results are presented with respect to the load proportionality factor (LPF), rather than the load intensity itself.

\subsection{Linear analysis of pinched cylinder}

\begin{figure}
	\includegraphics[width=\linewidth]{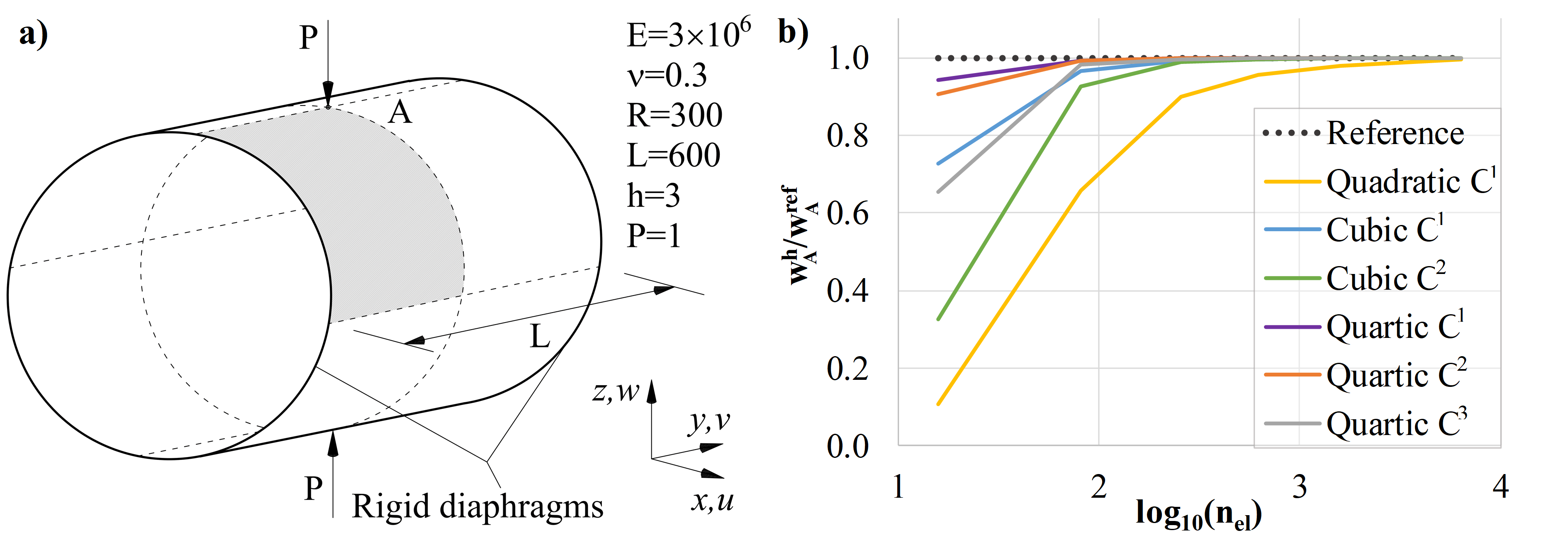}
	\caption{Linear analysis of the pinched cylinder. a) Geometry, material characteristics, and applied load. b) Convergence of the radial displacement at point A for different NURBS orders and interelement continuities.}
	\label{fig:pinched cylinder disposition and convergence}
\end{figure}

The first example deals with the linear analysis of a pinched cylinder with rigid diaphragms at its ends as shown in \fref{fig:pinched cylinder disposition and convergence}a. The deflection at point $A$ where the force is applied commonly serves as a benchmark test for shell elements. Since the structure has three planes of symmetry, only one eighth of the shell is considered. The analytical solution based on the Fourier series approximation exists, but it involves some numerical issues, as discussed in [12]. The reference value often used in the literature is $1.82488\times10^{-5}$ and it stems from the approximation of the solution with $80\times80$ Fourier terms. Here, the reference value of $1.827158\times10^{-5}$ obtained with $8192 \times 8192$ Fourier terms is adopted, \cite{2017duong}. The analysis results for different IGA meshes are shown in \fref{fig:pinched cylinder disposition and convergence}b where the expected convergence behavior can be observed. 

\begin{figure}
\centering
	\includegraphics[width=8cm]{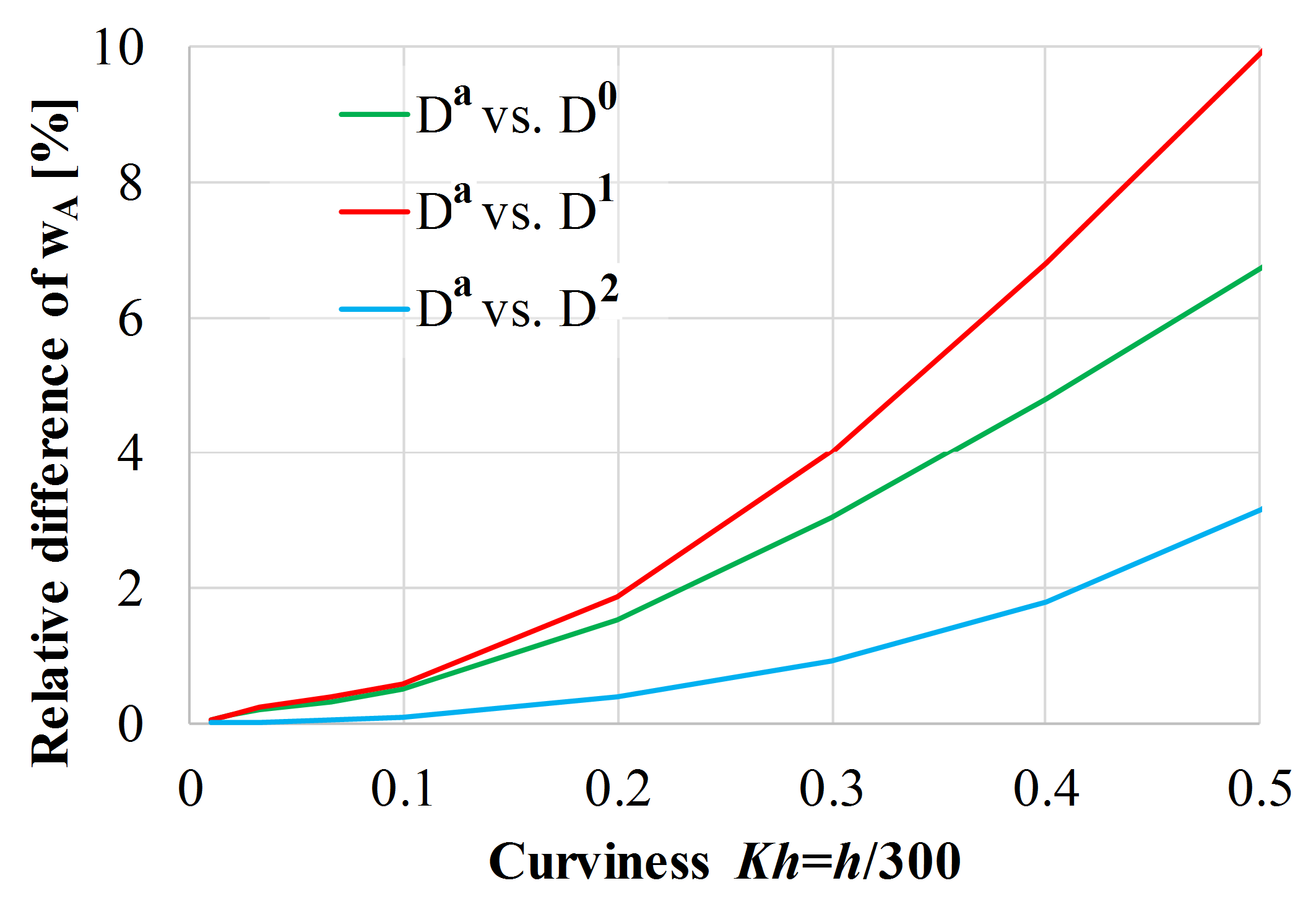}
	\caption{The relative difference of displacement for different constitutive material models with respect to the full analytical one.}
	\label{fig:pinched cylinder different constitutive models}
\end{figure}

\textcolor{black}{Furthermore, the influence of the different constitutive models on the deflection at point $A$ is examined. The midsurface is represented by $36\times36$ cubic $C^1$ elements and the curviness $Kh$ is varied. Since the geometry of midsurface is kept constant, the curviness becomes a function of thickness only, $Kh=h/300$.} In \fref{fig:pinched cylinder different constitutive models}, the obtained values of displacement component $w_A$ for the three reduced constitutive models are related to the one computed by the full analytical $D^a$ model. As expected, these relative differences increase with the curviness and the results of the $D^2$ model are the closest to the $D^a$ model. Note that for strongly curved shells, the $D^0$ model is more accurate than the $D^1$ model, which implies that the simplified coupling of membrane and bending actions introduced in the $D^1$ model does not improve the accuracy in this example. 

\subsection{Snap behavior of a shallow shell}

\begin{figure}
	\includegraphics[width=\linewidth]{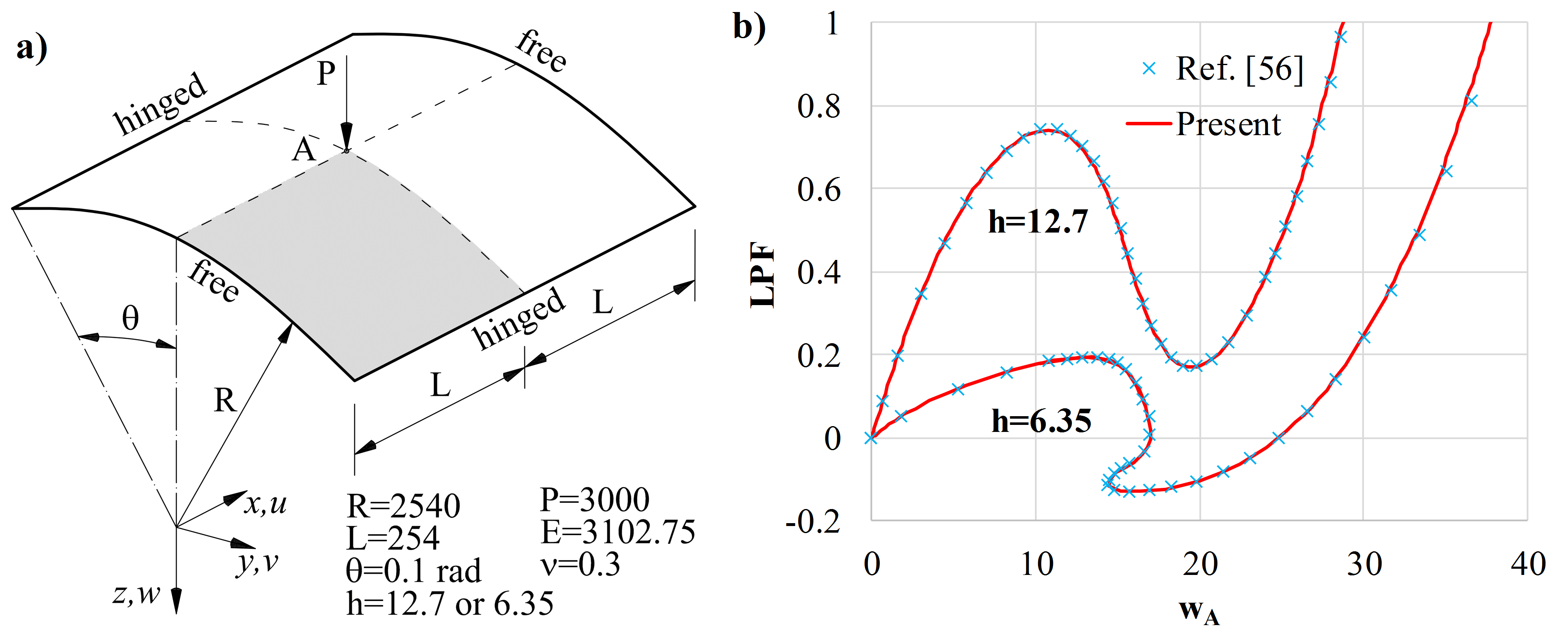}
	\caption{Shallow shell: a) Geometry, material characteristics, and applied load. b) Comparison of results obtained with the present formulation and those from \cite{2004sze}.}
	\label{fig:shallow shell}
\end{figure}

The second problem deals with the shallow shell illustrated in \fref{fig:shallow shell}a. This is a standard example which exhibits both snap-through and snap-back behavior. It serves as a rigorous test for algorithms that solve complex nonlinear behavior of shells. Several variants of this example exist, \cite{1981crisfield}, and we are referring to the ones discussed in \cite{2004sze}. \textcolor{black}{By enforcing the symmetry boundary conditions, only one fourth of the structure is analyzed.} A mesh with $4\times4$ quartic elements with $C^3$ interelement continuity is employed and the results are compared with those from \cite{2004sze} in \fref{fig:shallow shell}b. Almost full agreement is observed for the two different values of thickness. The results for all four constitutive models are practically the same due to the small curviness of this shell at all configurations. Overall, these results confirm that the presented formulation can deal with complex equilibrium paths, but it should be noted that these shallow shells do not undergo large rotations.

The following efficiency comparison of the different constitutive models is performed on the shallow shell with thickness $t=12.7$. Three series of calculations are done for each constitutive model using the Intel Core i7-8700 CPU. The results are linearly averaged and shown in Table \ref{tab:efficiency}. In each test, the prescribed accuracy of a solution is reached after 21 load increments and a total of 52 arc-length iterations. Evidently, the full analytical constitutive model $D^a$ requires significantly more computational time than the reduced models. This is expected since the expressions for the $D^a$ model are considerably more involved, see Appendix B, Eqs.~\eqref{B1}, \eqref{B2} and \eqref{B3}. On the other hand, the computational time of all three reduced models is quite similar.

\begin{table}
	\caption{Computational times for different constitutive models of the shallow shell example. For all models, the solution procedure requires 21 load increments and 52 arc-length iterations.} 
	\centering 
	\begin{tabular}{c c c c c} 
		\hline  
		\rule{0pt}{2.5ex}  
		Model & $D^0$ & $D^1$ & $D^2$ & $D^a$ \\ 
		\hline 
		\rule{0pt}{2.5ex} 
		Time [s] & 28.70 & 29.54 & 30.08 & 98.36 \\ 
		\hline 
	\end{tabular}
	\label{tab:efficiency} 
\end{table}

\subsection{Pinched semi-cylinder}

\begin{figure}
	\includegraphics[width=\linewidth]{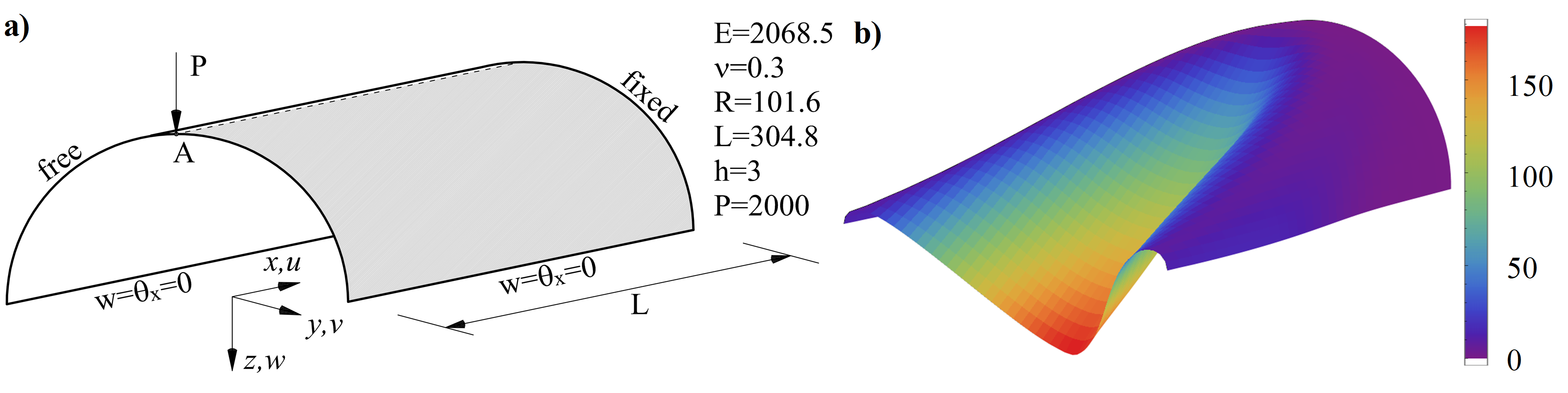}
	\caption{Pinched semi-cylinder. a) Geometry, material characteristics, and applied load. b) Deformed configuration for $\textrm{LPF}=1$. The colors show the total value of the displacement vector.}
	\label{fig:pinched semi cyl disposition and def config}
\end{figure}

\begin{figure}
	\includegraphics[width=\linewidth]{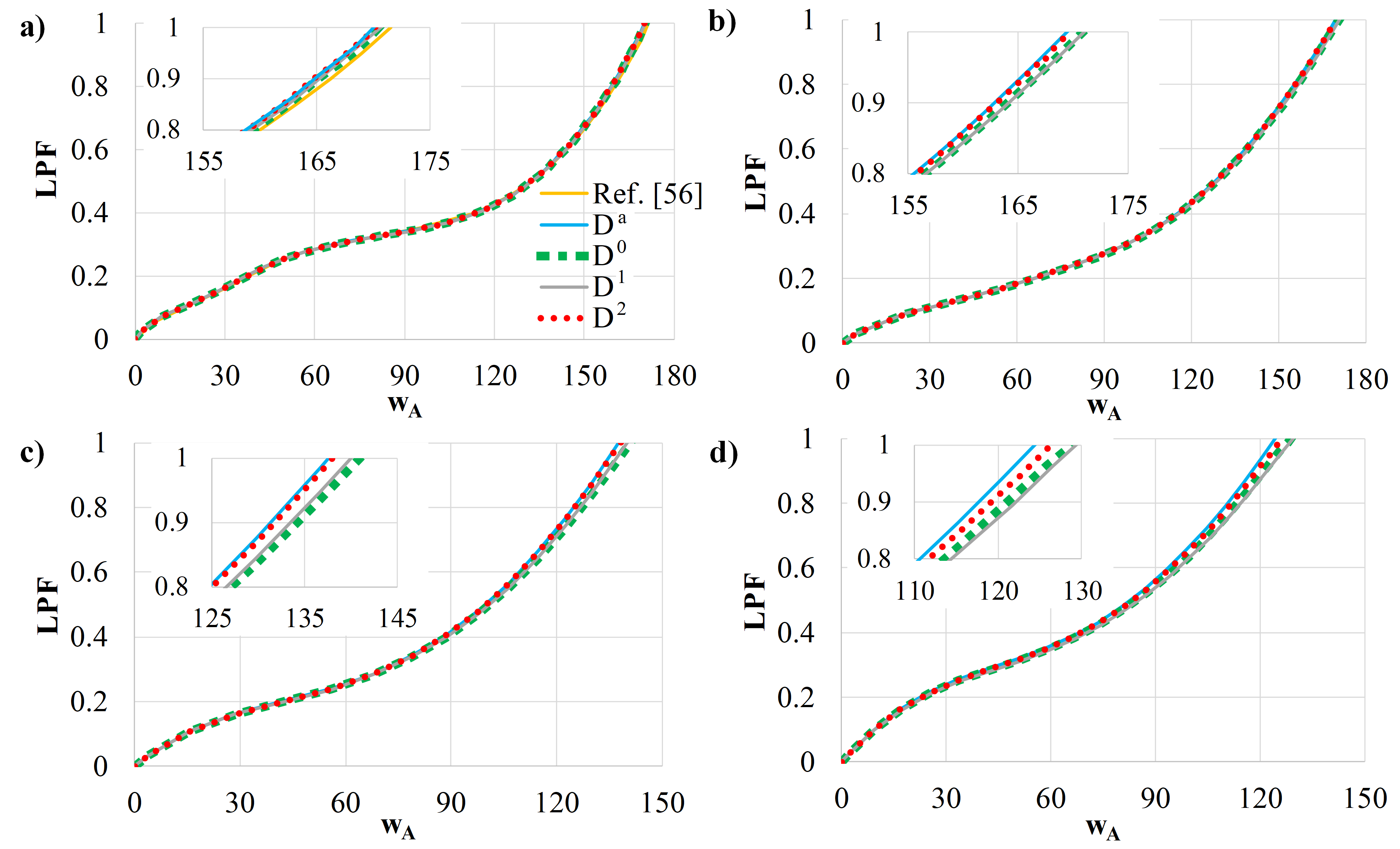}
	\caption{Pinched semi-cylinder. a) Comparison of present formulation with \cite{2004sze} for thickness $h$ and load $P$; Comparison of deflection for different constitutive models for: b) $2h$ and $7P$; c) $4h$ and $24P$; d) $8h$ and $95P$. Insets show zoomed parts for $\textrm{LPF}>0.8$.}
	\label{fig:pinched semi cyl comparison of displacement}
\end{figure}

The pinched semi-cylinder as in \fref{fig:pinched semi cyl disposition and def config}a is the well-documented nonlinear benchmark for thin shells which exhibit both large displacements and rotations. A mesh with $20\times20$ quadratic $C^1$ elements is utilized for the discretization of the one half of the symmetric structure and the deformed configuration is displayed in \fref{fig:pinched semi cyl disposition and def config}b. \textcolor{black}{A comparison of the deflection at point $A$ is shown in \fref{fig:pinched semi cyl comparison of displacement}a. All results are in good agreement, but small differences can be detected for LPF=1. Models $D^a$ and $D^2$ give virtually indistinguishable results, similarly as models $D^0$ and $D^1$. The largest relative difference between these models is close to 0.3 \%. Comparison with the reference solution given in \cite{2004sze} reveals that the largest relative difference is less than 0.9 \%.}

\begin{figure}
	\includegraphics[width=\linewidth]{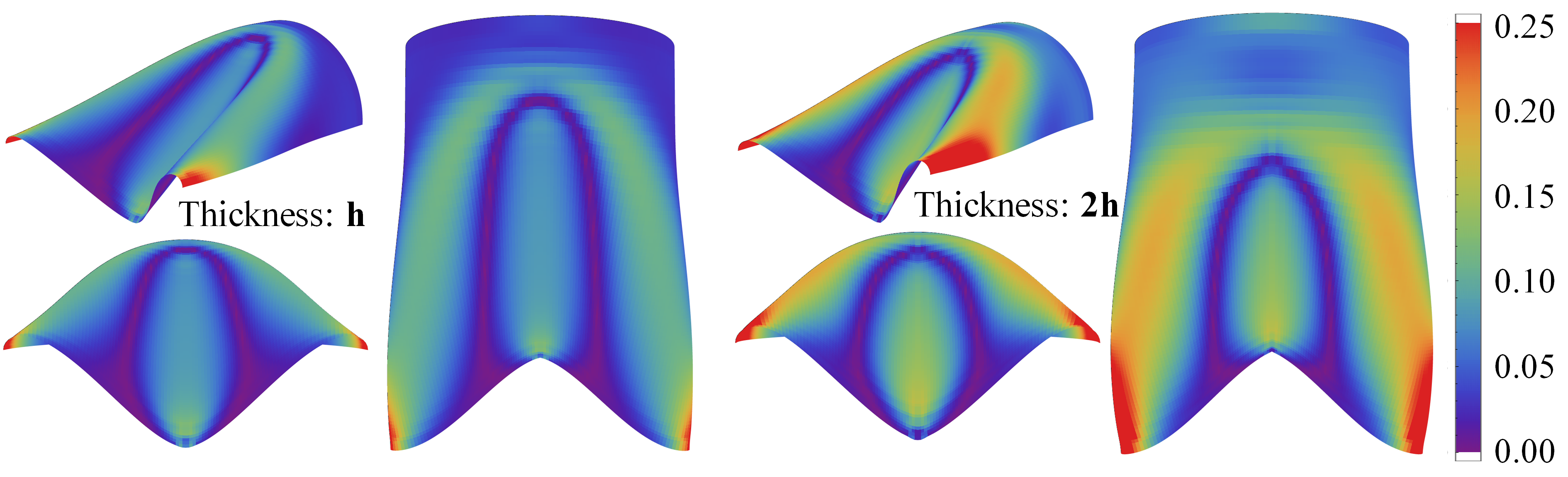}
	\caption{Pinched semi-cylinder. Distribution of curviness for $\textrm{LPF}=1$ from three points of view for a shell with two different values of thickness, i.e., $h$ and $2h$. For a clearer representation, a threshold of $Kh=0.25$ is applied.}
	\label{fig:pinched semi cyl curvniess}
\end{figure}

\textcolor{black}{In order to examine the influence of constitutive relation on this characteristic displacement component, the shell thickness and the load intensity are increased and three cases are considered: (i) $h/P=6/14000$, (ii) $h/P=12/48000$, and (iii) $h/P=24/190000$. The differences between the results obtained with the $D^a$ constitutive model and the reduced ones are studied. The corresponding results are given in \fref{fig:pinched semi cyl comparison of displacement}b, \fref{fig:pinched semi cyl comparison of displacement}c, and \fref{fig:pinched semi cyl comparison of displacement}d.} Evidently, all of the constitutive model types return virtually indistinguishable results for the first part of the equilibrium path, approximately $\textrm{LPF}<0.3$, regardless of the thickness. The discrepancies become visible for the other part of the equilibrium path, especially when $w_A>R$. This is due to the cumulative nonlinear effect and also because the curviness value $Kh$ becomes significant at some parts of the shell. It is emphasized that the curviness is measured at an arbitrary configuration. Note that the $D^2$ model starts to visibly deviate from the analytical $D^a$ model only for the shell with reference curviness $Kh = 24/101.6 \approx 0.236$ which corresponds to initially strongly curved shell. The distribution of curviness for two values of thickness and $\textrm{LPF}=1$ is shown in \fref{fig:pinched semi cyl curvniess} where the threshold value of $Kh=0.25$ is utilized in order to depict this characteristic more clearly. To be precise, all areas with $Kh \geq 0.25$ are colored with the same color. It is evident how some parts of the shell become strongly curved.

\subsection{Pullout of a cylinder}

\begin{figure}
	\includegraphics[width=\linewidth]{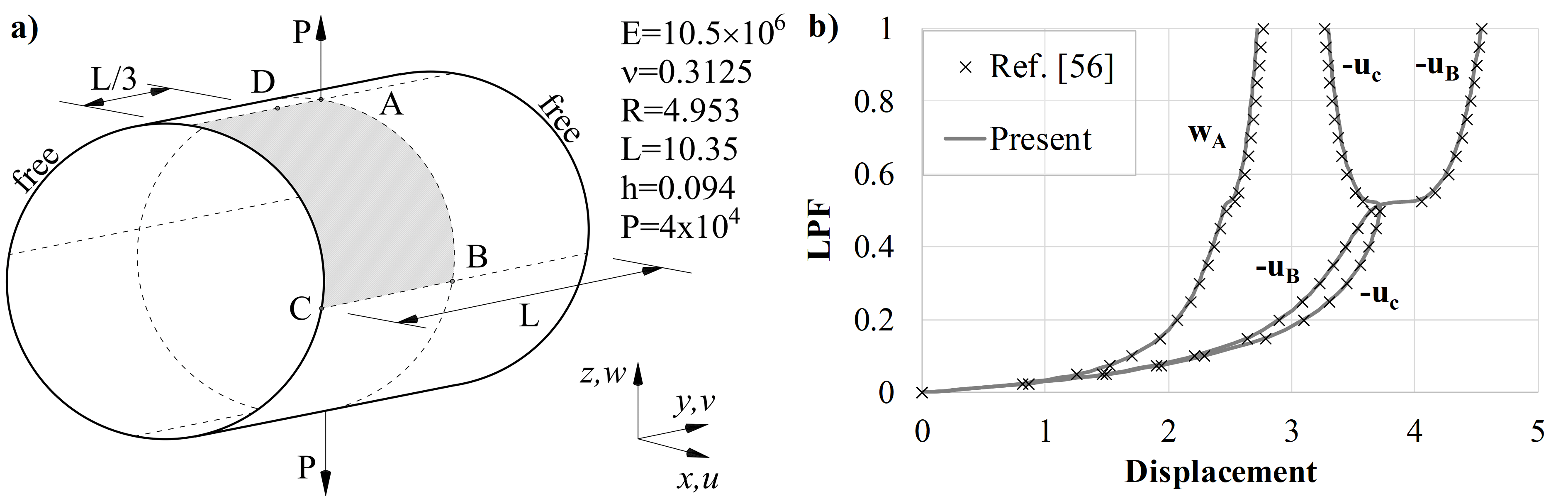}
	\caption{Pullout of a cylinder. a) Geometry, material characteristics, and applied load. b) Comparison of displacements between two present models and \cite{2004sze} for three displacement components. }
	\label{fig: pullout: disposition and displacement comparison}
\end{figure}

\begin{figure}
	\centering
	\includegraphics[width=9cm]{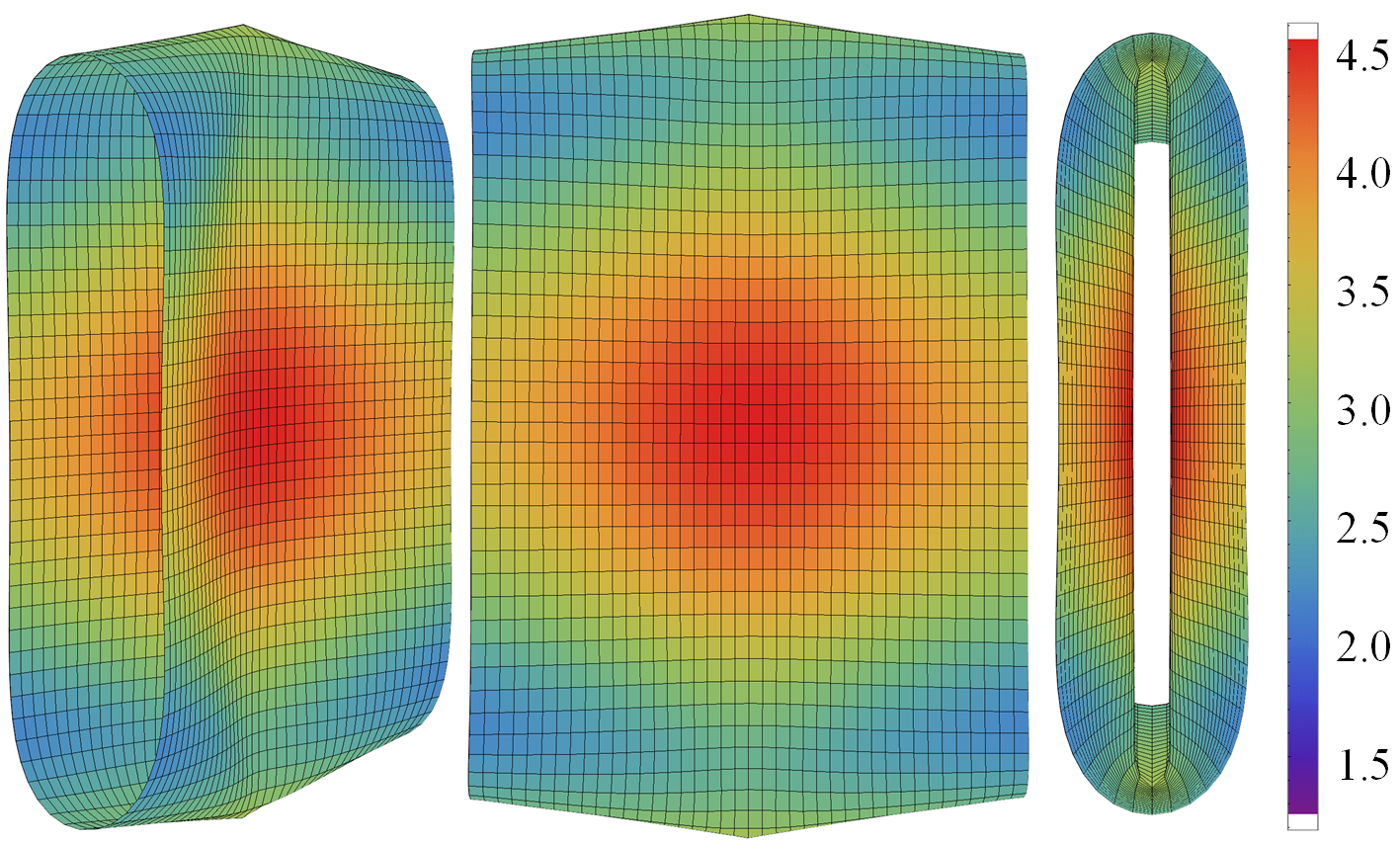}
	\caption{Pullout of a cylinder. Deformed configuration for $\textrm{LPF}=1$ from three different perspectives. The color-coding represents the total value of the displacement vector.}
	\label{fig:pullout def config}
\end{figure}

The geometric nonlinear analysis of a cylinder with free ends, as in \fref{fig: pullout: disposition and displacement comparison}a, is a well-established test for shell formulations. The shell is pulled out by the two opposing forces and three displacement components are compared, \cite{2004sze}. Due to the symmetry of the problem, only one eighth of the shell is discretized using a mesh with $60\times60$ cubic $C^2$ elements. The comparison with the reference solutions is displayed in \fref{fig: pullout: disposition and displacement comparison}b where the good agreement is apparent. \textcolor{black}{Since all constitutive models return virtually indistinguishable results, their comparison is omitted.} The deformed configuration from three different points of view is displayed in \fref{fig:pullout def config}.

\begin{figure}
	\includegraphics[width=\linewidth]{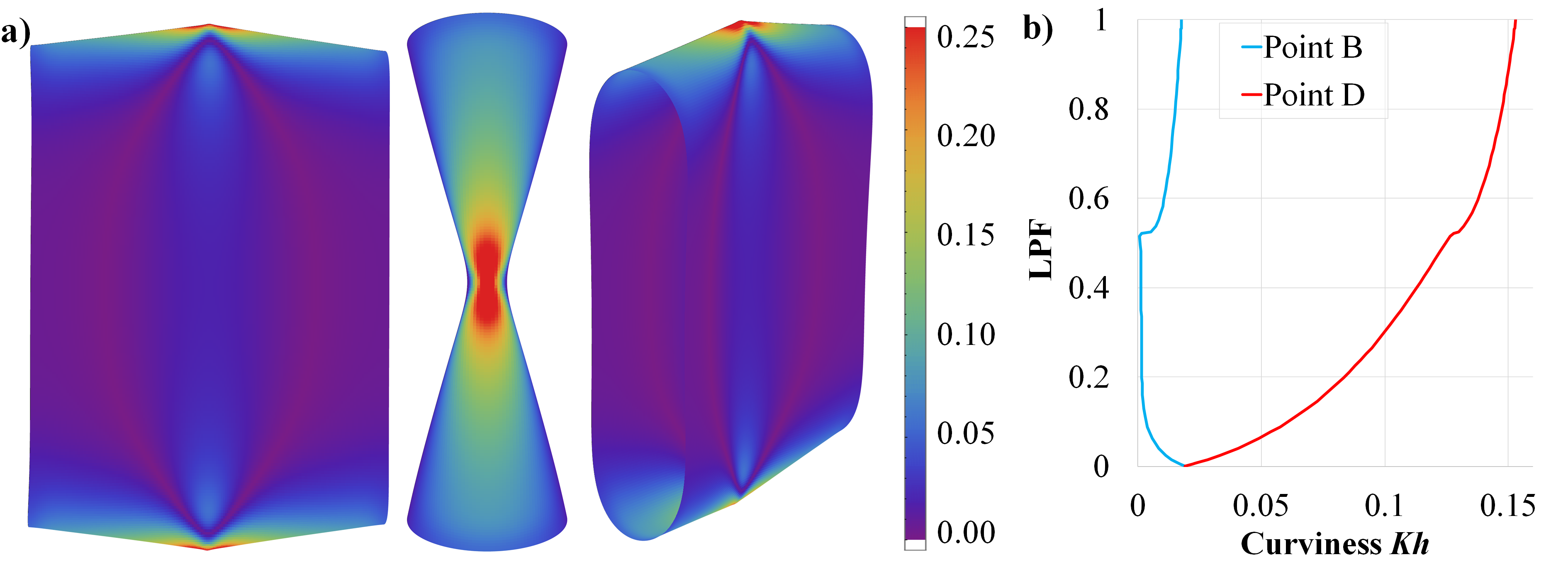}
	\caption{Pullout of a cylinder. a) The distribution of the curviness for $\textrm{LPF}=1$ from three points of view. The threshold value of $Kh=0.25$ is utilized. b) Curviness at points $B$ and $D$ vs. LPF.}
	\label{fig: pullout: curviness}
\end{figure}

The initial curviness of the observed shell is $Kh \approx 0.019$, but for $\textrm{LPF}=1$ the maximum curviness of $Kh=1.42$ develops at point $A$. The distribution of the curviness at the final configuration is given in \fref{fig: pullout: curviness}a using the threshold value of $Kh=0.25$. It is evident that some parts of this shell become strongly curved locally. The change of the curviness with respect to the load proportionality factor is shown in \fref{fig: pullout: curviness}b for two characteristic points.

\begin{figure}
	\includegraphics[width=\linewidth]{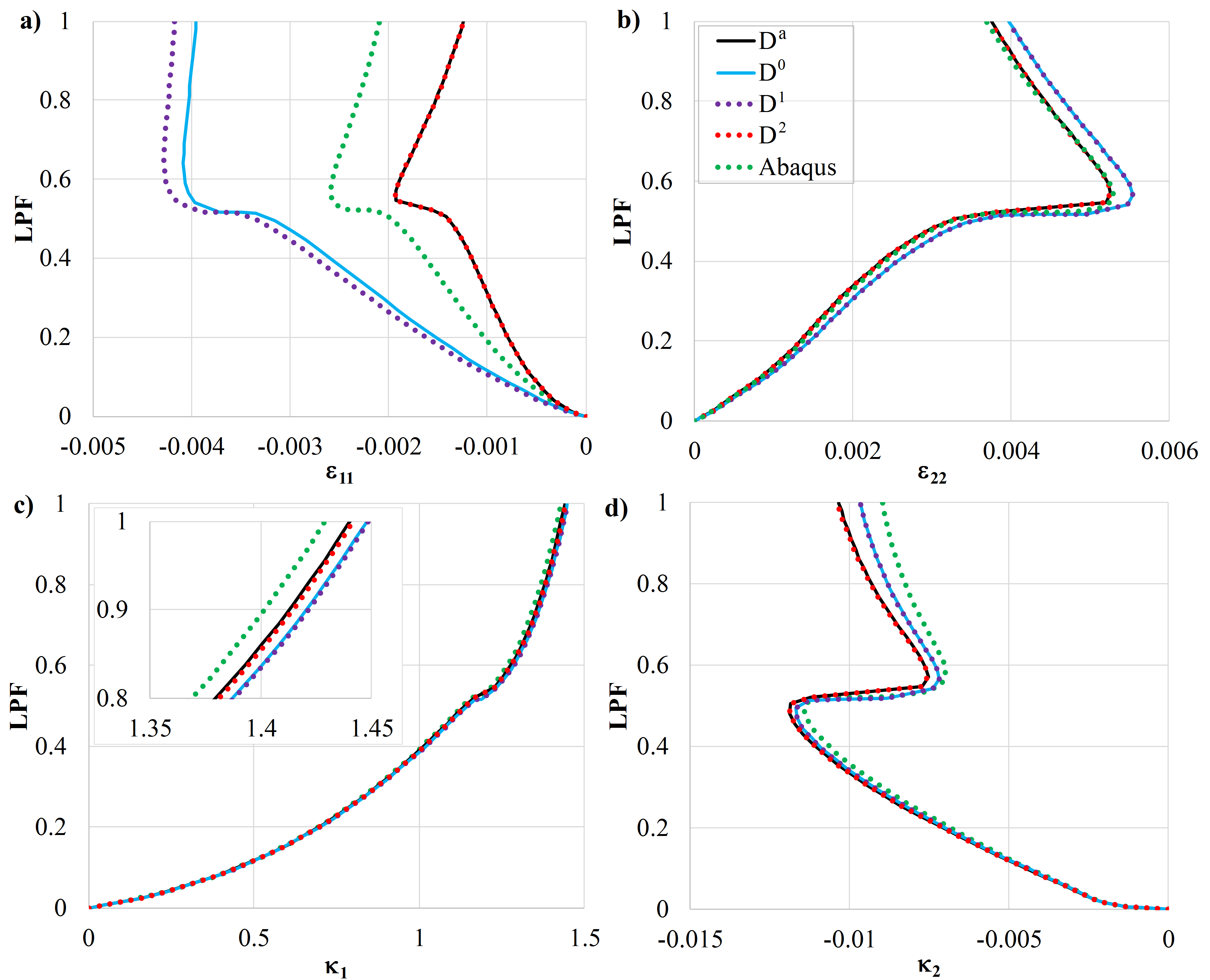}
	\caption{Pullout of a cylinder. Comparison of reference strains at point $D$ obtained with the four different constitutive models and Abaqus: a) $\ii{\epsilon}{}{11}$; b) $\ii{\epsilon}{}{22}$; c) $\ii{\kappa}{}{1}$; d) $\ii{\kappa}{}{2}$.}
	\label{fig: pullout: reference strains}
\end{figure}

In order to examine the influence of the curviness on the structural response, strains at point $D$ are considered in detail. This point exhibits strong local curviness, i.e., $Kh\approx0.154$, and it is far enough away from the applied force to avoid singularity issues. The four reference strains are calculated using the four different constitutive models. Furthermore, Abaqus results using a mesh with 15700 S4R elements are used for reference. All the results are compared in \fref{fig: pullout: reference strains} and interesting complex equilibrium paths are detected. 

Most importantly, the dominant reference strain is $\ii{\kappa}{}{1}$ and its equilibrium path appears to be practically identical for all the considered models, \fref{fig: pullout: reference strains}c. However, a close inspection of the zoomed part for $\textrm{LPF}>0.8$ reveals discrepancies, \fref{fig: pullout: reference strains}c. \textcolor{black}{The models $D^a$ and $D^2$ return practically identical results, while the same observation holds for the models $D^0$ and $D^1$. The relative difference of results for all reduced models and Abaqus, with respect to $D^a$ model, is less than 1 \% for LPF=1.} The other change of curvature, $\ii{\kappa}{}{2}$, is smaller for an order of magnitude and the results obtained with different models start to visibly deviate for $\textrm{LPF}>0.3$, \fref{fig: pullout: reference strains}d. Regarding the membrane reference strains, all models return different equilibrium paths for $\ii{\epsilon}{}{11}$. For $\ii{\epsilon}{}{22}$, the reduced models $D^0$ and $D^1$ differ from the full analytical $D^a$ solution. This is expected since the deformed curviness at point $D$ of this shell is dominated by the curvature component along the $\xi$ direction $\left(\ii{b}{1}{1} \gg \ii{b}{2}{2} \right)$, which reflects on these membrane strain components, see \eqqref{eq:def: strain rate wrt Cartesian}. The differences in comparison with Abaqus are evident and they are mainly attributed to the application of the more rigorous constitutive models within the present formulation. Additional influences are mesh density, continuity, and strain averaging at the interelement boundaries, which are standard issues in FEA and which are not of primary interest here.
 
\begin{figure}
	\includegraphics[width=\linewidth]{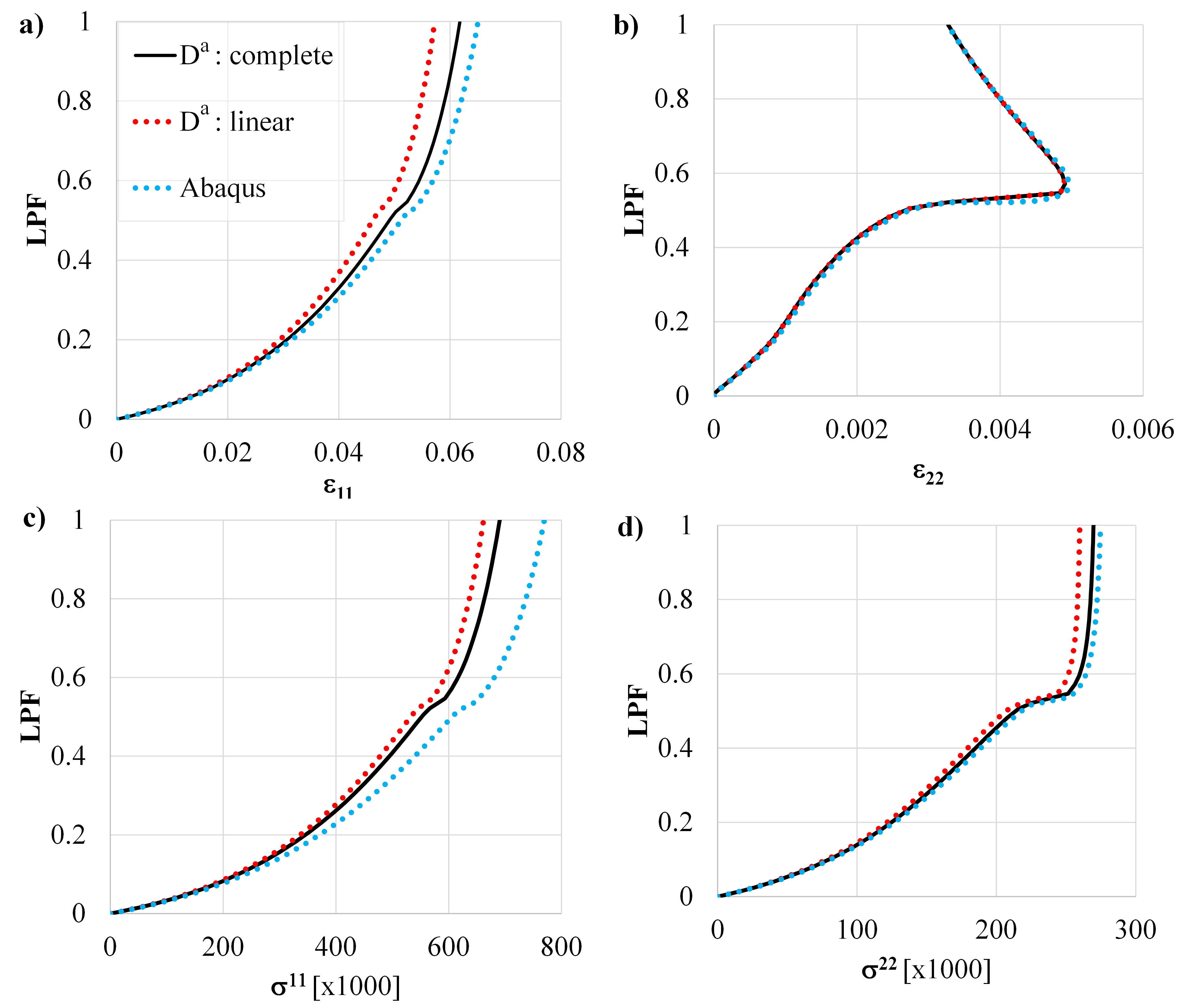}
	\caption{Pullout of a cylinder. Comparison of strains and stresses at the outer surface at point $D$ obtained with simplified linear and complete nonlinear relations between the reference and equidistant strains, and Abaqus: a) $\ii{\epsilon}{out}{11}$; b)  $\ii{\epsilon}{out}{22}$; c)  $\ii{\sigma}{11}{out}$; d)  $\ii{\sigma}{22}{out}$.}
	\label{fig:pullout: strains at outer point}
\end{figure}

\begin{figure}
	\includegraphics[width=\linewidth]{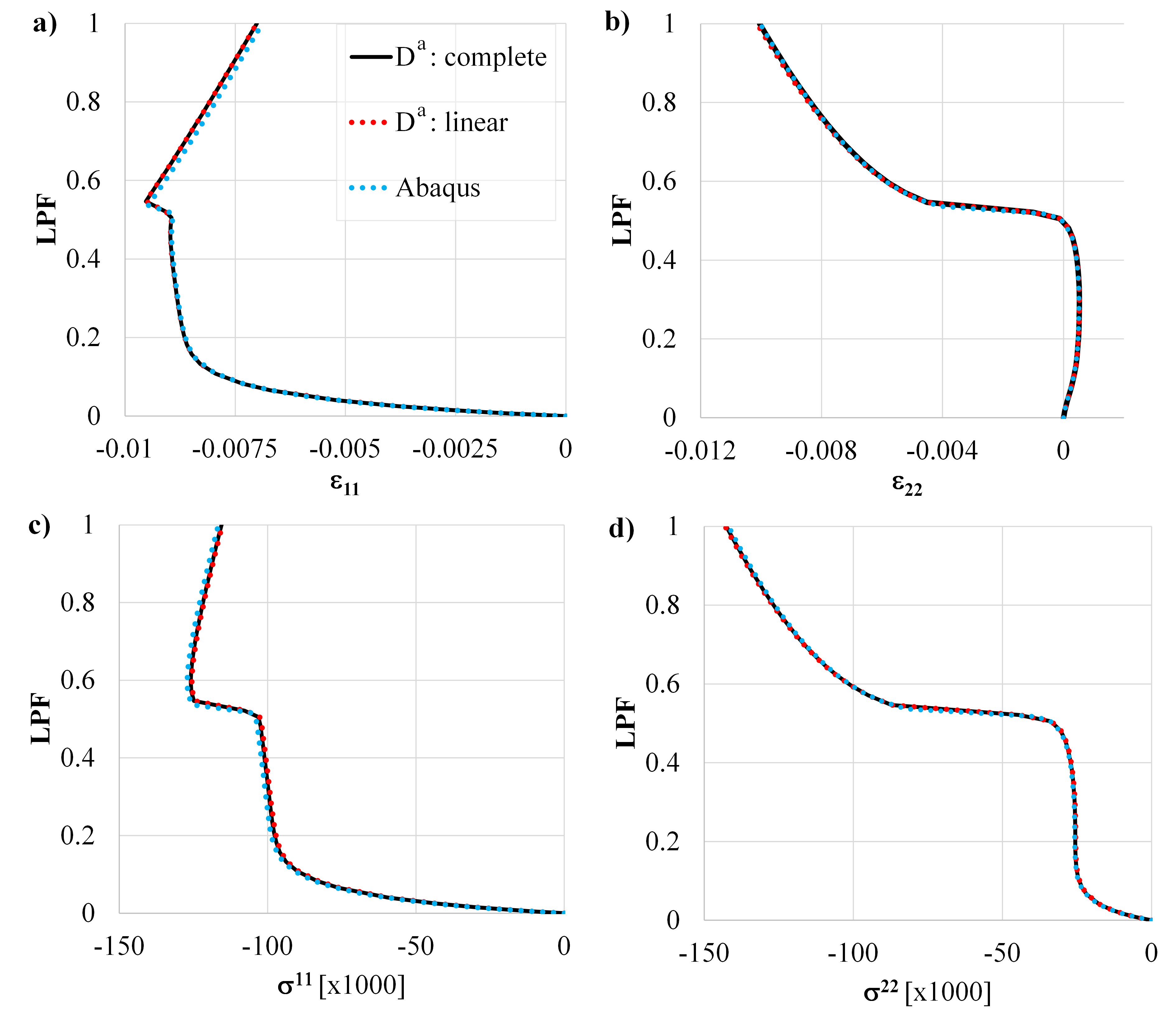}
	\caption{Pullout of a cylinder. Comparison of strains and stresses at the outer surface at point $B$ obtained with simplified linear and complete nonlinear relations between the reference and equidistant strains, and Abaqus: a) $\ii{\epsilon}{out}{11}$; b)  $\ii{\epsilon}{out}{22}$; c)  $\ii{\sigma}{11}{out}$; d)  $\ii{\sigma}{22}{out}$.}
	\label{fig:pullout: strains at outer point B}
\end{figure}

Moreover, the $D^a$ constitutive model is employed for the calculation of the strains and stresses at the outer surface at point $D$. The strains are calculated using the approximative linear \eqref{eq:linear form of eq strain} and full nonlinear \eqref{eq: components of eq strain rate final} relations between reference and equidistant strains. The stresses are calculated from these strains using \eqref{eq:constitutive relation at equidistant surface}. Comparison with strains and stresses from Abaqus is given in \fref{fig:pullout: strains at outer point}. The dominant strain component is $\ii{\epsilon}{out}{11}$ and the results obtained with Abaqus and \eqqref{eq:linear form of eq strain} differ from the rigorous ones, obtained with \eqqref{eq: components of eq strain rate final}. The observed difference with respect to Abaqus is mainly due to the difference in curvature change $\ii{\kappa}{}{1}$, \fref{fig: pullout: reference strains}c. The other component of strain is smaller for an order of magnitude and all models return similar results, which is again a consequence of the fact that $\ii{b}{1}{1} \gg \ii{b}{2}{2}$ at this point. Regrading the stresses, discrepancies are observed for both components. The noted differences with respect to Abaqus exist mainly due to the differences in reference strains and assumed strain distributions over thickness. Another influence comes from the calculation of the physical components, which are here determined rigorously by the \eqqref{eq:physical strain}. Remarkably, these differences practically vanish in the areas of shell that are not strongly curved. As an example, let us observe the equidistant strains and stresses on the outer surface at point $B$, \fref{fig:pullout: strains at outer point B}. In this area, the curviness for $\textrm{LPF}=1$ is small, $Kh\approx0.018$, and the results for strains and stresses agree well. This underscores that the observed discrepancies at the point $C$ are mainly due to the large curviness of shell at that particular point.

\subsection{Pinched cylinder}

\begin{figure}
	\includegraphics[width=\linewidth]{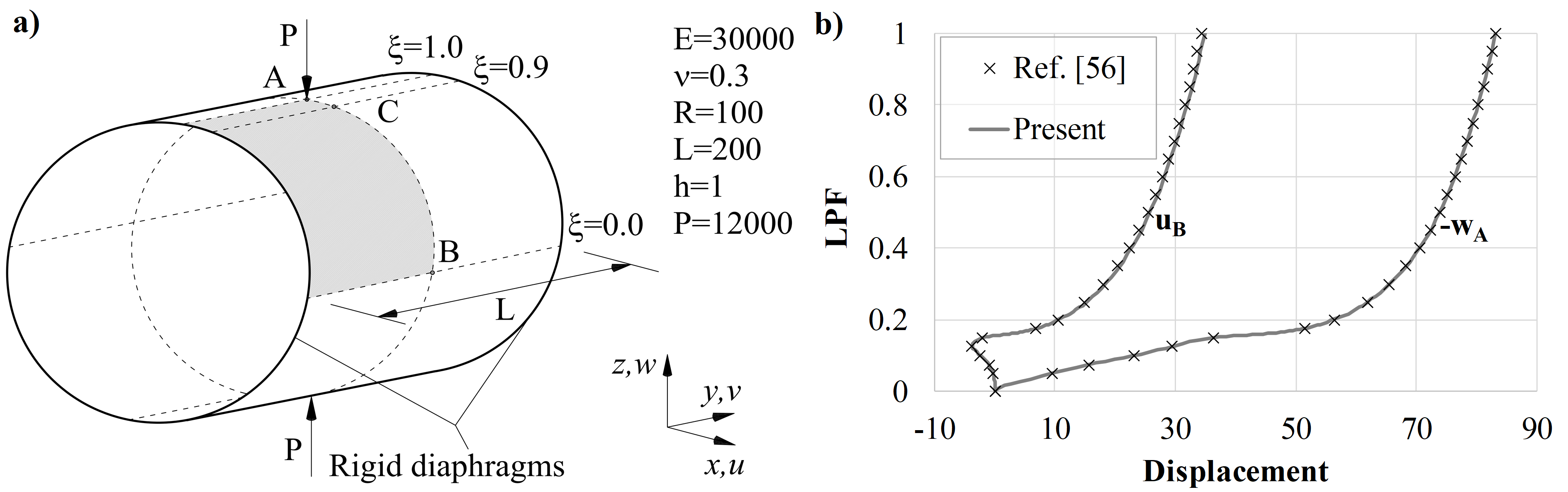}
	\caption{Pinched cylinder. a) Geometry, material characteristics, and applied load. b) Comparison of characteristic displacements components between present formulation and \cite{2004sze}. }
	\label{fig:pinched: disposition and displacements}
\end{figure}

As the final example, the pinched cylinder with rigid diaphragms is considered again, \fref{fig:pinched: disposition and displacements}a. However, the emphasis is now on the geometric nonlinear response, \cite{2004sze}. The reference curviness is $Kh=0.01$, as in the linear example. \textcolor{black}{A mesh with $50\times50$ quadratic $C^1$ elements is employed for the discretization of one eighth of the shell, similar to that in \cite{2017duong}. Two characteristic displacement components are compared in \fref{fig:pinched: disposition and displacements}b.} \textcolor{black}{For the complete equilibrium path, there are no practical differences between the results obtained with the different constitutive models, so they are not displayed separately.} Furthermore, the IGA results are in agreement with the reference solution given in \cite{2004sze}. 

It is already noted that this example is very challenging for the solution using the arc-length method. The tangent stiffness matrix has one negative eigenvalue for $\textrm{LPF}>0.089$ which results with the erroneous negative sign of the predictor solution. As noted, several types of arc-length procedure are tested here, and the one implemented in the Abaqus proved to be the most robust, \cite{2009smith}.

\begin{figure}
	\includegraphics[width=\linewidth]{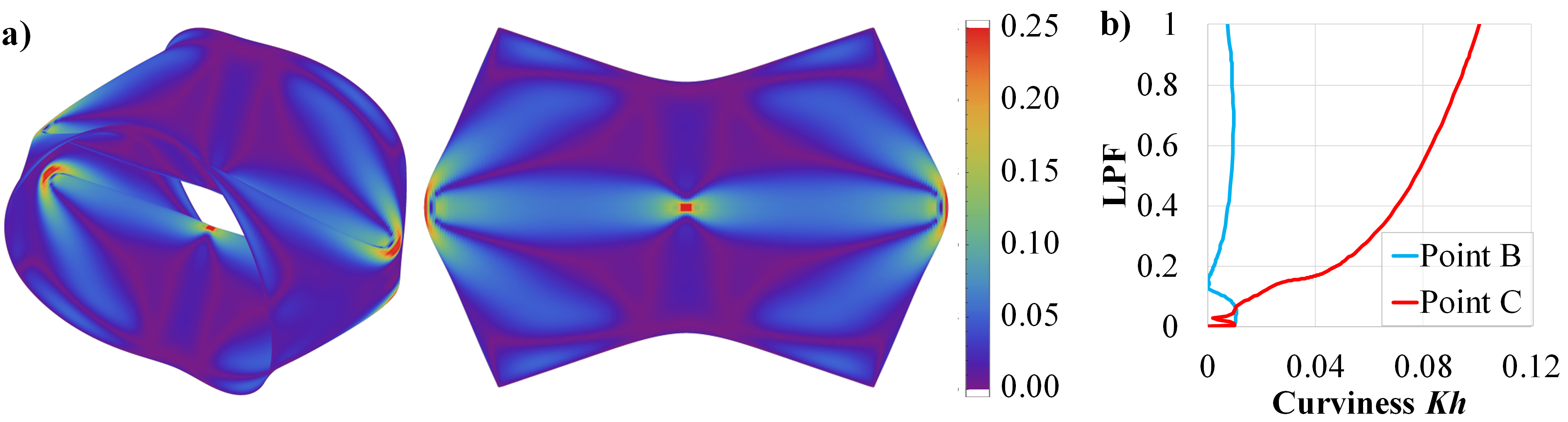}
	\caption{Pinched cylinder. a) Distribution of curviness for $\textrm{LPF}=1$ from two points of view. The threshold value of $Kh=0.25$ is utilized. b) Curviness at points $B$ and $C$ vs. LPF.}
	\label{fig:pinched: curviness}
\end{figure}

Although the initial curviness of the observed shell is $Kh=0.01$, the finite deformation resulted in the final configuration with a maximum curviness at point $A$ of $Kh=0.346$. The distribution of the curviness for $\textrm{LPF}=1$ is given in \fref{fig:pinched: curviness}a. It is evident that some parts of this shell become strongly curved locally. The change of curviness with respect to the load proportionality factor is shown in \fref{fig:pinched: curviness}b for two characteristic points.

\begin{figure}
	\includegraphics[width=\linewidth]{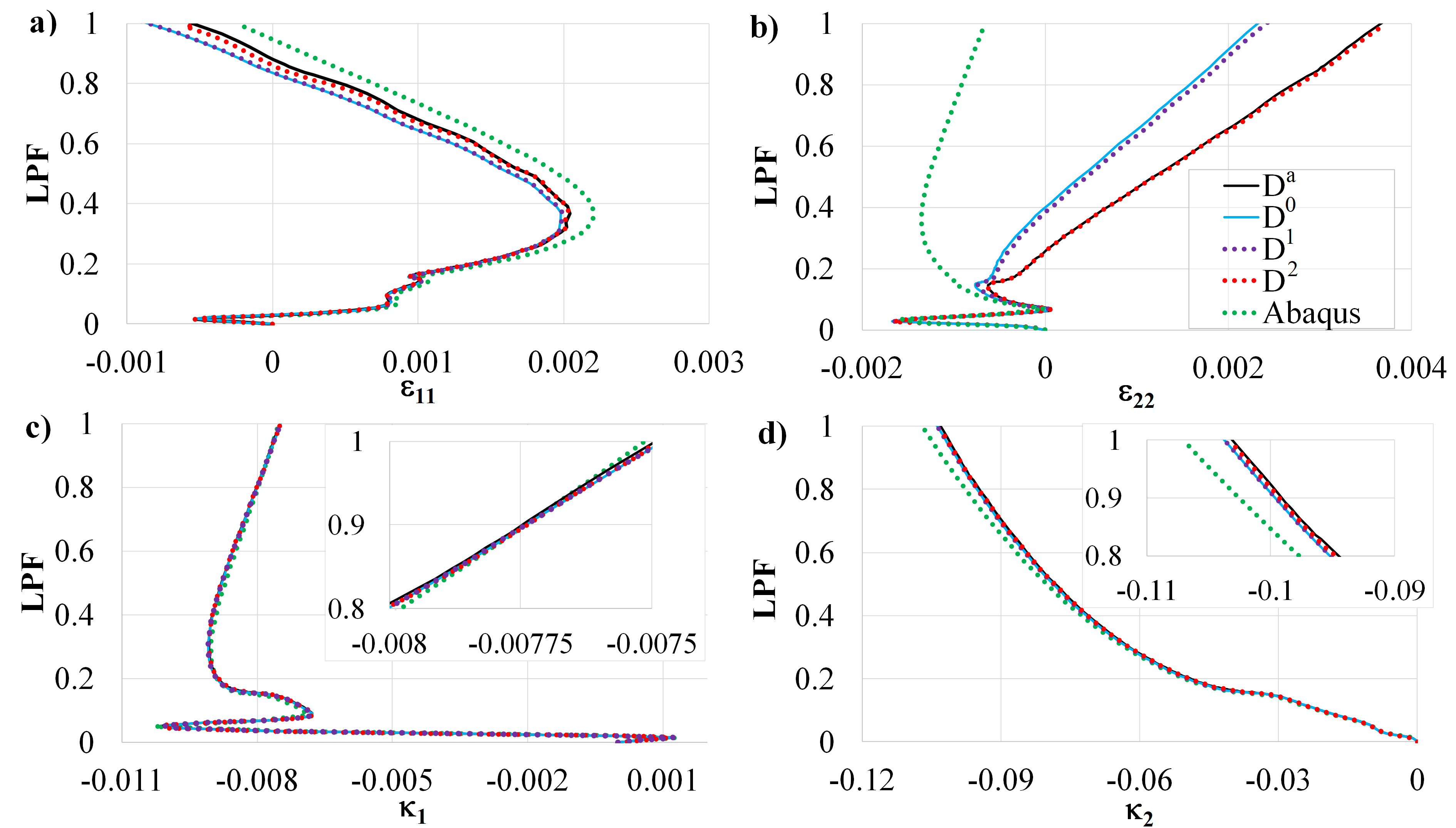}
	\caption{Pinched cylinder. Comparison of reference strains at point $C$ obtained with four constitutive models and Abaqus: a) $\ii{\epsilon}{}{11}$; b) $\ii{\epsilon}{}{22}$; c) $\ii{\kappa}{}{1}$; d) $\ii{\kappa}{}{2}$.}
	\label{fig:pinched: ref strains}
\end{figure}

Next, the strains at point $C$ are considered in detail. There, the shell is strongly curved, i.e., $Kh=0.101$ for $\textrm{LPF}=1$. The four reference strains are calculated using the four different constitutive models and compared in \fref{fig:pinched: ref strains} with the results from Abaqus, which employs a mesh with 62800 S4R elements. In comparison to the previous example, even more complex equilibrium paths are observed. The dominant reference strain is $\ii{\kappa}{}{2}$ and it is almost identical for all studied models, \fref{fig:pinched: ref strains}d. \textcolor{black}{An exception is the Abaqus model which gives different results for $\textrm{LPF}>0.4$, with a maximum relative difference of nearly $3 \% $ for LPF=1. Considering the other component of the change of curvature, all models return similar results, \fref{fig:pinched: ref strains}c.} Regarding the reference membrane strains, they are in agreement for $\textrm{LPF}<0.1$ but clear differences exist for the greater part of equilibrium path which is especially pronounced for $\ii{\epsilon}{}{22}$, \fref{fig:pinched: ref strains}b. The constitutive models $D^a$ and $D^2$ return practically the same results, while there are small differences of membrane strains for $\textrm{LPF}>0.7$. Additionally, the models $D^0$ and $D^1$ are in close agreement. 

\begin{figure}
	\includegraphics[width=\linewidth]{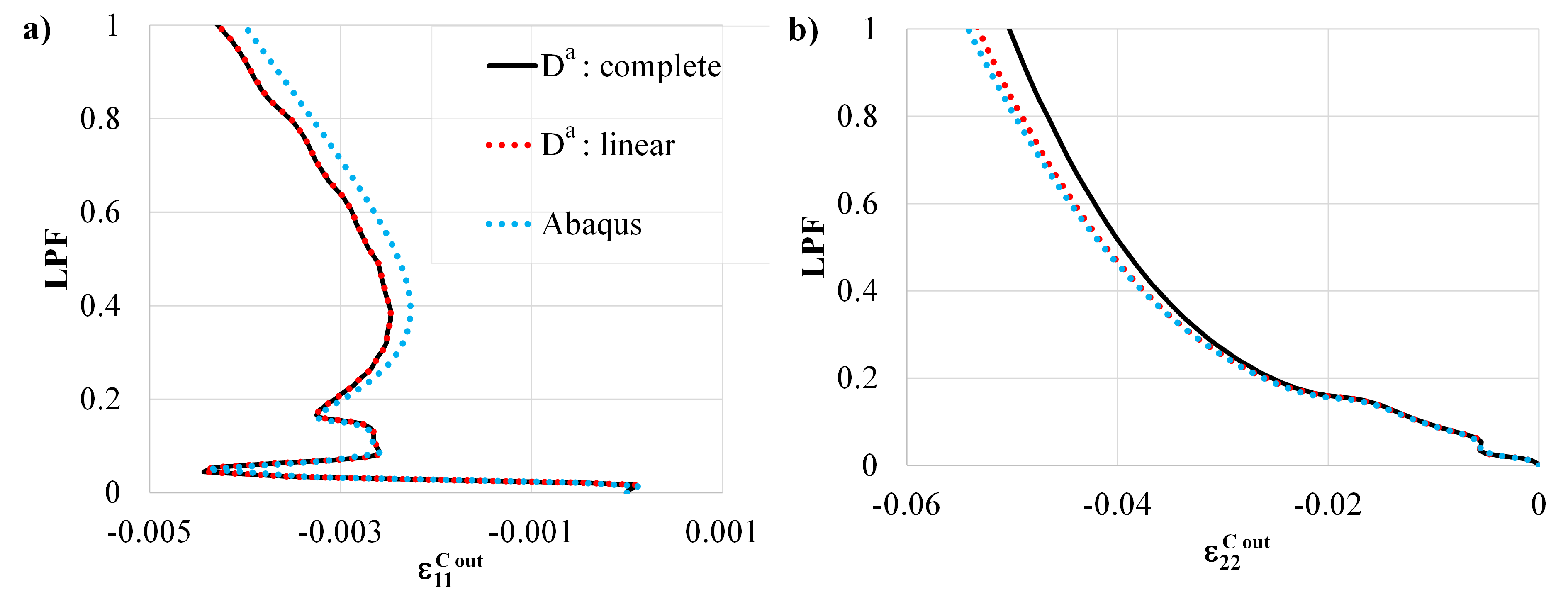}
	\caption{Pinched cylinder. Comparison of strains at outer surface at point $C$ obtained with linear and complete nonlinear relations between the reference and equidistant strains, and Abaqus: a) $\ii{\epsilon}{out}{11}$; b) $\ii{\epsilon}{out}{22}$.}
	\label{fig:pinched: eq strains}
\end{figure}

Finally, the $D^a$ constitutive model is employed and the strains at the outer surface are calculated using approximative \eqref{eq:linear form of eq strain} and exact \eqref{eq: components of eq strain rate final} relations between reference and equidistant strains, \fref{fig:pinched: eq strains}. Good agreement of all models for $\textrm{LPF}<0.2$ is observed for the dominant strain, $\ii{\epsilon}{out}{22}$. However, as the load and curviness increase, the differences become noticeable. The results obtained with approximative relation and Abaqus are close since the dominant reference strain is $\ii{\kappa}{}{2}$, \fref{fig:pinched: ref strains}c. The results of the other equidistant strain component are indifferent with respect to the Eqs.~\eqref{eq: components of eq strain rate final} and \eqref{eq:linear form of eq strain} due to the fact that $\ii{b}{1}{1} \approx 0$ at point $C$.

\begin{figure}
	\includegraphics[width=\linewidth]{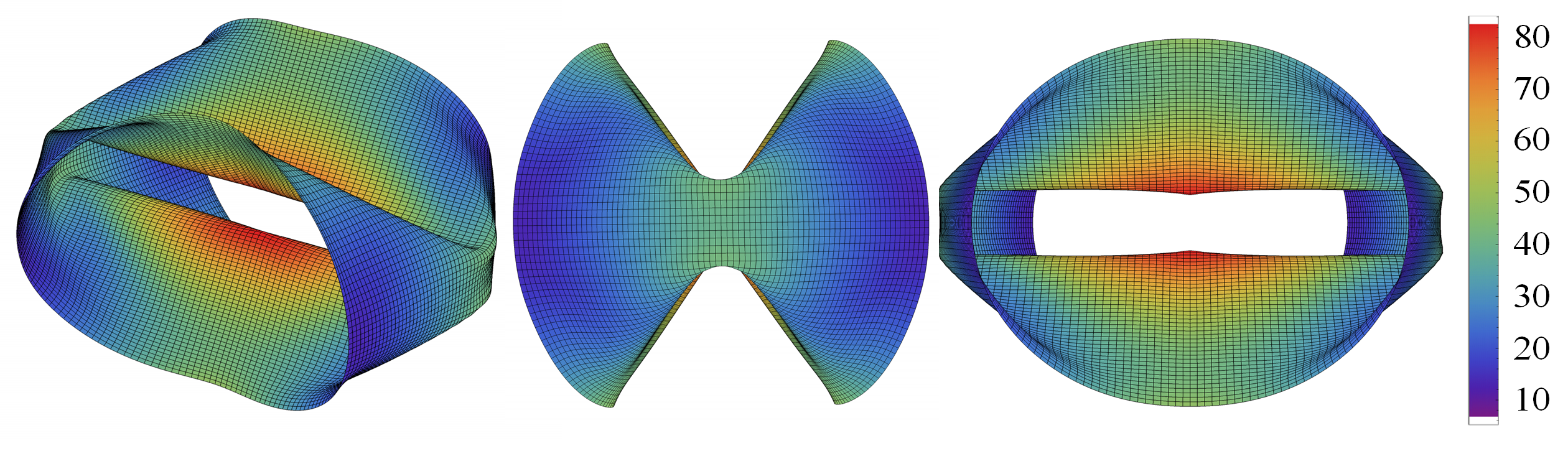}
	\caption{Pinched cylinder. Deformed configuration for $\textrm{LPF}=1$ from three perspectives. Color shows the total value of displacement vector. }
	\label{fig:pinched: def config}
\end{figure}

\section{Conclusions}

A rigorous geometric nonlinear IGA shell formulation is presented and validated by thorough numerical experiments. A key finding is that the initial configuration is not enough to determine the slenderness of the structure; all subsequent deformed configurations must be considered. As a result, some standard approximations of applied shell theories lead to erroneous results for strongly curved shells.

We strictly derive the nonlinear equation of equilibrium by considering the geometrically exact relations between the equidistant and reference strains of a KL shell. The isogeometric approach is applied for the spatial discretization and the arc-length method is employed as a solution procedure. The geometric stiffness matrix is derived by the variation of strains with respect to the metric. The symmetry of this matrix is maintained which proves that the selected force and strain quantities are energetically conjugated.

The main contribution of the paper is the full analytical integration of the virtual power along the thickness direction which is performed using the reciprocal shift tensor. The resulting $D^a$ model is complex and thus, computationally expensive. However, this constitutive model is exact, and allows us to derive different reduced constitutive models $D^0$, $D^1$, and $D^2$. Based on the present numerical experiments, it follows that the approximative constitutive model $D^2$ is well-suited for the nonlinear modeling of strongly curved shells. It is as efficient as the simple $D^0$ model and returns results which are virtually indistinguishable from the exact ones. On the other hand, the displacement field of thin, small-curvature shells, is practically indifferent to the more rigorous constitutive relations, no matter how large the displacements and rotations are and the simple decoupled equations of the $D^0$ model can return reasonably accurate results. Regardless of the constitutive model used, nonlinear distribution of strain through the thickness should be considered during the post-processing phase in order to obtain accurate strain field of strongly curved shells. 

A detailed inspection of the observed differences between IGA and Abaqus was beyond the scope of the present research and it should be a subject of further research. The development of an adaptive refinement scheme that applies higher regularity on smooth parts of the deformed structure, and lower regularity at the kinks and similar irregularities would be an interesting improvement. Additionally, the different constitutive relations could be applied, depending on the curviness of the observed area, like the blended shell approach suggested in \cite{2013benson}. More accurate shell models should also include the change of thickness as well, since this directly influences the curviness.

\section*{Acknowledgments}

The first author expresses his gratitude to Professor Sava Vukelić (1938-1992) for the inspiring discussions and infectious enthusiasm regarding the mechanics.

\textcolor{black}{During this work, our beloved colleague and friend, Professor Gligor Radenković (1956-2019), passed away. The second author acknowledges that his unprecedented enthusiasm and love for mechanics were crucial for much of his previous, present, and future research. His students and colleagues will remember him by his incredibly infectious energy for both life and science.}

The presented research is partially supported by the Ministry for Scientific-Technological Development, Higher Education and Information Society of Republic of Srpska through the project \textit{Dynamic stability of thin-walled structures using the Isogeometric finite strip method}. This support is gratefully acknowledged. 

We also acknowledge the support of the Austrian Science Fund (FWF): M 2806-N.

\section*{Appendix A. Strain at an equidistant surface}

\setcounter{equation}{0}
\renewcommand\theequation{A\arabic{equation}}

If we insert Eqs.~\eqref{eq: velocity gradient equidistant} and \eqref{eq: derivative of v3} into \eqqref{eq:def: strain rate equidistant wrt Cartesian}, the components of $\ieq{d}{}{\alpha\beta}$ can be written as: 
\begin{equation}
\label{A1}
\begin{aligned}
	2 \ieq{d}{}{\alpha\beta} &= \left( \idef{x}{n}{,\alpha} - \zeta \idef{b}{\nu}{\alpha} \idef{x}{n}{,\nu}\right) \left( \ii{v}{}{n,\beta} + \zeta \ii{v}{}{n,3\beta}\right) + \left( \idef{x}{n}{,\beta} - \zeta \idef{b}{\nu}{\beta} \idef{x}{n}{,\nu}\right) \left( \ii{v}{}{n,\alpha} + \zeta \ii{v}{}{n,3\alpha}\right) \\
	&= \idef{x}{n}{,\alpha} \ii{v}{}{n,\beta} + \idef{x}{n}{,\beta} \ii{v}{}{n,\alpha} - \zeta \left( \idef{b}{\nu}{\alpha} \idef{x}{n}{,\nu} \ii{v}{}{n,\beta} + \idef{b}{\nu}{\beta} \idef{x}{n}{,\nu} \ii{v}{}{n,\alpha} - \idef{x}{n}{,\alpha} \ii{v}{}{n,3\beta} - \idef{x}{n}{,\beta} \ii{v}{}{n,3\alpha}\right) \\
	&\;\;\;\; - \zeta^2 \idef{x}{n}{,\nu} \left( \idef{b}{\nu}{\alpha} \ii{v}{}{n,3\beta} + \idef{b}{\nu}{\beta} \ii{v}{}{n,3\alpha}\right) \\
	&= 2 \ii{d}{}{\alpha\beta} - \zeta \left( \idef{b}{\nu}{\alpha} \idef{x}{n}{,\nu} \ii{v}{}{n,\beta} + \idef{b}{\nu}{\beta} \idef{x}{n}{,\nu} \ii{v}{}{n,\alpha}  \right.\\
	& \left.\;\;\;\; - \idef{x}{n}{,\alpha} \ii{v}{}{n,3\beta} - \idef{x}{n}{,\beta} \ii{v}{}{n,3\alpha}\right) - \zeta^2 \idef{x}{n}{,\nu} \left( \idef{b}{\nu}{\alpha} \ii{v}{}{n,3\beta} + \idef{b}{\nu}{\beta} \ii{v}{}{n,3\alpha}\right).
\end{aligned}
\end{equation}
\noindent By utilizing \eqqref{eq: derivative of v3}, the terms that are multiplied by $\zeta$ can be transformed to: 
\begin{equation}
\label{A2}
\begin{aligned}
	& \idef{x}{n}{,\nu} \left( \idef{b}{\nu}{\alpha} \ii{v}{}{n,\beta} + \idef{b}{\nu}{\beta} \ii{v}{}{n,\alpha}\right) - \idef{x}{n}{,\alpha} \left( \idef{\Gamma}{\mu}{\beta\gamma} \idef{x}{,\gamma}{n} \ii{v}{}{k,\mu} \idef{x}{k}{,3} - 	\idef{b}{\mu}{\beta} \idef{x}{,3}{n} \ii{v}{}{k,\mu} \idef{x}{k}{,3} - \idef{x}{,\mu}{n} \ii{v}{}{k,\mu\beta} \idef{x}{k}{,3} \right. \\ 
	& \left. + \idef{x}{,\mu}{n} \ii{v}{}{k,\mu} \idef{b}{\nu}{\beta} \idef{x}{k}{,\nu}\right) - \idef{x}{n}{,\beta} \left( \idef{\Gamma}{\mu}{\alpha\gamma} \idef{x}{,\gamma}{n} \ii{v}{}{k,\mu} \idef{x}{k}{,3} - \idef{b}{\mu}{\alpha} \idef{x}{,3}{n} \ii{v}{}{k,\mu} \idef{x}{k}{,3} - \idef{x}{,\mu}{n} \ii{v}{}{k,\mu\alpha} \idef{x}{k}{,3} \right. \\
	& \left. + \idef{x}{,\mu}{n} \ii{v}{}{k,\mu} \idef{b}{\nu}{\alpha} \idef{x}{k}{,\nu} \right) = \idef{x}{n}{,\nu} \left( \idef{b}{\nu}{\alpha} \ii{v}{}{n,\beta} + \idef{b}{\nu}{\beta} \ii{v}{}{n,\alpha}\right) - \idef{\Gamma}{\mu}{\beta\alpha} \ii{v}{}{k,\mu} \idef{x}{k}{,3} + \ii{v}{}{k,\alpha\beta} \idef{x}{k}{,3} - \ii{v}{}{k,\alpha} \idef{b}{\nu}{\beta} \idef{x}{k}{,\nu} \\
	&- \idef{\Gamma}{\mu}{\beta\alpha} \ii{v}{}{k,\mu} \idef{x}{k}{,3} + \ii{v}{}{k,\alpha\beta} \idef{x}{k}{,3} - \ii{v}{}{k,\beta} \idef{b}{\nu}{\alpha} \idef{x}{k}{,\nu} = 2 \idef{x}{k}{,3} \left( \ii{v}{}{k,\alpha\beta} - \idef{\Gamma}{\mu}{\alpha\beta} \ii{v}{}{k,\mu}\right).
\end{aligned}
\end{equation}
\noindent Analogously, the expression in \eqqref{A1} that is multiplied by $\zeta^2$ can be written as: 
\begin{equation}
\label{A3}
\begin{aligned}
	&\idef{x}{n}{,\nu} \left( \idef{b}{\nu}{\alpha} \ii{v}{}{n,3\beta} + \idef{b}{\nu}{\beta} \ii{v}{}{n,3\alpha} \right) \\
	&= \idef{b}{\nu}{\alpha} \left( \idef{\Gamma}{\mu}{\beta\nu} \ii{v}{}{k,\mu} \idef{x}{k}{,3} - \ii{v}{}{k,\nu\beta} \idef{x}{k}{,3} + \ii{v}{}{k,\nu} \idef{b}{\mu}{\beta} \idef{x}{k}{,\mu}\right) \\
	&\;\;\;\; + \idef{b}{\nu}{\beta} \left( \idef{\Gamma}{\mu}{\alpha\nu} \ii{v}{}{k,\mu} \idef{x}{k}{,3} - \ii{v}{}{k,\nu\alpha} \idef{x}{k}{,3} + \ii{v}{}{k,\nu} \idef{b}{\mu}{\alpha} \idef{x}{k}{,\mu}\right) \\
	&= \left( \idef{b}{\nu}{\alpha} \idef{b}{\mu}{\beta} + \idef{b}{\nu}{\beta} \idef{b}{\mu}{\alpha}\right) \ii{v}{}{k,\nu} \idef{x}{k}{,\mu} - \idef{b}{\nu}{\alpha} \idef{x}{k}{,3} \left( \ii{v}{}{k,\nu\beta} - \idef{\Gamma}{\mu}{\nu\beta} \ii{v}{}{k,\mu}\right) - \idef{b}{\nu}{\beta} \idef{x}{k}{,3} \left( \ii{v}{}{k,\nu\alpha} - \idef{\Gamma}{\mu}{\nu\alpha} \ii{v}{}{k,\mu} \right).
\end{aligned}
\end{equation}
\noindent If we introduce designations for the following quantities:
\begin{equation}
\label{A4}
	\imd{\kappa}{}{\alpha\beta} = \idef{x}{k}{,3} \left( \ii{v}{}{k,\alpha\beta} - \idef{\Gamma}{\mu}{\alpha\beta} \ii{v}{}{k,\mu} \right),
\end{equation}
\begin{equation}
\label{A5}
	2 \imd{\chi}{}{\alpha\beta} = \left( \idef{b}{\nu}{\alpha} \idef{b}{\mu}{\beta} + \idef{b}{\nu}{\beta} \idef{b}{\mu}{\alpha}\right) \ii{v}{}{k,\nu} \idef{x}{k}{,\mu} - \idef{b}{\nu}{\alpha} \imd{\kappa}{}{\nu\beta} - \idef{b}{\nu}{\beta} \imd{\kappa}{}{\nu\alpha},
\end{equation}
\noindent the strain rate at an arbitrary point can be written as:
\begin{equation}
\label{A6}
	\ieq{d}{}{\alpha\beta} = \ii{d}{}{\alpha\beta} - \zeta \imd{\kappa}{}{\alpha\beta} - \ii{\zeta}{2}{} \imd{\chi}{}{\alpha\beta}.
\end{equation}
\noindent It is evident that $\imd{\chi}{}{\alpha\beta}$ is a complicated quantity which consists of the products of curvatures, curvature change rates, and gradients of velocity. If we notice that $\nu$ and $\mu$ are dummy indices, \eqqref{A5} reduces to:
\begin{equation}
\label{A7}
	2 \imd{\chi}{}{\alpha\beta} =  \idef{b}{\nu}{\alpha} \idef{b}{\mu}{\beta} \left( \ii{v}{}{k,\nu} \idef{x}{k}{,\mu} + \ii{v}{}{k,\mu} \idef{x}{k}{,\nu} \right) - \idef{b}{\nu}{\alpha} \imd{\kappa}{}{\nu\beta} - \idef{b}{\nu}{\beta} \imd{\kappa}{}{\nu\alpha} = 2 \idef{b}{\nu}{\alpha} \idef{b}{\mu}{\beta} \ii{d}{}{\nu\mu} - \idef{b}{\nu}{\alpha} \imd{\kappa}{}{\nu\beta} - \idef{b}{\nu}{\beta} \imd{\kappa}{}{\nu\alpha}.
\end{equation}
\noindent Using this relation, the strain rate at an arbitrary point can be expressed as a function of the strain and curvature change rates of the shell midsurface as:
\begin{equation}
\label{A8}
	\ieq{d}{}{\alpha\beta} = \langle \ieqdef{C}{\mu}{\alpha} \left[ \left( \ii{\delta}{\nu}{\beta} + \zeta \idef{b}{\nu}{\beta} \right) \ii{d}{}{\mu\nu} - \zeta \imd{\kappa}{}{\mu\beta} \right]\rangle^{sym}, \quad \ieqdef{C}{\mu}{\alpha} = \ii{\delta}{\mu}{\alpha} - \zeta \idef{b}{\mu}{\alpha},
\end{equation}

\noindent where the symmetry of the expression refers to the symmetry by indices $\alpha$ and $\beta$, and $\mu$ and $\nu$, that is:
\begin{equation}
\label{A9}
	2 \ieq{d}{}{\alpha\beta} = \ieqdef{C}{\mu}{\alpha} \left[ \left( \ii{\delta}{\nu}{\beta} + \zeta \idef{b}{\nu}{\beta} \right) \ii{d}{}{\mu\nu} - \zeta \imd{\kappa}{}{\mu\beta} \right] + \ieqdef{C}{\nu}{\beta} \left[ \left( \ii{\delta}{\mu}{\alpha} + \zeta \idef{b}{\mu}{\alpha} \right) \ii{d}{}{\mu\nu} - \zeta \imd{\kappa}{}{\nu\alpha} \right].
\end{equation}
\noindent By introducing: 
\begin{equation}
\label{A10}
\begin{aligned}
	\ii{A}{\mu\nu}{\alpha\beta} &= \frac{1}{2} \left[ \ieqdef{C}{\mu}{\alpha} \left( \ii{\delta}{\nu}{\beta} + \zeta \idef{b}{\nu}{\beta} \right) + \ieqdef{C}{\nu}{\beta} \left( \ii{\delta}{\mu}{\alpha} + \zeta \idef{b}{\mu}{\alpha} \right) \right] = \ii{\delta}{\mu}{\alpha} \ii{\delta}{\nu}{\beta} - \zeta^2 \idef{b}{\mu}{\alpha} \idef{b}{\nu}{\beta},\\
	\ii{B}{\mu\nu}{\alpha\beta} &= \frac{1}{2} \left( \ii{\delta}{\nu}{\alpha} \ieqdef{C}{\mu}{\beta} + \ii{\delta}{\mu}{\beta} \ieqdef{C}{\nu}{\alpha} \right),
\end{aligned}
\end{equation}
\noindent we finally obtain the compact form for the strain rate at an equidistant surface:
\begin{equation}
\label{A11}
	\ieq{d}{}{\alpha\beta} = \ii{A}{\mu\nu}{\alpha\beta} \ii{d}{}{\mu\nu} - \zeta \ii{B}{\mu\nu}{\alpha\beta} \imd{\kappa}{}{\mu\nu}.
\end{equation}

\section*{Appendix B. Closed-form of the constitutive matrix}

\setcounter{equation}{0}
\renewcommand\theequation{B\arabic{equation}}

The elements of the constitutive matrix given by \eqqref{eq: definition of parts of constitutive tensor} are now presented in an expanded form. 
\begin{equation}
\label{B1}
\begin{aligned}
	\ii{D}{\gamma\lambda\chi\omega}{M} &= \left(\ii{I}{}{0} - 2 T \ii{I}{}{1} + \ii{T}{2}{} \ii{I}{}{2}\right) \ii{D}{\gamma\lambda\chi\omega}{} \\
	&\;\;\;\; + \frac{1}{2} \left( \ii{b}{\lambda}{\nu} \ii{K}{1}{M} + \ii{b}{\lambda}{\mu} \ii{b}{\mu}{\nu} \ii{K}{2}{M} \right) \ii{D}{\gamma\nu\chi\omega}{}  + \frac{1}{2} \left( \ii{b}{\omega}{\varphi} \ii{K}{1}{M} + \ii{b}{\omega}{\mu} \ii{b}{\mu}{\varphi} \ii{K}{2}{M} \right) \ii{D}{\gamma\lambda\chi\varphi}{} \\
	&\;\;\;\; + \frac{1}{2} \left( \ii{b}{\chi}{\varphi} \ii{K}{1}{M} + \ii{b}{\chi}{\mu} \ii{b}{\mu}{\varphi} \ii{K}{2}{M} \right) \ii{D}{\gamma\lambda\varphi\omega}{}
	+ \frac{1}{2} \left( \ii{b}{\gamma}{\nu} \ii{K}{1}{M} + \ii{b}{\gamma}{\mu} \ii{b}{\mu}{\nu} \ii{K}{2}{M} \right) \ii{D}{\nu\lambda\chi\omega}{}\\
	&\;\;\;\; + \frac{1}{4} \left[ \ii{b}{\omega}{\varphi} \ii{b}{\lambda}{\nu} \ii{K}{3}{M} + \left( \ii{b}{\omega}{\mu} \ii{b}{\mu}{\varphi} \ii{b}{\lambda}{\nu} + \ii{b}{\omega}{\varphi} \ii{b}{\lambda}{\mu} \ii{b}{\mu}{\nu}\right) \ii{K}{4}{M} + \ii{b}{\omega}{\epsilon} \ii{b}{\epsilon}{\varphi} \ii{b}{\lambda}{\mu} \ii{b}{\mu}{\nu} \ii{I}{}{4} \right] \ii{D}{\gamma\nu\chi\varphi}{}\\
	&\;\;\;\; + \frac{1}{4} \left[ \ii{b}{\omega}{\varphi} \ii{b}{\gamma}{\nu} \ii{K}{3}{M} + \left( \ii{b}{\omega}{\mu} \ii{b}{\mu}{\varphi} \ii{b}{\gamma}{\nu} + \ii{b}{\omega}{\varphi} \ii{b}{\gamma}{\mu} \ii{b}{\mu}{\nu}\right) \ii{K}{4}{M} + \ii{b}{\omega}{\epsilon} \ii{b}{\epsilon}{\varphi} \ii{b}{\gamma}{\mu} \ii{b}{\mu}{\nu} \ii{I}{}{4} \right] \ii{D}{\nu\lambda\chi\varphi}{}\\
	&\;\;\;\; + \frac{1}{4} \left[ \ii{b}{\chi}{\varphi} \ii{b}{\lambda}{\nu} \ii{K}{3}{M} + \left( \ii{b}{\chi}{\mu} \ii{b}{\mu}{\varphi} \ii{b}{\lambda}{\nu} + \ii{b}{\chi}{\varphi} \ii{b}{\lambda}{\mu} \ii{b}{\mu}{\nu}\right) \ii{K}{4}{M} + \ii{b}{\chi}{\epsilon} \ii{b}{\epsilon}{\varphi} \ii{b}{\lambda}{\mu} \ii{b}{\mu}{\nu} \ii{I}{}{4} \right] \ii{D}{\gamma\nu\varphi\omega}{}\\
	&\;\;\;\; + \frac{1}{4} \left[ \ii{b}{\chi}{\varphi} \ii{b}{\gamma}{\nu} \ii{K}{3}{M} + \left( \ii{b}{\chi}{\mu} \ii{b}{\mu}{\varphi} \ii{b}{\gamma}{\nu} + \ii{b}{\chi}{\varphi} \ii{b}{\gamma}{\mu} \ii{b}{\mu}{\nu}\right) \ii{K}{4}{M} + \ii{b}{\chi}{\epsilon} \ii{b}{\epsilon}{\varphi} \ii{b}{\gamma}{\mu} \ii{b}{\mu}{\nu} \ii{I}{}{4} \right] \ii{D}{\nu\lambda\varphi\omega}{},
\end{aligned}
\end{equation}
\begin{equation}
\label{B2}
\begin{aligned}
	- \ii{D}{\lambda\gamma\chi\omega}{MB} &= \left(\ii{I}{}{1} - 2 T \ii{I}{}{2} + \ii{T}{2}{} \ii{I}{}{3}\right) \ii{D}{\lambda\gamma\chi\omega}{} + \frac{1}{2} \ii{K}{2}{M} \left( \ii{b}{\gamma}{\nu} \ii{D}{\lambda\nu\chi\omega}{} + \ii{b}{\lambda}{\nu} \ii{D}{\nu\gamma\chi\omega}{} \right) \\
	&\;\;\;\; + \frac{1}{2} \left( \ii{b}{\omega}{\varphi} \ii{K}{1}{MB} + \ii{b}{\omega}{\mu} \ii{b}{\mu}{\varphi} \ii{K}{2}{MB} \right) \ii{D}{\lambda\gamma\chi\varphi}{} + \frac{1}{2} \left( \ii{b}{\chi}{\varphi} \ii{K}{1}{MB} + \ii{b}{\chi}{\mu} \ii{b}{\mu}{\varphi} \ii{K}{2}{MB} \right) \ii{D}{\lambda\gamma\varphi\omega}{} \\
	&\;\;\;\; + \frac{1}{4} \ii{b}{\gamma}{\nu} \left( \ii{b}{\omega}{\varphi} \ii{K}{3}{MB} + \ii{b}{\omega}{\mu} \ii{b}{\mu}{\varphi} \ii{I}{}{4} \right) \ii{D}{\lambda\nu\chi\varphi}{} + \frac{1}{4} \ii{b}{\lambda}{\nu} \left( \ii{b}{\omega}{\varphi} \ii{K}{3}{MB} + \ii{b}{\omega}{\mu} \ii{b}{\mu}{\varphi} \ii{I}{}{4} \right) \ii{D}{\nu\gamma\chi\varphi}{} \\
	&\;\;\;\; + \frac{1}{4} \ii{b}{\gamma}{\nu} \left( \ii{b}{\chi}{\varphi} \ii{K}{3}{MB} + \ii{b}{\chi}{\mu} \ii{b}{\mu}{\varphi} \ii{I}{}{4} \right) \ii{D}{\lambda\nu\varphi\omega}{} +\frac{1}{4} \ii{b}{\lambda}{\nu} \left( \ii{b}{\chi}{\varphi} \ii{K}{3}{MB} + \ii{b}{\chi}{\mu} \ii{b}{\mu}{\varphi} \ii{I}{}{4} \right) \ii{D}{\nu\gamma\varphi\omega}{},
\end{aligned}
\end{equation}
\begin{equation}
\label{B3}
\begin{aligned}
	\ii{D}{\lambda\gamma\chi\omega}{B} &= \left(\ii{I}{}{2} - 2 T \ii{I}{}{3} + \ii{T}{2}{} \ii{I}{}{4}\right) \ii{D}{\lambda\gamma\chi\omega}{} \\
	&\;\;\;\; + \frac{1}{2} \left( \ii{I}{}{3} - T \ii{I}{}{4} \right) \left( \ii{b}{\lambda}{\nu} \ii{D}{\nu\gamma\chi\omega}{} + \ii{b}{\gamma}{\nu} \ii{D}{\lambda\nu\chi\omega}{} + \ii{b}{\chi}{\nu} \ii{D}{\lambda\gamma\nu\omega}{} + \ii{b}{\omega}{\nu} \ii{D}{\lambda\gamma\chi\nu}{}\right) \\
	&\;\;\;\; + \frac{1}{4} \ii{I}{}{4} \left( \ii{b}{\lambda}{\nu} \ii{b}{\chi}{\varphi} \ii{D}{\nu\gamma\varphi\omega}{} + \ii{b}{\lambda}{\nu} \ii{b}{\omega}{\varphi} \ii{D}{\nu\gamma\chi\varphi}{} + \ii{b}{\gamma}{\nu} \ii{b}{\chi}{\varphi} \ii{D}{\lambda\nu\varphi\omega}{} + \ii{b}{\gamma}{\nu} \ii{b}{\omega}{\varphi} \ii{D}{\lambda\nu\chi\varphi}{}\right).
\end{aligned}
\end{equation}
\noindent Here, the following designations are introduced:
\begin{equation}
\label{B4}
\begin{aligned}
	\ii{K}{1}{M} &= 2 \ii{I}{}{1} - 3 T \ii{I}{}{2} + \ii{T}{2}{} \ii{I}{}{3}, \quad
	\ii{K}{2}{M} =\ \ii{I}{}{2} - T \ii{I}{}{3}, \; \ii{K}{3}{M} = 4 \ii{I}{}{2} - 4 T \ii{I}{}{3} + \ii{T}{2}{} \ii{I}{}{4}, \\
	\ii{K}{4}{M} &= 2 \ii{I}{}{3} - T \ii{I}{}{4}, \quad \ii{K}{1}{MB} =  2 \ii{I}{}{2} - 3 T \ii{I}{}{3} + \ii{T}{2}{} \ii{I}{}{4}, \quad \ii{K}{2}{MB} =  \ii{I}{}{3} - T \ii{I}{}{4}, \\
	\ii{K}{3}{MB} &= 2 \ii{I}{}{3} - T \ii{I}{}{4},
\end{aligned}
\end{equation}
\noindent while $T$ denotes the trace of the curvature tensor. Finally, utilizing symbolic capabilities of Wolfram Mathematica, following closed-form solutions for the integrals $I_0-I_4$ are obtained:
\begin{equation}
\label{B5}
\begin{aligned}
	p &= \sqrt{4 b - \ii{T}{2}{}}, \quad m = \frac{b h - T}{p}, \quad n = \frac{b h + T}{p}, \\
	r &= 4 + b \ii{h}{2}{} + 2 h T, \quad s = 4 + b h^2 - 2 h T,  \\
	\ii{I}{}{0} &= \int_{-h/2}^{h/2} \frac{1}{1-\zeta T + \zeta^2 b} \dd{\zeta} = \frac{2}{p} \left[ \arctan(m) + \arctan(n)\right],\\
	\ii{I}{}{1} &= \int_{-h/2}^{h/2} \frac{\zeta}{1-\zeta T + \zeta^2 b} \dd{\zeta} \\
	&= \frac{1}{2bp} \left\{ 2 T \left[ \arctan(m) + \arctan(n)\right] + p \ln(\frac{8 + 2 b h^2}{r} - 1) \right\},\\
	\ii{I}{}{2} &= \int_{-h/2}^{h/2} \frac{\zeta^2}{1-\zeta T + \zeta^2 b} \dd{\zeta} \\
	&= \frac{1}{b^2 p} \left\{ b h p - \left( 2 b - T^2 \right) \left[ \arctan(m) + \arctan(n)\right] - T p \arctan(\frac{2 h T}{4 + b h^2} ) \right\},\\
	\ii{I}{}{3} &= \int_{-h/2}^{h/2} \frac{\zeta^3}{1-\zeta T + \zeta^2 b} \dd{\zeta} \\
	&= \frac{1}{2 b^3 p} \left\{ \left( 6 b T - 2 T^3 \right) \arctan(-n) + \left( 2 T^3 - 6 b T \right) \arctan(m)  \right.\\ 	 
	&\left. \;\;\;\; + p \left[ 2 b h T + \left( T^2 - b \right) \ln(s)  + \left( b-T^2 \right)  \ln(r)  \right]\right\}, \\
	\ii{I}{}{4} &= \int_{-h/2}^{h/2} \frac{\zeta^4}{1-\zeta T + \zeta^2 b} \dd{\zeta} \\
	&= \frac{1}{12 b^4 p} \left\{ -12 \left( 2 b^2 - 4 b T^2 + T^4 \right) \arctan(-n) + 12 \left( 2 b^2 - 4 b T^2 + T^4 \right) \arctan(m) \right. \\
	&\left. \;\;\;\; + p \left[ b h \left( -12 b + b^2 h^2 + 12 T^2 \right) + 6 T \left( T^2 - 2 b \right) \ln(s) + 6 \left( 2 b T - h^3 \right) \ln(r) \right]\right\}.
\end{aligned}
\end{equation}
\noindent However, it should be noted that these closed-form integrals can become unstable for small-curvature shells and an appropriate switch in the algorithm is required in order to include reduced expressions for which the standard thin shell assumptions are valid. Therefore, the presented expressions are ideally suited for strongly curved parts of a shell, but the exact limit between small and large curvature shells is arbitrary and some criteria must be introduced, e.g., $Kh>0.05$.

\section*{Appendix C. Different reduced constitutive models}

\setcounter{equation}{0}
\renewcommand\theequation{C\arabic{equation}}

Three simple reduced models are utilized in this paper, $D^0$, $D^1$, and $D^2$. $D^0$ model is analogous to the flat plate model for which:
\begin{equation}
\label{C1}
\begin{gathered}
	\frac{1}{1-\zeta T + \zeta^2 b} \approx 1 \Rightarrow \ii{I}{}{0} = h, \quad \ii{I}{}{1} = 0, \ii{I}{}{2} = \frac{h^3}{12}, \\
	\ii{D}{\lambda\gamma\chi\omega}{M} = h \ii{D}{\lambda\gamma\chi\omega}{}, \quad \ii{D}{\lambda\gamma\chi\omega}{MB} = 0, \quad \ii{D}{\lambda\gamma\chi\omega}{B} = \frac{h^3}{12} \ii{D}{\lambda\gamma\chi\omega}{}.	
\end{gathered}
\end{equation}
\noindent The other two models are based on the first-order Taylor approximation:
\begin{equation}
\label{C2}
	\frac{1}{1-\zeta T + \zeta^2 b} \approx 1 + \zeta T \Rightarrow \ii{I}{}{0} = h, \quad \ii{I}{}{1} = \frac{h^3}{12} T = T \ii{I}{}{2}, \quad  \ii{I}{}{2} = \frac{h^3}{12},
\end{equation}
\noindent where the $\ii{D}{}{M}$ and $\ii{D}{}{B}$ terms are the same as those of the $\ii{D}{0}{}$ model, but the coupling term $\ii{D}{}{MB}$ differs:
\begin{equation}
\label{C3}
\begin{aligned}
	\textrm{for} \; D^1 \; \textrm{model:} \quad - \ii{D}{\lambda\gamma\chi\omega}{MB} &= -\frac{h^3}{12} T \ii{D}{\lambda\gamma\chi\omega}{}, \\
	\textrm{for} \; D^2 \; \textrm{model:} \quad - \ii{D}{\lambda\gamma\chi\omega}{MB} &= -\frac{h^3}{12} T \ii{D}{\lambda\gamma\chi\omega}{} \\
	&\;\;\;\; + \frac{h^3}{24} \left( \ii{b}{\gamma}{\nu} \ii{D}{\lambda\nu\chi\omega}{} + \ii{b}{\lambda}{\nu} \ii{D}{\nu\gamma\chi\omega}{} + 2 \ii{b}{\omega}{\varphi} \ii{D}{\lambda\gamma\chi\varphi}{} + 2 \ii{b}{\chi}{\varphi} \ii{D}{\lambda\gamma\varphi\omega}{} \right).
\end{aligned}
\end{equation}

	\bibliography{shellbib} 
	\bibliographystyle{ieeetr}

\end{document}